\documentstyle[12pt]{article}
\input amssym.def
\input amssym.tex
\newcommand{\nc}{\newcommand}

\newcommand{\Hom}{\,{\rm Hom}\,}

\newcommand{\bla}{\phantom{bbbbb}}
\newcommand{\onebl}{\phantom{a} }
\newcommand{\eqdef}{\;\: {\stackrel{ {\rm def} }{=} } \;\:}

\newcommand{\half}{ {\frac{1}{2} } }
\newcommand{\vol}{ \,{\rm vol}\, }

\newcommand{\beq}{\begin{equation}}
\newcommand{\eeq}{\end{equation}}
\newcommand{\beqst}{\begin{equation*}}
\newcommand{\eeqst}{\end{equation*}}
\newcommand{\barr}{\begin{array}}
\newcommand{\earr}{\end{array}}
\newcommand{\beqar}{\begin{eqnarray}}
\newcommand{\eeqar}{\end{eqnarray}}
\newtheorem{theorem}{Theorem}[section]

\newtheorem{corollary}[theorem]{Corollary}

\newtheorem{lemma}[theorem]{Lemma}
\newtheorem{prop}[theorem]{Proposition}
\newtheorem{definition}[theorem]{Definition}
\newtheorem{remit}[theorem]{Remark}

\newcommand{\matr}[4]{\left \lbrack \begin{array}{cc} #1 & #2 \\
     #3 & #4 \end{array} \right \rbrack}

\newenvironment{rem}{\begin{remit}\rm}{\end{remit}}

\newcommand{\RR}{{\Bbb R }}
\newcommand{\CC}{{\Bbb C }}
\nc{\FF}{ {\Bbb F} } 
\newcommand{\ZZ}{{\Bbb Z }}

\newcommand{\UU}{{\Bbb U }}

\newcommand{\cala}{{\mbox{$\cal A$}}}

\newcommand{\cald}{{\mbox{$\cal D$}}}

\newcommand{\calf}{{\mbox{$\cal F$}}}
\newcommand{\calg}{{\mbox{$\cal G$}}}
\newcommand{\calh}{{\mbox{$\cal H$}}}

\newcommand{\call}{{\mbox{$\cal L$}}}
\newcommand{\calm}{{\mbox{$\cal M$}}}

\def\a{\alpha}
\def\b{\beta}
\def\g{\gamma}
\def\d{\delta}
\def\e{\epsilon}
\def\z{\zeta}

\def\l{\lambda}

\def\r{\rho}

\def\G{\Gamma}

\def\L{\Lambda}

\def\S{\Sigma}

\def\Ph{\Phi}

\def\O{\Omega}

\newcommand{\renorm}{{ \setcounter{equation}{0} }}

\nc{\Imm}{ {\rm Im} }
\nc{\Si}{\Sigma}
\nc{\si}{\sigma}

\nc{\liek}{  \mbox{\bf k}  }
\nc{\liet}{  {\mbox{\bf t}  } } 
\nc{\lieg}{  {\mbox{\bf g}  } } 
\nc{\lieb}{  {\mbox{\bf b}  } } 

\nc{\lietp}{{ {\liet}^\perp  }}
\nc{\liets}{ {\liet}^* }
\nc{\lieks}{ {\liek}^* }
\nc{\hk}{H^*_K}
\nc{\hht}{H^*_T}
\nc{\lietpl}{{ \liet_+} }

\nc {\Om}{\Omega}
\nc{\om}{\omega}

\nc{\diag}{ {\rm diag} }

\nc{\lrar}{\longrightarrow}
\nc{\Proof}{ \noindent{\em Proof:} }
\nc{\Cok}{ {\rm Cok} }

\nc{\inpr}[1]{ \langle #1 \rangle }
\nc{\zloc}{{ \mu^{-1}(0) } }
\nc{\fk}{F_K}
\nc{\ft}{F_T}

\nc{\eva}[2]{ { #1(#2)   } }
\nc{\evab}[2]{ { #2(#1)   } }

\nc{\epin}[1]{e^{ \isq \inpr{ #1} }  }
\nc{\emin}[1]{e^{ -\isq \inpr{ #1} }  }
\nc{\epinev}[2]{e^{ \isq \eva{ #1}{#2  } }  }
\nc{\eminev}[2]{e^{ -\isq \eva{ #1}{#2  } }  }
\nc{\epinevb}[2]{e^{ \isq \evab{ #1}{#2  } }  }
\nc{\eminevb}[2]{e^{ -\isq \evab{ #1}{#2  } }  }

\nc{\intk}{\int_{\phi \in \liek} }
\nc{\intt}{\int_{\phi \in \liet} }
\nc{\inttperp}{\int_{\phi \in \lietp} }
\nc{\inttpl}{\int_{\phi \in \lietpl} }
\nc{\tintk}{\int_{z \in \lieks} }
\nc{\tintt}{\int_{z\in \liets} }
\nc{\bom}{ {\bar{\omega} } }

\nc{\fism}{{F_I}}
\nc{\fisml}{{F_{I,\l}}}
\nc{\filg}{{\hat{F_I}}}
\nc{\filgl}{{\hat{F_{I,\l}}}}
\nc{\filgint}{{\filg \cap (\mu_\torc^{-1}(0)/\torc)    }}
\nc{\filglint}{{\filgl \cap (\mu_\torc^{-1}(0)/\torc)    }}
\nc{\mmr}{{(r(n-r))}}

\nc{\cent}{c}

\nc{\xred}{{M_{\rm red} } }
\nc{\mred}{\xred}
\nc{\xc}{{ M(c)} } 
\nc{\emtc}{ {M_{\liet} (c) } } 
\nc{\tran}{{\L_0}}
\nc{\stran}{s_{\L_0}}
\nc{\stranm}{s_{ m}}
\nc{\bh}{{ \underline{h} } }
\nc{\mcex}{ { M(c \exp (\xi) ) } }
\nc{\mexp} {\mcex}
\nc{\mtc}{ {M_\liet(c)} }
\nc{\mt}{M_t}
\nc{\mtcexp}{{M_\liet(c \exp (\xi) )}}
\nc{\mcexp}{{M (c \exp (\xi) )}}

\nc{\mtinv} {\mu_T^{-1}}
\nc{\minv} {\mu^{-1}}
\nc{\tf}{\tilde{f}}
\nc{\cexp}{c \exp(\xi)}

\nc{\kpr}{\kappa'}
\nc{\xvec}{X}

\nc{\res}{{\rm Res}}
\nc{\mredt}[1]{M_{\rm red}^T(#1)}
\nc{\mtred}{M^T_{\rm red} }
\nc{\normcon}{{}}

\nc{\sx}{s_\xi^*}

\nc{\he}[1]{{\hat{e_{#1}} } }
\nc{\hw}[1]{{\hat{w_{#1}} } }
\nc{\lietone}{{\hat{\liet}_1}}
\nc{\tone}{ { \hat{T}_1} } 
\nc{\lietl}{{\liet_n}}
\nc{\lietc}{{\liet_{n-1}}}
\nc{\torc}{ {T_{n-1} } }
\nc{\subgp}{ S\Bigl (U(n-1) \times U(1) \Bigr ) }
\nc{\tor}{{T_{n}}}
\nc{\nl}{{N_n}}
\nc{\ncom}{{N_{n-1}}}
\nc{\tnl }{{ \tilde{N_n}} }
\nc{\phil}{\Phi_n}
\nc{\phione}{\hat{\Phi}_1}
\nc{\phicom}{\Phi_{n-1}}

\nc{\mtlcexp}{ M_{\lietl} (c \exp (\xi) ) }
\nc{\hnl}{ {\widehat{N_n}} }
\nc{\efilg}{{e_\filg}}
\nc{\efn}{{e_\filg}}

\nc{\dbar}{{\bar{\nusym} }}
\nc{\tnlrt}{{\filg} }
\nc{\ccom}{{c_{n-1}}}
\nc{\tccom}{{\tilde{\ccom}}}
\nc{\tc}{\tilde{c}}
\nc{\cl}{{c_n}}
\nc{\tcl}{{\tilde{\cl}}}

\nc{\trivfac}{|W|}

\nc{\mredtorc}{{ M_{red,n-1} }}
\nc{\fred}{{F}}

\nc{\muklm}{{\mu_{K_{n-1}}^{-1}}}

\nc{\xicom}{ {\xi_{n-1}} }
\nc{\ncomc}{{N_{n-1}(\ccom)}}
\nc{\bpsi}{ {\bar{\Psi} }}

\nc{\tF}{\tilde{F}}
\nc{\ee}[1]{\exp 2 \pi i {#1}}
\nc{\eex}[1]{e^{ 2 \pi i {#1}}}
\nc{\xia}{\xi^{(I')}}
\nc{\xib}{\xi^{(J')}}
\nc{\epsr}[1]{{\epsilon_{#1}}}
\nc{\ccr}{c_r}
 \nc{\ccnr}{c_{n-r}}
 \nc{\ccn}{c_n}
\nc{\tpsi}{\tilde{\Psi}}
\nc{\ha}[1]{h_{#1}^{[I_1]} }
\nc{\hb}[1]{h_{#1}^{[I_2]} }
\nc{\tccr}{{\tilde{c_r}}}
\nc{\tccnr}{{\tilde{c_{n-r}}}}
\nc{\tccn}{{\tilde{c_n}}}
\nc{\chicom}{{M'(c \exp \xi)}}
\nc{\compform}{{\alpha}}
\nc{\Chi}{\Upsilon}
\nc{\bbeta}{\beta}
\nc{\pphi}{\phi}

\nc{\mnd}{\calm(n,d)}

\nc{\weightl}{\Lambda^w}
\nc{\intlat}{{\Lambda^I} } 
\nc{\mto}{\calm(2,1)}
\nc{\calgo}{\calg_0}
\nc{\mgh}{\mu_{G/H}} 
\nc{\ct}{ {\tilde{c}} }
\nc{\fd}{{F_\d}}
\nc{\nilp}{\e}
\nc{\fwc}{{\rm FWC}}
\nc{\hfeta}{{h^\eta_F}}

\nc{\frob}{\small \mid\!\times}

\nc{\mc}{{M(c)}}         
\nc{\bmc}{\bar{\alpha}}  

\nc{\yyb}{Y_c}

\nc{\epc} { {e_c} } 
\nc{\homfk}{ {\rm Hom} (\FF, K ) } 
\nc{\ktg} {K^{2g} } 
\nc{\proj} {{\rm pr} } 
\nc{\hl}{ {(h, \L)} } 

\nc{\inds}{ \delta}
\nc{\tuple} {(d_1, \dots, d_n) } 

\nc{\isq}{  } 
\nc{\homft}{ {\rm Hom} (\FF, T ) } 
\nc{\nusym}{{  \cald} }

\nc{\ev}[1]{ { {\rm ev}_{#1}}}
\nc{\hop}{\calh}

\nc{\tfr} { {\tilde{f}_r(X)} }
\nc{\tbrj}{ { {\tilde{b}_r}^j(X)} }
\nc{\tar}{{\tilde{a}_r(X)} }
\nc{\tbr}{{\tilde{b}_r^j(X)}}
\nc{\tfro} { {\tilde{f}_r(X)_1} }
\nc{\tbrjo}{ { {\tilde{b}_r}^j(X)_1} }
\nc{\tfrt} { {\tilde{f}_r(X)_2} }

\nc{\phiok}{\Phi_1^K}

\nc{\lw}{\Lambda^w}

\nc{\Ad} { { \rm Ad} } 
\nc{\ad} { { \rm ad} }

\begin{document}

\title{Intersection theory on   moduli spaces of holomorphic bundles 
of arbitrary rank on a Riemann surface}
\author{Lisa C. Jeffrey \\
Mathematics Department \\
Princeton University \\  Princeton, NJ 08544, USA 
\thanks{This material is based on work
supported by the National Science Foundation under Grant No.
DMS-9306029, and by grants from NSERC and FCAR.} 
  \\and\\
Frances C. Kirwan \\ Balliol College \\ Oxford OX1 3BJ, UK}
\date{August 1996; revised November 1997\\ $\phantom{a}$  \\ alg-geom/9608029} 
\maketitle

\renorm 
\section{Introduction}

Let $n$ and $d$ be coprime positive integers, and define $\mnd$ 
to be the moduli space of (semi)stable holomorphic
vector bundles of rank $n$, degree $d$ and fixed determinant
on a compact Riemann surface $\Sigma$. This moduli space
is a compact
K\"ahler manifold which has been studied
from many different points of view for more than three decades
(see for instance
Narasimhan and Seshadri 1965 \cite{AB,NS}). The subject of this
article is the characterization of the intersection pairings in the cohomology
ring\footnote{Throughout this paper all cohomology groups
will have complex coefficients, unless specified otherwise.} $H^*(\mnd)$.
A set of generators of this ring was described by Atiyah and Bott in their
seminal 1982
paper \cite{AB} on the Yang-Mills equations
on Riemann surfaces (where in addition inductive formulas for the Betti numbers of
$\mnd$ obtained earlier using number-theoretic
methods \cite{DR,HN} were rederived).  By Poincar\'{e} duality,
knowledge of the intersection pairings between products of these
generators (or equivalently knowledge
of the evaluation on the fundamental class of products of 
the generators) 
completely determines the structure of the cohomology ring.

In 1991 Donaldson \cite{Do} and Thaddeus \cite{T} 
gave formulas for the 
intersection pairings between products of these generators
in $H^*(\calm(2,1)) $ (in terms of Bernoulli
numbers). Then using physical methods, Witten \cite{tdgr}
found formulas for generating functions from which could be extracted
the intersection pairings between products  of these  generators
in $H^*(\mnd)$ for general rank $n$. These generalized his (rigorously 
proved) formulas \cite{qym} for the symplectic volume of $\mnd$:
 for instance, the symplectic 
volume of $\calm(2,1)$ is given by 
 \beq \label{svol} {\rm vol}(\calm(2,1)) = 
 \Bigl  (1 - \frac{1}{2^{2g-3}} \Bigr ) \frac{ \zeta(2g-2)}{2^{g-2}
\pi^{2g-2}} = \frac{2^{g-1}-2^{2-g}}{(2g-2)!} |B_{2g-2}| \eeq
where
$g$ is the genus of the Riemann surface, $\zeta$ is the Riemann 
zeta function and $B_{2g-2}$ is a Bernoulli number (see \cite{qym,T,Do}). 
The purpose of this paper is to obtain a mathematically
rigorous proof  of Witten's formulas for general rank $n$. Our announcement \cite{JK2} sketched
the arguments we shall use, concentrating mainly on the case of rank $n=2$.

The proof involves an application of the {\em nonabelian localization
principle} \cite{JK1,tdgr}. Let $K$ be a compact connected Lie
group with Lie algebra $\liek$, let
$(M,\omega)$ be a compact symplectic manifold equipped
with a Hamiltonian action of  $K$ 
and suppose that $0$ is a regular value of the moment map $\mu:M \to \lieks$
for this action.  One can use equivariant cohomology
on $M$ to study the cohomology ring of the reduced space, or symplectic quotient,
$\mred = \zloc / K$, which is an orbifold with an induced symplectic
form $\omega_0$. In particular it is shown in \cite{Ki1} that there is a
natural surjective homomorphism 
from the equivariant cohomology $H^*_K(M)$ of $M$ to the cohomology
$H^*(\mred)$ of the reduced space. For any
cohomology class $\eta_0 \in H^*(\mred)$
coming from $\eta \in H^*_K(M)$ via this map, we derived in \cite{JK1} a
formula (the residue formula, Theorem 8.1 of \cite{JK1})
for the evaluation $\eta_0[\mred]$ of  $\eta_0$
on the fundamental class of $\mred$. This formula involves the data that
enter the Duistermaat-Heckman formula \cite{DH}, and its generalization the
abelian localization formula \cite{ABMM,BV1,BV2} for the action of a maximal torus
$T$ of $K$ on $M$: that is, the set $\calf$ of connected
components $F$ of the fixed point set $M^T$ of
the action of $T$ on $M$, and the equivariant Euler classes $e_F$ of their normal
bundles in $M$. Let $\liet$ be the Lie algebra
of $T$; then the composition $\mu_T:M \to \liets$ of $\mu:M \to \lieks$ with the natural
map from $\lieks$ to $\liets$ is a moment map for the action of $T$ on $M$.
In the case when $K=SU(2)$ and  the order of the stabilizer in $K$ of a generic
point of
$\mu^{-1}(0)$ is $n_0$, the residue formula can be expressed in the form
\beq \label{rf2} \eta_0e^{\omega_0}[\mred] =  \frac{n_0}{2}{\rm Res}_{X=0} \Bigl( (2X)^2 \sum_{F\in\calf_+}
h^{\eta}_F(X) dX\Bigr) \eeq
where the subset $\calf_+$ of $\calf$ consists of those components $F$ of the fixed point set $M^T$
on which the value taken by the $T$-moment map $\mu_T:M \to \liets \cong \RR$ is
positive, and for $F\in \calf_+$ the inclusion of $F$ in $M$ is denoted by $i_F$ and the 
meromorphic function $h^{\eta}_F$ of $X\in\CC$ is defined
by
$$h_F^{\eta}(X) = \int_F \frac{i_F^* \eta(X) e^{\overline{\omega}(X)}}{e_F(X)} =
e^{\mu_T(F)(X)} \int_F \frac{i_F^*\eta(X) e^{\omega}}{e_F(X)}.$$
when $X\in \CC$ has been identified with ${\rm diag}(2\pi i,-2\pi i)X \in \liet\otimes\CC$.
Here $\overline{\omega}$ is the extension $\omega + \mu$ of the symplectic form
$\omega$ on $M$ to an equivariantly closed 2-form, while as before $\omega_0$ 
denotes the induced
symplectic form on $\mred$. 
Finally ${\rm Res}_{X=0}$ denotes the ordinary residue at $X=0$.

The moduli space $\mnd$ was described by 
Atiyah and Bott \cite{AB} as the symplectic reduction
of an infinite dimensional symplectic affine space 
$\cala$ 
with respect to 
the action of an infinite dimensional group $\calg$ 
(the 
{\em gauge group}).\footnote{To obtain his generating 
functionals, Witten formally applied his version of nonabelian 
localization to the action of the gauge group on the infinite dimensional space $\cala$.}
However $\mnd$ can also be exhibited
as the symplectic quotient of a finite dimensional
symplectic space $M(c)$ by the Hamiltonian action of the finite
dimensional group $K=SU(n)$. 
One
characterization of the space $M(c)$ is that it is the symplectic reduction 
of the infinite dimensional affine space $\cala$ 
 by the action of the {\em based} gauge group $\calgo$ (which
is the kernel
of the evaluation map $\calg \to K$ at a prescribed basepoint: 
see
\cite{ext}).  Now if a compact group $G$ containing a 
closed normal subgroup
$H$  acts in a Hamiltonian fashion on a symplectic manifold $Y$, then 
one may \lq\lq reduce in stages'': the space $\mu_H^{-1}(0)/H $ 
has a residual Hamiltonian action of the quotient group $G/H$ with
moment map $\mgh: \mu_H^{-1}(0)/H  \to {(\lieg/{\bf h} )}^* $, and 
$\mu_G^{-1}(0)/G$ is naturally identified as a symplectic 
manifold with 
$\mgh^{-1}(0)/(G/H). $ Similarly   $M(c)$ has a 
Hamiltonian action of $\calg/\calg_0 \cong K$, and the symplectic reduction 
with respect to this action is identified with the 
symplectic reduction of $\cala$ with respect to the full gauge group 
$\calg$.

Our strategy for obtaining Witten's formulas is to apply
 nonabelian localization
to this extended moduli space $M(c)$, 
which has a much more concrete (and entirely
finite-dimensional) characterization described in
 Section 4 below. Unfortunately
technical difficulties arise, because $M(c)$ 
is both singular and noncompact. The
noncompactness of $M(c)$ causes the more serious problems, the most immediate
of which is that there are infinitely many components $F$ of the fixed point
set $M(c)^T$. These, however, are easy to identify (roughly
speaking they correspond to bundles which are direct sums of line
bundles), and there are obvious
candidates for the equivariant Euler classes of their normal bundles, if the
singularities of $M(c)$ are ignored. In the case when $n=2$, for example, a
na\"{\i}ve application of the residue formula (\ref{rf2}), 
with some sleight of hand, 
would yield
$${\rm vol}(\calm(2,1)) = e^{\omega_0} [\calm(2,1)]
  = (-1)^g {\rm Res}_{X=0} \sum_{j=0}^{\infty} 
\frac{e^{(2j+1)X}}{2^{g-2}X^{2g-2}}$$
\beq \label{f3} \bla = (-1)^g {\rm Res}_{X=0} 
\frac{e^{X}}{2^{g-2}X^{2g-2}(1-e^{2X})}
  = (-1)^{g-1}{\rm Res}_{X=0} \frac{1}{2^{g-1}X^{2g-2}{\rm sinh}(X)}.\eeq
This does give the correct answer (it agrees with (\ref{svol}) above). 
However it is far from obvious how this calculation
might be justified, since the infinite sum does not converge in
a neighbourhood of $0$, where the residue is taken, and indeed the sum of the
residues at $0$ of the individual terms in the sum does not converge.

These difficulties can be overcome by making use of a different approach
to nonabelian localization given recently by Guillemin-Kalkman 
\cite{GK}
and independently by Martin \cite{Ma}. 
This is made up of two steps: the first is to reduce to the case of a torus
action, and the second, when $K=T$ is a torus, is to study the change in
the evaluation on the fundamental class of the reduced space
$\mu_T^{-1}(\xi)/T$ of the cohomology class induced by $\eta$,
as $\xi$ varies in $\liets$. It is in fact an immediate
consequence of the residue formula that if $T$ is a maximal torus
of $K$ and $\xi\in\liets$ is any regular value
sufficiently close to $0$
of the $T$-moment map $\mu_T:M\to \liets$, 
then the evaluation $\eta_0[\mred]$ of  $\eta_0\in H^*(\mred)$
on the fundamental class of $\mred = \zloc/K$ is equal to the evaluation
of a related element of $H^*(\mu_T^{-1}(\xi)/T)$
on the fundamental class of the $T$-reduced
space $\mu_T^{-1}(\xi)/T$. This was first observed by Guillemin and Kalkman
\cite{GK} and by Martin \cite{Ma}, who gave an independent proof which
showed that $\eta_0[\mred]$ is also equal to an evaluation on
$$\zloc/T = (M_{\liet} \cap \mu_T^{-1}(0))/T$$
where $M_{\liet} = \mu^{-1}(\liet)$. In our situation the space $M_{\liet}$
turns out to be \lq\lq periodic'' in a way which enables us to avoid working with infinite
sums except in a very trivial sense. This is done by comparing the results
of relating evaluations on
$(M_{\liet} \cap \mu_T^{-1}(\xi))/T$ for different values 
of $\xi$ in two ways: using
the periodicity and using Guillemin and Kalkman's arguments, 
which can be made to
work in spite of the noncompactness of $M(c)$. The 
singularities can be dealt with
because $M(c)$ is embedded naturally and equivariantly
in a nonsingular space, and integrals over $M(c)$ can be rewritten
as integrals over this nonsingular space.

In the case when $n=2$ our approach gives expressions for the pairings in
$H^*(\mto)$ as residues similar to those in 
(\ref{f3}) above. When $n>2$ we consider the 
action of a suitable one-dimensional subgroup 
$\hat{T}_1$ of $T$, with Lie algebra
$\hat{\liet}_1$ say, on the quotient of $\mu^{-1}(\hat{\liet}_1)$ 
by a subgroup of $T$ whose
Lie algebra is a complementary subspace to $\hat{\liet}_1$ 
in $\liet$. This leads
to an inductive formula for the pairings on $H^*(\mnd)$, 
and thus to expressions
for these pairings as iterated residues (see Theorems \ref{mainab} and
\ref{t9.6}
below, which are the central  results of this paper). Witten's
formulas, on the other hand, express the pairings as infinite sums over those
elements of the weight lattice of $SU(n)$ which lie in the interior of a
fundamental Weyl chamber (see Section 2). These infinite sums are
difficult to calculate in general, and there is apparently (see
\cite{tdgr} Section 5) no
direct proof even that they are always zero when the pairings they represent
vanish on dimensional grounds.  However, thanks to an argument of Szenes
(see Proposition \ref{p:sz} below), 
Witten's formulas can be identified with the
iterated residues which appear in our approach.

Over the moduli space $\mnd$ there is a natural line bundle $\call$ (the
Quillen line bundle \cite{Q}) whose fibre at any point representing
a semistable holomorphic bundle $E$ is the determinant line
$${\rm det} \bar{\partial} = {\rm det} H^1(\S,E) \otimes {\rm det} H^0(\S,E)^*$$
of the associated $\bar{\partial}$-operator. Our expressions for pairings
in $H^*(\mnd)$ as iterated residues, together with the Riemann-Roch
formula, lead easily (cf. Section 4 of \cite{Sz}) to a proof of the Verlinde
formula for
$$\dim H^0(\mnd,\call^k)$$
for positive integers $k$ (proved by Beauville and Laszlo in
\cite{BL}, by Faltings in \cite{F}, 
by Kumar, Narasimhan and Ramanathan in \cite{KNR} and by
Tsuchiya, Ueno and Yamada in \cite{TUY}).

This paper is organized as follows. In Section 2 we describe the generators
for the cohomology ring $H^*(\mnd)$ and Witten's formulas for the intersection
pairings among products of these generators. In Section 3 we outline tools from 
the Cartan model of 
equivariant cohomology, which will be used in later sections, and
the different versions of localization which will be relevant.
In Section 4 we recall properties of 
 the extended moduli space $M(c)$, and in Section 5 we construct the
equivariant differential forms representing equivariant Poincar\'{e}
duals which enable us to rewrite integrals over singular spaces
as integrals over ambient nonsingular spaces. Then Section 6 begins 
the application of nonabelian
localization to the extended moduli space, and Section 7 
analyses the fixed point
sets which arise in this application. Section 8 uses 
induction to complete the proof of 
Witten's formulas when the pairings are between cohomology classes of a
particular form, Section 9 extends the inductive argument to give
formulas for all pairings, and in Section 10 it is shown that these
agree with Witten's formulas. Finally as an application 
Section 11 gives a proof of the Verlinde formula for $\mnd$.

We would like to thank the Isaac Newton Institute in Cambridge, the Institute
for Advanced Study in Princeton, the Institut Henri Poincar\'{e}
and Universit\'{e} Paris VII, the Green-Hurst Institute for Theoretical
Physics in Adelaide and the Massachusetts Institute of
Technology for their hospitality during crucial phases in the evolution
of this paper. We also thank A. Szenes for pointing out an 
error in an earlier version of the paper: since the original 
version 
of this paper was written, Szenes has obtained new results \cite{Sz2}
which are closely related to the results given in 
Section 11 of our
paper.

\renorm 
\section{The cohomology of the moduli space $\mnd$ and Witten's formulas for
intersection pairings}

\nc{\lb}[1]{ {l_{#1} } }
\nc{\yy}[1]{Y_{#1}}
\nc{\liner}[1]{L_{#1} }
\nc{\lambdr}[1]{\Lambda_{#1} }
\nc{\linestd}{\liner{(\lb{1}, \dots, \lb{n-2} )} }
\nc{\lambstd}{\lambdr{(\lb{1}, \dots, \lb{n-2} )} }
\nc{\expsum}[1]{ { (e^{-{#1} } - 1 ) } }
\nc{\itwopi}{ { 2 \pi i }}
\nc{\Res}{\res}
\nc{\indset}{{(l_1, \dots, l_{n-1})}  }
\nc{\indsettwo}{{(l_1, \dots, \dots, l_{n-2})}  }

\nc{\laregint}{ {\lambstd^{reg} } }

In order to avoid exceptional cases, we shall assume 
throughout that the Riemann surface $\Sigma$ has genus $g\geq 2$.

A set of generators for
the cohomology\footnote{In this 
paper, all cohomology groups are assumed to be with complex
coefficients.}
  $H^*(\mnd)$
of the moduli space $\mnd$ of stable holomorphic
vector bundles of coprime rank $n$ and degree 
$d$ and fixed determinant on a compact 
Riemann surface $\Sigma$ of genus $g \ge 2 $ 
 is given in \cite{AB} by Atiyah and Bott.
It may be described as follows.
There is a universal rank $n$ vector  bundle
$$ \UU \to \Sigma \times \mnd $$
which is unique up to tensor product with the pullback of any holomorphic
line bundle on $\mnd$; for definiteness Atiyah and Bott
impose an extra normalizing condition which determines
the universal bundle up to isomorphism, but this is
not crucial to their argument (see \cite{AB}, p. 582). 
 Then by \cite{AB} Proposition 2.20 the following elements
of  $H^*(\mnd)$ for $2\leq r\leq n$ make up a set of generators: 
$$ f_r = ([\Sigma], c_r(\UU)), $$
$$ b_r^j = (\a_j, c_r(\UU)), $$
$$ a_r = (1, c_r(\UU)). $$
Here, $[\Sigma]$
 $ \in H_2(\Sigma)$ and $\a_j \in H_1(\Sigma)$ $(j = 1, \dots, 2g)$
form standard bases  of $H_2(\Sigma, \ZZ)$ and  $H_1(\Sigma, \ZZ)$, and
the bracket represents the slant product $H^N(\Sigma \times \mnd)
\otimes H_j(\Sigma) \to H^{N-j}(\mnd). $ 
 More generally if $K=SU(n)$ and $Q$ is an invariant polynomial of
degree $s$ on its Lie algebra $\liek = su(n)$ then there is an
associated element of $H^*(BSU(n))$ and hence an associated
element of $H^*(\Sigma \times \mnd)$ which is a characteristic
class $Q(\UU)$ of the universal bundle $\UU$. Hence the slant
product gives rise to classes 
$$ ([\Sigma],Q(\UU)) \in H^{2s-2}(\mnd), $$
$$ (\alpha_j,Q(\UU)) \in H^{2s-1}(\mnd), $$
and
$$ (1,Q(\UU)) \in H^{2s}(\mnd). $$
In particular, letting $\tau_r$ $\in S^r(\lieks)^K$ 
denote the invariant polynomial associated to the 
$r$-th Chern class, we recover 
\beq \label{1.2} \label{9}
f_r = ([\Sigma],\tau_r(\UU)), \eeq  
$$ b_r^j = (\a_j,\tau_r (\UU)) $$
and
$$ a_r =(1, \tau_r(\UU)). $$
A special role is played by the invariant polynomial 
$\tau_2 = -  \inpr{\cdot , \cdot}/2 $  on $\liek$
given by the Killing form or invariant 
inner product.  We normalize the inner product
as follows for $K = SU(n)$: 
\beq \label{1.02} \inpr{X, X} = - {\rm Trace} (X^2)/(4\pi^2). \eeq
The  class $f_2$ associated to 
$-\inpr{ \cdot, \cdot}/2 $ is the cohomology class
of the symplectic form on $\mnd$.

As was noted in the introduction, Atiyah and Bott identify $\mnd$
with the symplectic reduction of an infinite dimensional affine
space ${\cal A}$ of connections by the action of an infinite
dimensional Lie group ${\cal G}$ (the gauge group). They show that
associated to this identification there is a natural surjective
homomorphism of rings from the equivariant cohomology ring
$H^*_{\bar{{\cal G}}}({\cal A})$ to $H^*(\mnd)$, where
$\bar{{\cal G}}$ is the quotient of ${\cal G}$ by its central
subgroup $S^1$. There is a canonical ${\cal G}$-equivariant
universal bundle over $\Sigma \times {\cal A}$, and 
the slant products of its
Chern classes with $1\in H_0(\Sigma)$, $\a_j \in H_1(\S)$ for
$1\leq j\leq 2g$ and $[\S]\in H_2(\S)$ give generators of
$H^*_{{\cal G}}({\cal A})$ which by abuse of notation
we shall also call $a_r$, $b_r^j$ and $f_r$. (In fact
$H^*_{{\cal G}}({\cal A})$ is freely generated by 
$a_1,\dots,a_n$, $f_2,\dots,f_n$ and $b_r^j$ for 
$1<r\leq n$ and $1\leq j\leq 2g$, subject only to the
usual commutation relations.) The surjection from ${\cal G}$
to $\bar{{\cal G}}$ induces an inclusion from
$H^*_{\bar{{\cal G}}}({\cal A})$ to $H^*_{{\cal G}}({\cal A})$ 
such that
$$H^*_{{\cal G}}({\cal A}) \cong H^*_{\bar{{\cal G}}}({\cal A})
\otimes H^*(BS^1)$$
if we identify $H^*(BS^1)$ with the polynomial subalgebra
of $H^*_{{\cal G}}({\cal A})$ generated by $a_1$, and then the
generators $a_r$, $f_r$ and $b_r^j$ for 
$1<r\leq n$ determine generators of $H^*_{\bar{{\cal G}}}({\cal A})$
and thus of $H^*(\mnd)$. These are the generators we shall use in this
paper. The normalization condition imposed by Atiyah and Bott corresponds
to using the isomorphism
$$H^*_{{\cal G}}({\cal A}) \cong H^*_{\bar{{\cal G}}}({\cal A})
\otimes H^*(BS^1)$$
obtained by identifying $H^*(BS^1)$ with the polynomial subalgebra
of $H^*_{{\cal G}}({\cal A})$ generated by $2(g-1)a_1 + f_2$; they choose this
condition because it has a nice geometrical interpretation in
terms of a universal bundle over $\S \times \mnd$.

In 
 Sections 4 and 5
of \cite{tdgr}, Witten obtained formulas for generating functionals 
from which one may extract all intersection pairings
$$ \prod_{r = 2}^n a_r^{m_r} f_r^{n_r} \prod_{k_r = 1}^{2g} 
(b_r^{k_r})^{p_{r, k_r} } [\mnd].$$
Let us begin with pairings of the form
\beq \label{1.00001}
\prod_{r = 2}^n a_r^{m_r} \exp f_2 [\mnd]. \eeq
When $m_r$ is sufficiently small  to ensure convergence
of the sum,  Witten obtains\footnote{In fact 
$\mnd$ is an $n^{2g}$-fold cover of  the space for which
Witten computes pairings: this accounts for the factor $n^{2g}$ in 
our formula  (\ref{2.005}). Taking this into account, 
(\ref{1.1}) follows from a special case of Witten's formula
\cite{tdgr} (5.21).}
\beq \label{1.1}
\prod_{r = 2}^n a_r^{m_r} \exp f_2 [\mnd]
= 
 c^\rho \G (-1)^{n_+ (g-1)}   \Biggl ( 
\sum_{\l \in \weightl_{\rm reg} \cap \liet_+  } \frac{c^{- \l} 
\prod_{r = 2}^n \tau_r(\itwopi \l) ^{m_r} }
{\nusym^{2g-2}(\itwopi \l) } \Biggr ),  \eeq
where
\beq \label{2.005}
\G = \frac{n^{2g} }{\# \Pi_1(K')} (\frac{\vol(K')}{(2 \pi)^{\dim
    K'}})^{2g-2} \Bigl ((2 \pi)^{n_+} \cald(\rho)\Bigr ) ^{2g-2}
= n^g 
\eeq
 is a universal constant for $K=SU(n)$ and 
$K' = K/Z(K)$, and 
the Weyl odd polynomial 
$\nusym$ on $\liets$ is defined by 
$$\nusym(X) = \prod_{\g > 0 } \g(X) $$
where $\g$ runs over the positive roots.
Here, $\rho$ is half the sum of the positive roots, 
and $n_+$ $ = n(n-1)/2$ is the number of positive roots.
The sum over $\l$ in
(\ref{1.1}) 
runs over those elements of the weight lattice $ \weightl$ 
that are in the interior of the fundamental Weyl chamber.\footnote{The weight
lattice $\weightl \subset \liet^*$ is the dual lattice of the integer
lattice $\intlat = {\rm Ker} (\exp) $ in $\liet$.}
The 
element 
\beq \label{1.p1} c = e^{2 \pi   i   d/n}
{\rm diag}( 1, \dots, 1) \eeq
is a generator of the 
centre $Z(K)$ of $K$, so since $\l \in  \liets$ is in 
${\rm Hom} (T, U(1))$, we may evaluate $\l$ on $c$ as in (\ref{1.1}): 
$c^\l$ is defined as $\exp \l (\tilde{c}) $ where $\tilde{c}$ is 
any element of the Lie algebra of $T$ such that 
$\exp \tilde{c} = c$.  Note that in fact when $d$ is coprime
to $n$ (so that when $n$ is even $d$ is odd) we have
$c^{\rho} = (-1)^{n-1}. $

Witten's formula \cite{tdgr} (5.21) covers pairings 
involving the $f_r$ for $r > 2$ and the 
$b_r^j$ as well as $f_2$ and the $a_r$. He obtains it by
reducing to the special case of pairings of the form (\ref{1.00001})
above (see \cite{tdgr}  Section 5, in particular the 
calculations (5.11) - (5.20)) and then applying
\cite{tdgr} (4.74) to this special case. In this special case of
pairings
of the form (\ref{1.00001}), Witten's formula
\cite{tdgr} (5.21) follows from our Theorem \ref{mainab}
using Proposition \ref{p:sz} below. Moreover our formula 
(Theorem \ref{t9.6}) for pairings involving all the 
generators $a_r$, $b_r^j$ and $f_r$ reduces to the special 
case just as Witten's does (see Propositions 
\ref{p9.1} and \ref{p9.2}). Thus 
Witten's formulas are equivalent to ours, although they look
very different (being expressed in terms of infinite 
sums indexed by dominant weights instead of in terms of iterated
residues).

For the sake of concreteness it is worth examining the special case when the
rank $n = 2$ so that
 the degree $d$ is odd. In fact, since tensoring by a fixed line
bundle of degree $e$ induces a homeomorphism between $\mnd$ and
$M(n,d+ne)$, we may assume that $d= 1$. 
In this case the dominant weights $\l$ are just the positive integers.
The relevant generators of $H^*(\mto)$ are \beq f_2 \in H^2(\mto)\eeq
(which is the cohomology class of the symplectic form 
on $\mto$)
and 
\beq a_2 \in H^4(\mto): \eeq
 these arise from the invariant polynomial 
$\tau_2 = -\inpr{\cdot,\cdot}/2$ by 
$a_2 = \tau_2(1)$, $f_2 = \tau_2([\Sigma]) $ 
(see (\ref{1.2})). We find then that the formula (\ref{1.1}) reduces for 
$m \le g - 2 $ to\footnote{Here, we have identified $-a_2$ with Witten's class
$ \Theta$ and $f_2 $ with Witten's class $\omega$.}
(\cite{tdgr}, (4.44))
\beq \label{1.3}
a_2^j \exp (f_2 )[{\calm}(2,1)] 
= \frac{ 2^{2g}}{2 (8 \pi^2)^{g-1}}
 \Biggl ( \sum_{n = 1}^\infty 
(-1)^{n + 1}\frac{  \pi^{2j}  }{ n^{2g-2-2j }}
\Biggr ). \eeq 
 Thus one obtains the 
  formulas   found by
 Thaddeus in Section 5 of \cite{T} for the intersection pairings
$a_2^m f_2^n [\mto]$; these intersection pairings are given by 
Bernoulli numbers, or equivalently are given in terms of the 
Riemann zeta function
 $\zeta(s) = \sum_{n \ge 1} 1/n^s$. As Thaddeus shows in Section 4
of \cite{T}, this
 is enough to determine all the intersection pairings in the case
when the rank $n$ is two, because all the pairings 
$$ a_2^{m} f_2^{n} \prod_{k = 1}^{2g} 
(b_2^{k})^{p_{ k} } [\calm(2,1)]$$
are zero except those of the form
$$ a_2^{m} f_2^{n} 
b_2^{2i_1-1}b_2^{2i_1 } \ldots b_2^{2i_q-1} b_2^{2i_q } [\calm(2,1)]$$
where $m+2n+3q=3g-3$ and 
$1\leq i_1 < \ldots <i_q \leq g$, and this expression equals
the evaluation of
 $a_2^m f_2^n$ on the corresponding moduli space of rank 2 and degree 1
bundles over a Riemann
surface of genus $g-q$ if $q\leq g-2$, and equals 4 if $q=g-1$.


\newcommand{\bracearg}[1]{ { [[ #1 ]] } }
\newcommand{\tildarg}[1]{ { [[ #1 ]]  } }

 Szenes \cite{Sz} has proved  that the expression on the right hand
side of 
(\ref{1.1}) 
may be rewritten in a particular form. To state the result we must
introduce some notation.
The Lie algebra $\liet = \liet_n$ of the maximal torus $T$ of $SU(n)$ 
is $$ \liet = \{ (X_1, \dots, X_n) \in \RR^n: X_1 + \dots 
+ X_n = 0 \}. $$
Define coordinates
$\yy{j} = e_j(X) = X_j - X_{j +1}$ on $\liet$ for $j = 1, \dots, n-1$. 
The positive 
roots of $SU(n)$  are then $\g_{jk}(X) = X_j - X_k$ 
$ = Y_j + \dots + Y_{k-1}$ for $1 \le j < k \le n$.
The {\em integer lattice} $\intlat$ of $SU(n)$ is generated by the simple
roots $e_j, j = 1, \dots, n - 1$. The dual lattice  to $\intlat$ 
with respect to the inner  product  $\inpr{ \cdot, \cdot} $ introduced 
at (\ref{1.02}) is the 
{\em weight lattice}  $\weightl$ $\subset \liet$: in terms of the 
inner product $\inpr{ \cdot, \cdot} $, it is given by 
$\weightl = \{ X \in \liet:   \; Y_j \in \ZZ  ~\mbox{for 
$j = 1, \dots, n - 1 $} \}$. We define also 
$\weightl_{\rm reg} (\liet_n)  = \{ X \in  \weightl: 
\; Y_j  \ne 0  ~\mbox{ for $ j = 1, \dots, n-1$}  $
and $\gamma_{jk} (X) \ne 0 $ for any $j \ne k \}. $

\begin{definition} \label{bracedef}
Let $f: \liet \otimes \CC \to \CC $ be a meromorphic 
function of the  form 
\beq \label{fundef} 
f(X) = g(X) e^{-\g (X)} \eeq 
where
$\g (X) = \gamma_1 Y_1 +
\dots + \g_{n-1} Y_{n-1}$ for $(\gamma_1, \dots, 
\gamma_{n-1} )$ $\in \RR^{n-1}$. We define
$$ \tildarg{\gamma} = (\tildarg{\gamma}_1, 
\dots, \tildarg{\gamma}_{n-1} )$$ to be the 
element of $\RR^{n-1}$ for which
$0 \le \tildarg{\gamma}_j < 1$ for all $j = 1, \dots, n-1$
and $\tildarg{\gamma} = \gamma ~{\rm mod}~ \ZZ^{n-1}$. 
(In other words, $\bracearg{\gamma} =$ $\sum_{j = 1}^{n-1} 
\tildarg{\gamma}_j e_j$ is the unique element of $\liet$ $\cong 
\RR^{n-1}$ which is in the fundamental domain defined by the simple
roots for the translation action on 
$\liet_n$ of the integer lattice, and which is 
 equivalent to $\gamma$ under translation by the integer lattice.)

We also define the meromorphic function $\bracearg{f}: 
 \liet \otimes \CC \to \CC $ by 
$$ \bracearg{f}(X) = g(X) e^{- \tildarg{\g}(X)}. $$
\end{definition}

\begin{prop} \label{p:sz}{\bf [Szenes]} 
Let $f: \liet\otimes \CC \to \CC$ be defined by 
$$
f(X) = \frac{\prod_{r = 2}^n \tau_r(X)^{m_r} e^{ - \ct(X) } }
{\nusym(X)^{2g-2} }. $$
Provided that
the $m_r$ are sufficiently small to ensure 
convergence of the sum, we have
$$ \sum_{\l \in \weightl_{reg} (\liet_{n}) \cap \liet_+ } 
f(\itwopi\l)  = \frac{1}{n!}
\Res_{Y_1 = 0} \dots 
\Res_{Y_{n-1} = 0 } \Bigl ( 
\frac{\sum_{w \in W_{n-1}}
\bracearg{w(f)}(X)}{\expsum{\yy{n-1} } \dots \expsum{\yy{1}} } \Bigr ), $$ 
where $W_{n-1} \cong S_{n-1}$ is the Weyl group of $SU(n-1)$
embedded in $SU(n)$ in the standard way using the first $n-1$
coordinates $X_1,\dots,X_{n-1}$.
\end{prop}

\begin{rem} \label{r2.1}
Here, we  have  introduced coordinates $Y_j = e_j(X)$ 
on $\liet$ using the simple roots 
$$\{ e_j: j = 1, \dots
n-1 \} $$
 of $\liet$, and  $ \weightl_{\rm reg}$ denotes the regular part of
 the 
weight lattice $\weightl$ (see below). 
Also, we have introduced the   unique   element $\tc$ 
of  $\liet$  which satisfies $e^{2 \pi i \tc} = c$ and which
 belongs to the fundamental domain defined
by the simple roots for the translation action on
$\tor$ of the integer lattice $\L^I$: this simply means that
$\inpr{\tc,X} = \g_1 Y_1 + \dots + \g_{n-1} Y_{n-1}$ where
$0\leq \g_j <1$ for $1\leq j\leq n-1$. (In the
notation introduced
in Definition \ref{bracedef} this  says that
$ \tc = \bracearg{ ~(d/n, d/n, \dots, -(n-1)d/n) }$.)
Also, $\lietpl$ denotes the fundamental Weyl chamber, which is a
fundamental domain for the action of the Weyl group on $\liet$.

If $g(Y_k, \dots, Y_{n-1})$ is a meromorphic
function of $Y_k, \dots, Y_{n-1}$, we interpret 
$\res_{Y_k = 0 } g(Y_{k}, \dots, Y_{n-1})$ as the ordinary
one-variable residue of $g$ regarded 
as a function of $Y_k$ with $Y_{k+1}, \dots, Y_{n-1}$ held constant.
\end{rem}

The rest of this section will be devoted to a proof of Proposition
\ref{p:sz}.

We shall prove the following theorem:
\begin{theorem} \label{t2.04}
Let $f: \liet_{n} \otimes \CC \to \CC$ be a meromorphic function 
of the form $f(X) = g(X) e^{-\g (X)}$ where 
$\g (X) = \gamma_1 Y_1 +
\dots + \g_{n-1} Y_{n-1}$ with $0 \leq \g_{n-1} < 1$, 
and $g(X)$ is a rational function of $X$
with poles only on the zeros of the roots $\g_{jk}$ and 
decaying sufficiently fast at infinity.
Then
$$\sum_{\l \in \weightl_{reg} (\liet_{n})} 
f(\itwopi\l)=  \Res_{Y_1 = 0} \dots 
\Res_{Y_{n-1} = 0 } \Bigl ( 
\frac{\sum_{w \in W_{n-1}} \bracearg{w(f)}(X) }{\expsum{\yy{n-1} } 
\dots \expsum{\yy{1}} } \Bigr ) $$
where $W_{n-1}$ is the Weyl group of $SU(n-1)$ embedded in $SU(n)$ 
using the first $n-1$ coordinates.
\end{theorem}

\begin{rem} \label{r:intadd}
Notice that if $f$ is as in the hypothesis of the Theorem 
(but here one may omit the hypothesis that $0 \le \g_{n-1} <1$)
 then
$$\sum_{\l \in \weightl_{reg} (\liet_{n})} 
f(\itwopi\l) = 
\sum_{\l \in \weightl_{reg} (\liet_{n})} 
\bracearg{f}(\itwopi\l), $$ where 
$\bracearg{f}$ was defined in Definition 
\ref{bracedef}.
\end{rem}

\noindent{\em Proof of Proposition \ref{p:sz} given Theorem 
\ref{t2.04}:} 
The function 
\beq \label{2.fn}
f(X) = \frac{\prod_{r = 2}^n \tau_r(X)^{m_r} e^{ - \ct(X) } }
{\nusym(X)^{2g-2} } \eeq
satisfies the hypotheses of the Theorem, provided that the 
$m_r$ are small enough to ensure convergence of the sum. Notice 
that if $\lambda \in \weightl_{\rm reg} (\liet_n)$ then 
$e^{ - \itwopi \ct(\l)} = c^{ - \l} $ satisfies 
$c^{- \l} = c^{- w\l } $ for all elements $w$ 
of the Weyl group $W$. Thus for this particular $f$ we have that
$$\sum_{\l \in \weightl_{reg} (\liet_{n})} 
f(\itwopi\l)  = n! 
\sum_{\l \in \weightl_{reg} (\liet_{n}) \cap \liet_+ } 
f(\itwopi\l) . $$ So
$$
\sum_{\l \in \weightl_{reg} (\liet_{n}) \cap \liet_+ } 
f(\itwopi\l)  = \frac{1}{n!}
\Res_{Y_1 = 0} \dots 
\Res_{Y_{n-1} = 0 } \Bigl ( 
\frac{\sum_{w \in W_{n-1}} \bracearg{w(f)}(X)}{\expsum{\yy{n-1}
 } \dots \expsum{\yy{1}} } \Bigr ) $$
which is the statement of Proposition \ref{p:sz}.\hfill $\square$

It remains to prove Theorem \ref{t2.04}. By induction on 
$n$ it suffices to prove 
\begin{lemma} \label{l2.04} 
Let  $f  = f_{(n)} : \liet_n \to \CC$ be as in the statement of
Theorem \ref{t2.04}. 
Define 
$f_{(n-1)} : \liet_{n-1} \to \CC$  by 
$$f_{(n-1)}
(Y_1, \dots, Y_{n-2}) 
= \Res_{\yy{n-1} = 0} \frac{f(Y_1, 
\dots, Y_{n-1})}{e^{- \yy{n-1}} - 1}. $$
Then 
$$ \sum_{\l \in \weightl_{reg} (\liet_{n} )} f_{(n)} (\itwopi \l) 
= \sum_{\l \in \weightl_{reg} (\liet_{n-1}) } 
\sum_{j=1}^{n-1} (q_j f)_{(n-1)} (\itwopi \l), $$
where $q_j$ is the element of the Weyl group $W_{n-1} \cong S_{n-1}$
represented by swapping the coordinates $X_j$ and $X_{n-1}$.
\end{lemma}

\begin{rem} 
Note that by Remark \ref{r:intadd}, the sum
$\sum_{j=1}^{n-1} (q_j f)_{(n-1)} (\itwopi \l)$ is equal to 
$$\sum_{j=1}^{n-1} \bracearg{ (q_j f)_{(n-1)} } (\itwopi \l). $$
\end{rem}
\begin{rem} Note that 
the function $ \bracearg{ (q_j f)_{(n-1)} }$ satisfies the
hypotheses of Theorem \ref{t2.04}.
\end{rem}

\noindent{\em Proof of Lemma \ref{l2.04}:}
 Let 
$\lb{j} $ $(j  = 1, \dots, n - 2) $ be integers such that 
\beq \label{2.cond} \lb{j} + \lb{j+1}  + \dots + \lb{k} \ne 0
~\mbox{for any}~ 1 \le j \le k \le n - 2. \eeq
Define $\liner{(\lb{1}, \dots, \lb{n-2} )}$ 
to be the line $\{ (\itwopi \lb{1}, \dots, \itwopi \lb{n-2}, \yy{n-1} 
 ): \yy{n-1} \in \CC \}$. 
The condition (\ref{2.cond}) states that all the 
roots $\g_{jk}$ for $1 \le j < k \le n - 1 $ are nonzero on 
$\liner{(\lb{1}, \dots, \lb{n-2} )}$.

Let $f: \liet\otimes \CC \to \CC$ be a meromorphic function as
in the statement of Theorem \ref{t2.04}, having
poles only at the zeros of the roots 
$\g_{jk}$. 
We shall think of $f$ as a function $f(Y_1, \dots, Y_{n-1})$ of the 
coordinates $Y_1, \dots, Y_{n-1}. $ 
Define $\laregint$ to be 
$$\laregint = \{ X \in \linestd: ~Y_j(X) = X_j - X_{j+1} \in 
\itwopi \ZZ, $$
$$~\g_{jk}(X) \ne 0 ~\mbox{for any
$j \ne k $} \}. $$
The sum of all residues of the function 
$g_{\indsettwo} $ on $\CC$ given by $$g_\indsettwo (\yy{n-1}) = 
\frac{f(\itwopi \lb{1}, \itwopi \lb{2}, \dots, 
 \yy{n-1} ) }{ e^{- \yy{n-1} } - 1 } $$
is zero  and these residues occur when $\yy{n-1} \in \itwopi \ZZ$.
Therefore
we find that 
the sum 
$$ \sum_{ p  \in \laregint }  f (p) $$ is given by 
$$
- \sum_{(\itwopi\lb{1}, \dots,\itwopi \lb{n-2},
\itwopi \lb{n-1} ) \in \laregint }
 \Res_{\yy{n-1} = \itwopi l_{n-1} }
\frac{f(\itwopi \lb{1}, \itwopi \lb{2}, \dots, \itwopi \lb{n-2},
 \yy{n-1} ) }{ e^{- \yy{n-1} } - 1 } $$
\beq \label{2.001} = \sum_{j =1}^{n-2} 
\Res_{\yy{n-1} = -\itwopi (\lb{j}  + \dots + \lb{n-2} )  }
\frac{f(\itwopi \lb{1}, \dots, \itwopi \lb{n-2},
 \yy{n-1} ) }{ e^{- \yy{n-1} } - 1 } \: + \eeq
$$ \Res_{\yy{n-1} = 0  }
\frac{f(\itwopi \lb{1}, \dots, \itwopi \lb{n-2},
 \yy{n-1} ) }{ e^{- \yy{n-1} } - 1 }. $$

\nc{\lj}[1]{l^{(j)}_{#1} }

\begin{prop} \label{p2.3} Let $p_j$ be the point
$ X \in \linestd $ for which $Y_{n-1} = - \itwopi (l_j + 
\dots + l_{n-2})$, or equivalently
$ X_n = X_j . $ 
Then 
$$ \res_{Y_{n-1} = - \itwopi (l_j + \dots + l_{n-2} ) } \Biggl ( 
\frac{ f(\itwopi l_1, \dots, \itwopi l_{n-2}, Y_{n-1}) } 
{e^{-\yy{n-1} } - 1} \Biggr ) $$
$$ = 
\res_{Y_{n-1} = 0 } \Biggl ( \frac{ q_j(f)\Bigl (\itwopi \lj{1}, 
\dots, \itwopi \lj{n-2}, 
 Y_{n-1})  \Bigr ) }
{e^{-\yy{n-1} } - 1} \Biggr ).  $$
Here, we define an  involution
 $q_j: \liet \to \liet$ (for
$j = 1, \dots, n-1$) by 
$$q_j(X_1, \dots, X_j, \dots, X_{n-1}, X_n) = 
(X_1, \dots, X_{j-1},X_{n-1},X_{j+1}, \dots, 
X_{n-2},X_{j},X_n),$$
and the integers $\lj{1}, \dots, \lj{n-1} $ are defined by the 
equation
\beq \label{e:lj} q_j(X)|_{X = (\itwopi l_1, \dots, 
\itwopi l_{n-2}, Y_{n-1}) } 
= (\itwopi \lj{1}, \dots, \itwopi \lj{n-2}, 
\itwopi \lj{n-1} + Y_{n-1}) . \eeq
\end{prop}
\Proof For $j \le n-2$, the involution 
 $q_j$ is given in the 
coordinates $(Y_1, \dots, Y_{n-1})$ by 
$q_j: (Y_1, \dots, Y_{n-1}) \mapsto (Y'_1, \dots, Y'_{n-1}) $ where
$Y'_k = Y_k$ for $k \ne j-1, j, n-2, n-1$ and
\beq \label{2.009}Y'_{j-1} = Y_{j-1} + \dots + Y_{n-2}, \eeq
\beq \label{2.0010} Y'_j = - \sum_{j \le k \le n-2} Y_k, \eeq
\beq \label{2.0011} Y'_{n-2} = -  \sum_{j \le k \le n-3} Y_k, \eeq
\beq \label{2.0012} Y'_{n-1} =  Y_j + \dots + Y_{n-1}. \eeq
For $j = n-1$, $q_j$ is the identity map.  Notice that
$Y_{n-1}' $ is the only one of the transformed coordinates that
involves $Y_{n-1}.$
Notice also that $q_j$ takes $p_j$ to a point where 
$Y'_{n-1} = 0 . $ 

We now examine the image of $\linestd$ under $q_j$. The 
 integers $\lj{1}, \dots, \lj{n-1}$ were defined by the 
equation (\ref{e:lj}): in fact
 $\lj{k} = l_k$ for $k \ne j-1, j, n-2, n-1$ and 
\beq \label{2.0009}\lj{j-1} = l_{j-1} + \dots + l_{n-2}, \eeq
\beq \label{2.00010} \lj{j} = - \sum_{j \le k \le n-2} l_k, \eeq
\beq \label{2.00011} \lj{n-2} = - \sum_{j \le k \le n-3} l_k, \eeq
\beq \label{2.00012} \lj{n-1} =  l_j + \dots + l_{n-2}.  \eeq

\nc{\yj}{{ Y^{(j)}_{n-1} } }

We have that
$$
\res_{Y_{n-1} = - \itwopi (l_j + \dots + l_{n-2} ) } \Biggl ( 
\frac{ f(\itwopi l_1, \dots, \itwopi l_{n-2}, Y_{n-1}) } 
{e^{-\yy{n-1} } - 1} \Biggr ) = $$
$$
\res_{Y_{n-1} = - \itwopi (l_j + \dots + l_{n-2} ) }  \Biggl ( 
\frac{ f(\itwopi l_1, \dots, \itwopi l_{n-2}, Y_{n-1}) } 
{e^{-\yy{n-1} - \itwopi (l_j + \dots 
+ l_{n-2})  } - 1} \Biggr )  $$
(because $e^{\itwopi l_k } = 1$ for all $k = j, \dots, n-2$)
$$ 
= \res_{{\yj }= 0 } 
\Biggl ( 
\frac{ q_j(f)( \itwopi \lj{1}, \dots, \itwopi \lj{j-1} 
, \itwopi \lj{j}, 
\itwopi \lj{j+1}, \dots, 
\itwopi \lj{n-2}, \yj 
 }
{e^{-\yj } - 1} \Biggr ) $$
by the formulas (\ref{2.009} - \ref{2.0012})
where we have defined  $\yj = Y_{n-1} + \itwopi (l_j + \dots 
l_{n-2})$  so that $d\yj = dY_{n-1} $.
This completes the proof. \hfill $\square$

\begin{corollary} \label{c2.0}
We have
$$ \sum_{p \in \laregint} f(p) = $$
$$  \sum_{j = 1}^{n-2} 
\res_{Y_{n-1} = 0 } 
\Biggl ( 
\frac{ q_j(f)( \itwopi \lj{1}, 
\dots,
\itwopi \lj{n-2}, Y_{n-1} ) }
{e^{-Y_{n-1} } - 1} \Biggr ) +  $$
\beq \label{2.0055}
 \Res_{\yy{n-1} = 0  }
\frac{f(\itwopi \lb{1}, \dots, \itwopi \lb{n-2},
 \yy{n-1} ) }{ e^{- \yy{n-1} } - 1 },  \eeq
where the integers 
$\lj{1}, \dots, \lj{n-2} $ were defined by 
(\ref{2.0009} - \ref{2.00011}).
\end{corollary}
\Proof This follows by adding the results of Proposition 
\ref{p2.3} over all $j = 1, \dots, n-1$: on one side  this yields the 
 sum on the right hand side of (\ref{2.001})  (which 
according to (\ref{2.001}) is equal to 
$\sum_{p \in \laregint} f(p)$), and on the other side 
yields the sum on the right hand side of (\ref{2.0055}). 
\hfill $\square$

We shall complete the proof of  Lemma \ref{l2.04} by summing the equality
given in Corollary \ref{c2.0} over all possible 
$(l_1, \dots, l_{n-2})$ satisfying 
(\ref{2.cond}): the proof reduces to the following lemma.

\begin{lemma} 
In the notation of Proposition  \ref{p2.3}, 
$(\lj{1}, \dots, \lj{n-2} ) \in \weightl_{\rm reg} (\liet_{n-1})$. Moreover
for any $(l'_1, \dots, l'_{n-2}) \in \weightl_{\rm reg}
(\liet_{n-1})$ there is 
exactly one 
sequence of integers $(l_1, \dots, l_{n-2})$ satisfying
(\ref{2.cond}) such that 
$$( \lj{1}, \dots,  \lj{n-2}) 
 =( l'_1, \dots,  l'_{n-2} ).$$
\end{lemma}
\Proof This follows 
immediately from the proof of Proposition \ref{p2.3} and the fact that
the restriction of $q_j$ to $\liet_{n-1}$ is given by the action of 
an element of the Weyl group $W_{n-1}$ and hence maps 
$ \weightl_{\rm reg} (\liet_{n-1})$ to itself bijectively.\hfill $\square$

This completes the proof of Lemma \ref{l2.04} and hence
of Theorem \ref{t2.04} and Proposition \ref{p:sz}.

\renorm 
\section{Residue formulas and nonabelian localization}

Let $(M,\omega)$ be a compact symplectic manifold
with a Hamiltonian action of a compact connected Lie
group $K$ with 
Lie algebra $\liek$. Let $\mu:M \to \lieks$ be a moment map for this
action.

The $K$-equivariant cohomology with complex coefficients
$H^*_K(M)$ of $M$ may be identified with the cohomology of  the
chain complex 
\beq \label{1.001} \Omega^*_K(M) = (S(\lieks) \otimes \Omega^*(M))^K \eeq
of equivariant differential forms on $M$, equipped
 with the differential\footnote{This
definition of the equivariant cohomology 
differential differs by a factor of $i$ from
that used in \cite{tdgr} but is consistent with that used in \cite{JK1}.}
\beq \label{1.0002}
(D\eta)(\xvec) = d(\eta (\xvec) )  - \iota_{\xvec^{\#}} (\eta(\xvec) ) \eeq
where $\xvec^{\#}$ is the vector field on $M$ generated by the action of 
$\xvec$ (see Chapter 7 of \cite{BGV}). Here $(\Omega^*(M),d)$ is the
de Rham complex of
differential forms on $M$ (with
 complex coefficients), and $S(\lieks)$ denotes the algebra
of  polynomial functions on the Lie algebra $\liek$ of $K$. An element 
$\eta \in \Omega^*_K(M)$ may be thought of 
as a $K$-equivariant polynomial function from $\liek$ 
to $\Omega^*(M)$, or alternatively as a family of differential forms
on $M$ parametrized by $\xvec \in \liek$. 
The standard definition of degree is used on $\Omega^*(M)$ and
degree two is assigned to elements of $\lieks$.

In fact as a vector space, though not in general as a ring,
when $M$ is a compact symplectic manifold with a Hamiltonian action 
of $K$ then $H^*_K(M)$ is 
isomorphic to $H^*(M) \otimes H^*_K$ where $H^*_K=\Omega^*_K({\rm pt})
=S(\lieks)^K$ is the
equivariant cohomology of a point (see \cite{Ki1} Proposition 5.8).

The map $\Omega^*_K(M) \to 
\Omega^*_K({\rm pt}) = S(\lieks)^K$ 
given by integration over $M$ passes to $\hk(M)$. Thus  for any 
$D$-closed element 
$\eta \in \Omega^*_K(M)$ representing a cohomology class $[\eta]$, 
there is a corresponding element 
$\int_M \eta \in \Omega^*_K({\rm pt})$ which depends only on 
$[\eta]$. The same is true for any 
$D$-closed element $\eta = \sum_j \eta_j$ 
 which is a formal 
 series of elements $\eta_j$ in $\Omega^j_K(M)$  without polynomial 
dependence on $\xvec$: we shall in particular consider terms
of the form 
$$\eta (\xvec) e^{ (\bom(\xvec))} $$ 
where 
$\eta \in \Omega^*_K(M)$ and 
$$\bom(\xvec) = \omega + \mu(\xvec) \in \Omega^2_K(M). $$
Here $\mu:M \to \lieks$ is
 identified in the natural way with a linear function
on $\liek$ with values in $\Omega^0(M)$. It follows directly from the
definition of a moment map\footnote{We follow the convention
that $d\mu(X)=\iota_{X^{\sharp}}\omega$; some authors
have $d\mu(X)=-\iota_{X^{\sharp}}\omega$.} that $D \bom = 0 $.  

 If $\xvec$ lies in $\liet$, the Lie algebra of 
a chosen maximal torus $T$ of $K$, then there is a formula for 
$\int_M \eta(\xvec)$ (the {\em   abelian localization 
formula} \cite{AB,BGV,BV1,BV2})
which depends only on the fixed point set of 
$T$ in $M$. It tells us that 
\beq \label{1.002}
 \int_M \eta(\xvec) 
= \sum_{F \in \calf} \int \frac{i_F^* \eta(\xvec)}
{e_F(\xvec)} \eeq
where $\calf$ indexes the components $F$ of the fixed point set of 
$T$ in $M$, the inclusion of $F$ in $M$ is denoted
by $i_F$ and $e_F $ 
$\in H^*_T(M)$ is the equivariant Euler class of the normal 
bundle to $F$ in $M$. In particular, applying (\ref{1.002})
with  $\eta $ replaced by the formal equivariant cohomology class
 $\eta e^{ \bom}$ 
we have 
\beq \label{1.003} h^\eta(\xvec) \eqdef
 \int_M \eta(\xvec)e^{ \bom(\xvec)}  
= \sum_{F \in \calf} h^\eta_F(\xvec), \eeq
where 
\beq \label{1.004}   \hfeta(\xvec) = 
e^{ \mu(F)(\xvec)}\int_F  \frac{i_F^* \eta(\xvec)  e^{ \omega} 
 }{e_F (\xvec) }. \eeq
Note that the moment map $\mu$ takes a 
constant value $\mu(F) \in\liets$ 
 on each 
$F \in \calf$, and that the integral in (\ref{1.004}) 
is a rational function of $\xvec$. 

We shall assume throughout that $0$ is a regular value of the moment
map $\mu:M \to \lieks$; equivalently the action of $K$ on $\zloc$ has only
finite isotropy groups. The reduced space
$$\mred = \zloc/K$$
is then a compact symplectic
 orbifold. The cohomology (with complex coefficients,
as always in this paper) $H^*(\mred)$ of this reduced
space is naturally isomorphic to the equivariant cohomology
$H^*_K(\zloc)$ of $\zloc$, and by Theorem 5.4 of \cite{Ki1}
the inclusion of $\zloc$ in $M$ induces a surjection on equivariant cohomology
$$\hk(M) \to \hk(\zloc).$$
Composing we obtain a natural surjection
$$\Ph: \hk(M) \to H^*(\mred)$$
which we shall denote by
$$\eta \mapsto \eta_0.$$
When there is no danger of confusion 
we shall use the same symbol for $\eta\in\hk(M)$ and any
equivariantly closed differential form in $\Omega^*_K(M)$ which represents it.
Note that $(\bom)_0\in H^*(\mred)$ is represented by
the symplectic form $\omega_0$ induced on $\mred$ by $\omega$.

\noindent{\bf Remark} Later we shall be working with not only the 
reduced space $\xred = \mu^{-1}(0)/K$ with respect to the action of 
the nonabelian group $K$, but  also 
 $\mu^{-1}(0)/T$ and 
$\mredt{\xi}$  $ = \mu_T^{-1}(\xi)/T$ 
for regular values $\xi$ of the $T$-moment map $\mu_T$ which is the 
composition of $\mu$ with restriction from $\lieks$ to $\liets$. 
We shall use the same notation $\eta_0$ for the image of $\eta$
under the surjective homomorphism
$\Ph$ for whichever of the spaces
$\mu^{-1}(0)/K, $ $\mu^{-1}(0)/T $ or
$\mu_T^{-1}(0)/T $ we are considering, and the 
notation $\eta_{\xi}$ if we are working with $\mu_T^{-1}(\xi)/T$. 
It should be clear from the context which version of the map 
$\Ph$ is being used.

\bigskip

The main result (the residue formula, Theorem 
8.1) of \cite{JK1} gives a formula for the evaluation on the fundamental
class $[\mred]\in H_*(\mred)$, or equivalently
(if we represent cohomology classes by differential forms) the integral
over $\mred$, of the image
$\eta_0 e^{\omega_0}$ in $H^*(\mred)$ of any
formal equivariant cohomology class on $M$ of the type $\eta e^{\bom}$ where
$\eta \in \hk(M)$. 

\begin{theorem} \label{t4.1}{\bf (Residue 
formula, \cite{JK1} Theorem 8.1)} Let $\eta \in 
\hk(M) $ induce $\eta_0 \in H^*(\xred)$. Then 
\beq \label{jk81}   \eta_0 e^{{}\omega_0} [\xred] 
 = n_0 {C_K}  \res \Biggl ( 
\nusym^2 (X)
 \sum_{F \in \calf} \hfeta(X) [d X] \Biggr ), \eeq
where the constant\footnote{This constant differs by a factor of 
$(-1)^s (2\pi)^{s-l} $                                 
from that of \cite{JK1} Theorem 8.1. The reason 
for the factor of 
$ (2\pi)^{s-l} $                                 
is that in this paper we
shall adopt the
convention that  weights $\beta \in \liets$ send the integer lattice $\L^I =
{\rm Ker}(\exp:\liet \to T)$ to $\ZZ$ rather than to $2\pi \ZZ$, and that the
roots of $K$ are the nonzero weights of its complexified adjoint action. In
\cite{JK1} the roots send $\L^I$ to $2 \pi \ZZ$. The reason for the
factor of $(-1)^s$ is an error in Section 5 of \cite{JK1}.
In the last paragraph of p.307 of \cite{JK1} the appropriate form
to consider is $\prod_{j=1}^s (\theta^j dz'_j)$, and since 1-forms
anticommute this is $(-1)^s/i^s$ times the term in $\exp (idz'(\theta))$
which contributes to the integral (5.4) of \cite{JK1}.
 The constant also differs by
a factor of  $i^s$ from that of \cite{JK3} Theorem 3.1, because in that paper
the convention adopted on the equivariant cohomology differential is that of
\cite{tdgr}, not that of \cite{JK1}.}
$C_K$ is defined by 
\beq \label{4.001} C_K = \frac{(-1)^{s+n_+}}{ |W| \vol(T)}, \eeq
and $n_0$ is the order of the stabilizer in $K$ of a generic 
point\footnote{Note that in \cite{JK1} and \cite{JK2} $n_0$ is stated
incorrectly to be the order of the subgroup of $K$ which acts trivially
on $\zloc$ (i.e. the kernel of the action of $K$ on $\zloc$): see the
correction in Section 3 of \cite{JK3}. When $K=T$ is abelian, however,
the stabilizer in $K$ of a generic point of $\zloc$ is equal to the
kernel of the action of $K$ on $\zloc$. Moreover since the coadjoint action
of $T$ on $\liets$ is trivial, when $K=T$
this subgroup acts trivially on the normal
bundle to $\zloc$ in $M$ and hence is the kernel of the action of
$K$ on $M$.}
of $\zloc$.
\end{theorem}
In this formula $|W|$ is the order of the Weyl group $W$ of $K$, and we
have introduced $s = \dim K$ and 
$l  = \dim T$, while
$n_+ = (s-l)/2$ is the 
number of positive roots. The measure $[dX]$ on $\liet$ and volume $\vol(T)$ of
$T$ are obtained from the
 restriction of a fixed invariant inner product on $\liek$,
which is used to identify $\lieks$ with $\liek$ throughout.
Also, $\calf$ denotes the set of components of the fixed point
set of $T$, and if $F$ is one 
of these components then the meromorphic  function 
$\hfeta$ on $\liet \otimes \CC$ is defined by 
(\ref{1.004}). The polynomial
$\nusym: \liet \to \RR$ is defined by 
$$\nusym(X) = \prod_{\g > 0 } \g(X),$$
where $\g$ runs 
over the positive roots of $K$. Note that it would perhaps be more
natural to combine $(-1)^{n_+}$ from the constant $C_K$ with
$\nusym^2(X)$ and replace them by the product
$$\prod_{\g} \g(X)$$
of all the positive and negative roots of $K$.

The formula (\ref{jk81}) was
 called a residue formula in \cite{JK1} because the 
quantity $\res$  
(whose general definition 
was given in Section 8 of \cite{JK1})
can be expressed as a multivariable
residue\footnote{An alternative
 definition  in terms of iterated 1-variable residues
is given in Section 3 of \cite{JK3}.}, whose
domain is a class of  meromorphic differential forms on $\liet\otimes\CC$.
It is a linear map, but in order to apply it to individual terms in the residue
formula some choices have to be made which do not affect the residue
of the whole sum. Once the choices have been made one finds that many
of the terms in the sum contribute zero, and the formula can be rewritten as
a sum over a certain subset $\calf_+$ of the set $\calf $ of components of the
fixed point set $M^T$. When the rank of $K$ is one and $\liet$ is identified
with $\RR$, we can take
$$\calf_+ = \{ F \in \calf: \mu_T(F) > 0 \}. $$ 
In this paper we shall be particularly
interested in the case where $K$ has rank one, for which the results are as 
follows.  
\begin{corollary} \label{c4.2} {\bf \cite{Kalkman,Wu,JK1}}
In the 
situation of Theorem \ref{t4.1}, let $K = U(1)$.
 Then
$$ \eta_0 e^{{}\omega_0} [\xred] = - n_0
 \res_{X=0} \Bigl (  \sum_{F \in \calf_+} \hfeta(X) 
\Bigr )  $$ 
where $n_0$ is the order of the stabilizer in $K$ of a
generic point in $\zloc$.
Here, the meromorphic  function $\hfeta $ on  $ \CC$ is defined 
by (\ref{1.004}), and $\res_{X=0} $ denotes the coefficient of $1/X$,  where 
$X \in \RR$ has been identified with $2 \pi i X \in \liek$.  
\end{corollary}

\noindent{\bf Remark} The notation $\res_{X=0}$ is intended
to indicate the variable $X$ with respect to which the residue is
calculated, as well as the point 0 at which the residue is
taken, so that, for example, $\res_{-X=0} f(X) = -\res_{X=0} f(X)$.
It would perhaps be more natural to use the notation
$\res_{X=0} f(X)dX$, but we shall have numerous formulas
involving iterated residues of this type which would then become
too long and unwieldy.

\begin{corollary} \label{c4.3}{\bf(cf. \cite{JK1}, Corollary 8.2)}
In the situation of Theorem \ref{t4.1},
let $K = SU(2)$. Then 
$$ \eta_0 e^{{}\omega_0} [\xred] = 
\frac{n_0}{2} \res_{X=0} \Bigl ( (2 X)^2  \sum_{F \in \calf_+} \hfeta(X)
 \Bigr )  . $$ 
Here, $n_0$, $\res_{X=0}$, $\hfeta$ and $\calf_+$ are as in Corollary 
\ref{c4.2}, and $X \in \RR $ has been identified with
${\rm diag} (2 \pi i, - 2 \pi i) X \in \liet$. 
\end{corollary}

\begin{rem} \label{omiteo} Note that if the degree of $\eta$ is equal to the
dimension of $\mred$ then
$$\eta_0 e^{\omega_0} [\mred] = \eta_0 [\mred].$$
Alternatively for $K=U(1)$ or $K=SU(2)$, if
we multiply $\omega$ and $\mu$ by a real
scalar $\e >0$ and let $\e$ tend to $0$ we obtain
$$ \eta_0 [\xred] = - n_0
 \res_{X=0} \Bigl (  \sum_{F \in \calf_+} \int_F \frac{i_F^*\eta(X)}{e_F(X)}
\Bigr )  $$ 
or
$$ \eta_0 [\xred] = \frac{n_0}{2}
 \res_{X=0} \Bigl
 ((2X)^2  \sum_{F \in \calf_+} \int_F \frac{i_F^*\eta(X)}{e_F(X)}
\Bigr ).  $$ 
\end{rem}

The results we have stated so far require the symplectic manifold $M$
to be compact, and this condition is not satisfied in the situation
in which we would like to apply them (in order to obtain formulas for
the intersection pairings in the cohomology of moduli spaces of
bundles over compact Riemann surfaces). Luckily there are other related
results due to Guillemin and Kalkman 
\cite{GK}, and independently Martin \cite{Ma},
which as we shall see can be generalized to noncompact symplectic manifolds.

Guillemin and Kalkman and Martin 
have approached the problem of finding a formula for 
$$\eta_0 [\mred] = \int_\mred \eta_0$$
 in terms
of data on $M$ localised near $M^T$ in a slightly different way
from that described above. As Guillemin and
Kalkman observe, it follows immediately from the residue formula that
if $\xi \in\liets$
 is a regular value of the $T$-moment map $\mu_T:M \to \liets$
which is sufficiently close to $0$ then
\beq \label{rtmt} \eta_0 [\mred] = \frac{(-1)^{n_+}n_0(\nusym^2\eta)_{\xi}}
{n_0^T|W|} [\mu_T^{-1}(\xi)/T]\eeq
where $n_0$ (respectively $n_0^T$) is the order of the stabilizer
in $K$ (respectively $T$) of a generic point of $\zloc$ (respectively
$\mu_T^{-1}(0)$) and
$\mu_T^{-1}(\xi)/T$ is the reduced space for the action of $T$ on $M$ 
with respect to the shifted moment map $\mu_T - \xi$. Also
$(\nusym^2 \eta)_{\xi}  \in H^*(\mu_T^{-1}(\xi)/T)$
is the image of $ \nusym^2 \eta$ under the surjection $\Ph:H^*_T(M) \to
H^*(\mred)$. Here $\eta\in \hk(M)$ and $\nusym \in S(\liets) = H^*_T$
are regarded as elements of  $H^*_T(M)$ via the natural identification
of $H^*_K(M)$ with the Weyl invariant part $(H^*_T(M))^W$ of
$H^*_T(M)$ and the natural inclusion $H_T^* \to H^*_T(M)$. 
Martin gives a direct proof of (\ref{rtmt}) without appealing to the residue 
formula, which shows also that for any $\xi$ sufficiently close to $0$
\beq \label{rtmt2} \eta_0 [\mred] = \frac{n_0(\nusym \eta)_{\xi}}
{n'_0|W|} [\mu^{-1}(\xi)/T]\eeq
where $n'_0$ is the order of the stabilizer in $T$ of a generic
point in $\zloc$, provided that $\mu^{-1}(\xi)/T$ is oriented
appropriately.

\begin{rem}\label{orient} The symplectic form 
$\omega$ induces an orientation on $M$,
and the induced symplectic
 forms on $\mred = \zloc / K$ and on $\mu_T^{-1}(\xi)/T$
induce orientations on these quotients. We have made a choice of positive Weyl
chamber for $K$ in $\liet$; this determines
 a Borel subgroup $B$ (containing $T$)
of the complexification $G$ of $K$, such
 that the weights of the adjoint action of
$T$ on the quotient $\lieg / \lieb$
 of the Lie algebra $\lieg$ of $G$ by the Lie algebra
$\lieb$ of $B$ are the positive roots of $K$. We
 then get an orientation of the flag
manifold $K/T$ by identifying it with the complex
 space $G/B$. Modulo the action of finite
isotropy groups we have a fibration 
$$\zloc/T \to \zloc/K$$
with fibre $K/T$; thus the symplectic
 orientation of $\zloc/K$ and the orientation of $K/T$
determined by the choice of Weyl chamber induce an
 orientation of $\zloc /T$. Since 0
is a regular value of $\mu$, if $\xi$ is sufficiently close to 0
there is a homeomorphism from $\zloc/T$ to $\mu^{-1}(\xi)/T$
induced by a $T$-equivariant isotopy of $M$, so
 we get an induced orientation of
$\mu^{-1}(\xi)/T$. This is the orientation of
 $\mu^{-1}(\xi)/T$ which we shall use.

Note that given a positive Weyl chamber we have another choice of orientation
on $\mu^{-1}(\xi)/T$ which is compatible with the symplectic orientation on 
$\mu_T^{-1}(\xi)/T$ and the orientation of the normal bundle 
to $\mu^{-1}(\xi)/T$
in $\mu_T^{-1}(\xi)/T$ induced by identifying it in the natural way with 
the kernel of the restriction map $\lieks \to \liets$, thence via the
fixed invariant inner product on $\liek$ with $\liek/\liet$
and thus finally with the complex vector space $\lieg/\lieb$ as above. 
Because we have used the inner product to identify
$\liek/\liet$ with its dual here, this orientation differs
from the one chosen above by a factor of $(-1)^{n_+}$ where $n_+$ is the number
of positive roots.
\end{rem}

\begin{prop}  \label{p:sm}
{\bf (Reduction to the abelian case)} {\sc [S. Martin] 
\cite{Ma}} 
If $T$ is a maximal torus of $K$ and $K$ acts effectively on $M$, then 
for any regular value $\xi$ of $\mu_T$ sufficiently
close to $0$ we have that 
$$ \int_{\mu^{-1}(0)/K} 
(\eta e^{ \bom} )_0  = \frac{n_0}{n'_0|W|}
\int_{\mu^{-1}(0)/T} 
(\nusym \eta e^{ \bom} )_0 
= \frac{n_0}{n'_0|W|}
\int_{\mu^{-1}(\xi)/T}
(\nusym \eta e^{ \bom} )_{\xi} $$
$$ = \frac{(-1)^{n_+}n_0}{n_0^T|W|}
\int_{\mu_T^{-1}(\xi)/T}
(\nusym^2 \eta e^{ \bom} )_{\xi} $$
where $n_0$ is the order of the stabilizer in $K$ of a generic
point of $\zloc$ and $n_0^T$ (respectively $n'_0$) is the order of
the stabilizer in $T$ of a generic point of $\mu_T^{-1}(0)$
(respectively $\zloc$). 
\end{prop}

\begin{rem} Note that $(-1)^{n_+}\nusym^2$ is the product of all the roots of
$K$, both positive and negative.
\end{rem}

Martin proves this result by considering the diagram
$$\begin{array}{cccccc}
 & \mu^{-1}(0)/T & \cong & \mu^{-1}(\xi)/T
 & \hookrightarrow & \mu_T^{-1}(\xi)/T \\
  & \downarrow  &   &   &   &  \\
\mred = &  \mu^{-1}(0)/K &  &  & &
\end{array}$$
where the homeomorphism from $\mu^{-1}(0)/T $ to $\mu^{-1}(\xi)/T $ is
induced by a $T$-equivariant isotopy
of $M$ (for $\xi$ sufficiently close to $0$). 
For simplicity we shall consider the case when $n_0=n'_0=n_0^T =1$.
As before we use a fixed invariant inner product
on $\liek$ to identify $\lieks$ with
 $\liek$, which splits $T$-equivariantly as the
direct sum of $\liet$ and its orthogonal
  complement $\liet^{\perp}$. The projection
of $\mu:M \to \lieks\cong\liek$ onto 
$\liet^{\perp}$ then defines a $T$-equivariant section of
the bundle $M \times \liet^{\perp}$ on $M$, which has equivariant Euler class 
$(-1)^{n_+}\nusym$
if we orient $\liet^{\perp} \cong \liek/\liet$ by identifying it with the 
dual of the complex vector space 
$\lieg/\lieb$ as in Remark \ref{orient}. Hence 
if $\xi$ is a regular value of $\mu_T$ then 
$\mu^{-1}(\xi)/T$ is a zero-section of the induced 
orbifold bundle
$\mu_T^{-1}(\xi) \times_T \liet^{\perp}$ on $ \mu_T^{-1}(\xi)/T$, 
whose Euler class is $(-1)^{n_+}\nusym_{\xi}$. 
Thus under the conventions for orientations described in Remark
\ref{orient}, evaluating the restriction to $\mu^{-1}(\xi)/T$ of an element of
$H^*(\mu_T^{-1}(\xi))/T$ on the fundamental class 
$[\mu^{-1}(\xi)/T]$ gives the same
result as multiplying by $(-1)^{n_+}\nusym_{\xi}$ and 
evaluating on the fundamental class $[\mu_T^{-1}(\xi)/T]$. 

Now Martin observes that since the natural map
$$\Pi:\mu^{-1}(0)/T \to \mu^{-1}(0)/K = \mred$$
is a fibration with fibre $K/T$, modulo the action of finite isotropy groups
which act trivially on cohomology with complex coefficients, and since the
Euler characteristic of $K/T$ is
 nonzero (in fact it is the order $|W|$ of the Weyl
group of $K$), the evaluation of a cohomology class $\eta_0 \in H^*(\mred)$
on $[\mred]$ is given by the 
evaluation of an associated cohomology class on $[\zloc/T]$. More precisely we 
have
\beq \label{martin} \eta_0 [\mred] = \frac{e(V)}{|W|} 
\Pi^*(\eta_0) [\mu^{-1}(0)/T] \eeq
where $e(V)$ is the Euler class of the vertical subbundle of the tangent bundle
to $\mu^{-1}(0)/T$ with respect to the fibration $\Pi$. As this Euler class is
induced by $\nusym$ under the orientation conventions of Remark
\ref{orient},  this completes the proof.

\begin{rem} \label{noncom} \label{leg}
In this proof we saw that $\nusym_{\xi}$ is the cohomology class 
in $H^*(\mu_T^{-1}(\xi)/T)$ which is Poincar\'{e} dual to the homology
class represented by $\mu^{-1}(\xi)/T$. Thus $\nusym_{\xi}$ may be represented
by a closed differential form on $\mu_T^{-1}(\xi)/T$ with support in an
arbitrarily small neighbourhood of $\mu^{-1}(\xi)/T$. If we interpret
$\nusym_{\xi}$ in this way, Martin's proof of Proposition \ref{p:sm}
is valid even when $M$ is noncompact and has singularities, provided that
for $\xi$ near $0$ the subset $\mu^{-1}(\xi)$ is compact and does not
meet the singularities of $M$.

Note also that $K$ and hence $T$ act with at most finite isotropy groups on a
neighbourhood of $\mu^{-1}(0)$ in $\mu_T^{-1}(0)$, and so $\mu_T^{-1}(0)/T$
has at worst orbifold singularities in a neighbourhood of $\mu^{-1}(0)/T$. This
means that in Proposition \ref{p:sm}
 we do not need to perturb the value of the 
$T$-moment map $\mu_T$ from $0$ to a nearby regular value $\xi$ if, as above,
we represent $\nusym_{0}$ by a differential form on $\mu_T^{-1}(0)/T$ with
support in a sufficiently small neighbourhood of $\zloc/T$.
\end{rem}

This result reduces the problem
 of finding a formula for $\eta_0 [\mred]$ in terms
of data on $M$ localized near $M^T$ to the case when $K=T$ is itself a torus.
Guillemin and Kalkman, and independently
 Martin, then follow essentially the same
line. This is to consider the change in
$$\eta_{\xi} [\mu_T^{-1}(\xi)/T],$$
for fixed $\eta \in H^*_T(M)$, 
as $\xi$ varies through the regular values of $\mu_T$.
This is sufficient, if $M$ is a compact symplectic
 manifold, because the image $\mu_T(M)$ is bounded, so if $\xi$ is 
far enough from $0$ then $\mu_T^{-1}(\xi)/T$ is empty and thus
$\eta_{\xi} [\mu_T^{-1}(\xi)/T]=0.$

More precisely, the convexity theorem of Atiyah \cite{Aconv} and Guillemin and
Sternberg \cite{GSconv} tells us that the image $\mu_T(M)$ is a convex
polytope; it is the convex hull in $\liets$ of the set
$$\{ \mu_T(F) : F\in\calf\}$$
of the images $\mu_T(F)$ (each a single 
point of $\liets$) of the connected components $F$ of the fixed point
set $M^T$. This convex polytope is divided by codimension-one \lq\lq walls''
into subpolytopes, themselves convex hulls of subsets of
$\{ \mu_T(F) : F\in\calf\}$, whose interiors consist entirely of regular
values of $\mu_T$. When $\xi$ varies in the interior of one of these
subpolytopes there is no change in $\eta_{\xi} [\mu_T^{-1}(\xi)/T],$
so it suffices to understand what happens as $\xi$ crosses a
codimension-one wall.

Any such wall is the image $\mu_T(M_1)$ of a connected component
$M_1$ of the fixed point set of a circle subgroup $T_1$ of $T$. The
quotient group $T/T_1$ acts on $M_1$, which is a symplectic submanifold of $M$,
and the restriction of the moment map $\mu_T$ to $M_1$ has an orthogonal
decomposition
$$\mu_T|_{M_1} = \mu_{T/T_1} \oplus \mu_{T_1}$$
where $\mu_{T/T_1}: M_1 \to (\liet/\liet_1)^*$ is a moment map for the action
of $T/T_1$ on $M_1$ and $\mu_{T_1}:M_1 \to \liets_1$ is constant (because
$T_1$ acts trivially on $M_1$). If
 $\xi_1$ is a regular value of $\mu_{T/T_1}$ then
we have a reduced space
$$(M_1)_{{\rm red}} = \mu_{T/T_1}^{-1}(\xi_1 )/ (T/T_1).$$
Guillemin and Kalkman show that if $T$ acts effectively on $M$ (or 
equivalently if $n_0^T =1$; see Footnote 9) 
then, 
for an appropriate choice
of $\xi_1$, the change in $\eta_{\xi} [\mu_T^{-1}(\xi)/T]$
 as $\xi$ crosses the wall
$\mu_T(M_1)$ can be expressed as
$$({\rm res}_{M_1} (\eta))_{\xi_1} [(M_1)_{{\rm red}}]$$
for a certain residue operation (see Footnote 11 below)
$${\rm res}_{M_1} : H^*_T(M) \to H^{*-d_1}_{T/T_1}(M_1)$$
where $d_1 = {\rm codim} M_1 -2$. (Of course care is needed here about
the direction in which the wall is crossed; this can be resolved by a careful
analysis of orientations). By induction on the dimension of $T$ this gives
a method for calculating  $\eta_{\xi} [\mu_T^{-1}(\xi)/T]$ in terms of data
on $M$ localized near $M^T$.

It is easiest to see how this version of localization is related
to the residue formula of \cite{JK1}
in the special case when $K=T=U(1)$. In this case
$$\Omega_T^*(M) \cong \CC[X] \otimes \Omega^*(M)^T$$
is the tensor product of a polynomial ring in one variable $X$ (representing
a coordinate function on the Lie algebra $\liet$) with the algebra
of $T$-invariant de Rham forms on $M$. 
The Guillemin-Kalkman residue operation
$${\rm res}_{M_1} : H^*_T(M) \to H^{*-d_1}_{T/T_1}(M_1)$$
is then given in terms of the ordinary residue on $\CC$ by
$${\rm res}_{M_1} (\eta) = \res_{X=0} \frac{\eta|_{M_1}(X)}{e_{M_1}(X)}$$
where $\eta|_{M_1}(X)$  and
the equivariant Euler class $e_{M_1}(X)$ 
of the normal bundle
to $M_1$ in $M$ are regarded as polynomials
in $X$ with coefficients in $H^*(M_1)$. More precisely  we formally decompose
this normal bundle (using the splitting principle if necessary) as a sum
of complex line bundles $\nu_j$ on which $T$ acts with nonzero weights
$\beta_j \in \liets \cong \RR$, and because
$c_1(\nu_j) \in H^*(M_1)$ is nilpotent  we can express
$$\frac{\eta|_{M_1}(X)}{e_{M_1}(X)} = \frac{\eta|_{M_1}(X)}{\prod_j 
(c_1(\nu_j) + \beta_j X)} =  
\frac{\eta|_{M_1}(X)}{\prod_j (\b_j X)} \prod_j \Bigl(1+
\frac{c_1(\nu_j)}{\b_j X}\Bigr)^{-1}$$
as a finite Laurent series in $X$ with coefficients in $H^*(M_1)$. Then
${\rm res}_{M_1}(\eta)$ is simply the
 coefficient of $1/X$ in this expression\footnote{When
the dimension $l$ of $T$ is greater than one the
Guillemin-Kalkman residue operation
$${\rm res}_{M_1} : H^*_T(M) \to H^{*-d_1}_{T/T_1}(M_1)$$
is defined in almost exactly the same way, by choosing a coordinate system
$X=(X_1,\ldots,X_l)$ on $\liet$ where $X_1$ is a coordinate on $\liet_1$, and
taking the coefficient of $1/X_1$ in
 $\frac{\eta|_{M_1}(X)}{e_{M_1}(X)} $ expanded
formally as a Laurent series in
 $X_1$ with coefficients in $\CC[X_2, \ldots, X_l]
\otimes \Omega^*(M)^T$.}.
Since $T_1 = T$ acts trivially on $M_1$, we have $M_{1,{\rm red}} = M_1$ and
$M_1$ is a connected component of the fixed point set $M^T$, i.e.
$M_1 \in \calf$. Therefore
$$({\rm res}_{M_1} (\eta))_{\xi_1} [(M_1)_{{\rm red}}]
= \res_{X=0} \int_{M_1} \frac{\eta|_{M_1}(X)}{e_{M_1}(X)}.$$
Of course as $K=T=U(1)$ the convex polytope $\mu_T(M)$ in $\liets \cong \RR$
is a closed interval, divided into subintervals by the points
$\{\mu_T(F) : F \in \calf\}$. Thus the argument of Guillemin and Kalkman just
described, amplified by some careful consideration of orientations, tells us
that if $\xi >0$ is a regular value of $\mu_T$ and $n_0^T =1$ 
then the difference
$$ \eta_{\xi}[\mu_T^{-1}(\xi)/T] - \eta_0[\mu_T^{-1}(0)/T]  $$
can be expressed as
\beq \label{gk} \sum_{M_1 \in \calf: 0<\mu_T(M_1)<\xi} 
{\rm res}_{M_1} (\eta) [M_1]
=  \res_{X=0} \sum_{F\in\calf: 0<\mu_T(F) < \xi}  
\int_F \frac{i_F^*\eta(X)}{e_F(X)}. \eeq
If we take $\xi > {\rm sup} (\mu_T(M))$ then this gives the same result
as Corollary \ref{c4.2} (cf. Remark \ref{omiteo}).

\begin{prop}  \label{p:gkm}
{\bf ( Dependence  of symplectic quotients on parameters)} 
{\sc Guillemin-Kalkman \cite{GK} ; S. Martin \cite{Ma} }
If $K=T = U(1)$ and $n_0^T$ is the order of the stabilizer
in $T$ of a generic point of $\mu_T^{-1}(0)$ then\footnote{The
convention of Guillemin and Kalkman for the sign of the moment
map differs from ours (see Footnote 8). This accounts for a difference
in sign between their formula and ours.}
$$ 
\int_{\mu_T^{-1}(\xi_1)/T } (\eta e^{ \bom} )_{\xi_1} -
\int_{\mu_T^{-1}(\xi_0)/T } (\eta e^{ \bom} )_{\xi_0} = 
n_0^T \sum_{F \in \calf: \xi_0 < \mu_T(F) < \xi_1} 
{\rm Res}_{X=0} e^{ \mu_T(F)X}
\int_F \frac{\eta(X) e^{ \omega}  }{e_F(X) } . $$
where $X \in \CC$  has been identified with $2 \pi i X \in \liet
\otimes \CC$ and $\xi_0 < \xi_1$ are two regular values of  the moment map.

\end{prop}

\begin{rem} \label{arm} As we have already noted
these results can be deduced easily
 from the residue formula of \cite{JK1} when $M$ is a
compact symplectic manifold. However the proof of Proposition \ref{p:gkm},
just like that of Proposition \ref{p:sm} (see Remark \ref{noncom}),
can be adapted to apply in circumstances when
   $M$ is not compact and the residue formula of \cite{JK1}  is not valid.
   Indeed, as Guillemin and Kalkman observe, in the case when
$K=T=U(1)$ the basis of their argument  applies to any compact
oriented $U(1)$-manifold  $Y$ with boundary such that the action of $T=U(1)$
on the boundary $\partial Y$ is locally free. 
Let us suppose for simplicity that $T$
acts effectively on $M$ (i.e. that $n_0^T=1$; see Footnote 9)
and let $\zeta$ be a $U(1)$-invariant 
de Rham one-form on $Y-Y^T$ with the property
that $\iota_{v}(\zeta) = 1$, where the vector field $v$ is the 
infinitesimal generator of the $U(1)$-action. Guillemin and Kalkman 
showed that, at the level of forms, the map 
$\Ph:H_T^*(Y)\to H^*(\partial Y/T)$ which is the composition of
the restriction map from $H_T^*(Y)$ to $H_T^*(\partial Y)$ with the inverse of
the canonical isomorphism $H_T^*(\partial Y) \to H^*(\partial Y/T)$
is given by
$$\Ph(\eta) = \res_{X=0} \iota_{v}(\frac{\zeta\eta}{X-d\zeta})$$
(see (1.18) of \cite{GK}, noting that Guillemin and Kalkman have
a different convention for the equivariant cohomology differential,
which accounts for the minus sign). If tubular neighbourhoods
$U_1,\ldots,U_N$ of the components
 $F_1,\ldots,F_N$ of the fixed point set $Y^T$
are removed from $Y$, then Stokes' theorem can be applied to the
manifold with boundary $Y-\bigcup_{j=1}^N U_j$ using the formal
identity
$$D(\frac{\zeta\eta}{X-d\zeta}) = \eta$$
on $Y - \bigcup_{j=1}^N U_j$ to give, after using the fact that
$\int_{\partial Y} \alpha = \int_{\partial Y/T} \iota_v(\alpha)$
and taking residues at $X=0$, the formula
$$\int_{\partial Y/T} \Ph(\eta) 
= \res_{X=0} \sum_{j=1}^N \int_{F_j} \frac{\eta|_{F_j}(X)}{e_{F_j}(X)}$$
where 
$e_{F_j}$ is the equivariant Euler class of the normal
bundle to $F_j$ in $Y$. 

The formula  of Proposition \ref{p:gkm}
comes directly from this when the manifold with boundary $Y$ is
$\mu_T^{-1}[\xi_0,\xi_1]$ for a moment map $\mu_T:M \to \liets\cong\RR$
with regular values
 $\xi_0<\xi_1$, but there is no need for $\mu_T$ to be a moment
map or for $M$ to have a symplectic structure for the formula to be valid.
It is enough for $\mu_T:M \to \RR$ to be a smooth $T$-invariant map with
regular values $\xi_0 < \xi_1$ such that $T$ acts freely on the intersections
of $\mu_T^{-1}(\xi_0)$
and $\mu_T^{-1}(\xi_1)$ with the 
support of the equivariant differential form $\eta$.
There is also no need to assume that $M$ is compact; it suffices to
suppose  that
$\mu_T:M\to  \RR$ is a proper map. Indeed, the assumption
that $\mu_T$ is proper can itself be weakened; the same proof applies
provided only
that the intersection of $\mu_T^{-1}[\xi_0,\xi_1]$ with the support
of the equivariant differential form $\eta$ is compact.

\end{rem}

\renorm 
\section{Extended moduli spaces}

In \cite{ext} certain 
\lq\lq extended moduli spaces'' of flat connections on a compact
Riemann surface with one boundary component are studied. They
have natural symplectic structures, and
can be used to exhibit the moduli spaces $\mnd$ of interest to
us as finite-dimensional symplectic quotients or reduced spaces.
Our aim is to obtain Witten's formulas for intersection pairings
on $H^*(\mnd)$ by applying nonabelian localization to these extended
moduli spaces. They have a gauge-theoretic description (cf. the introduction
to this paper), but
we shall use a more
 concrete (and  entirely finite dimensional) characterization given 
in \cite{ext}.  

The space with which we want to work is defined by 
\beq \label{4.1} 
\mc = (\epsr{K} \times \epc)^{-1} (\bigtriangleup) \subset \Hom (\FF, K) 
\times \liek,
\eeq
where $\FF$ is the free group on $2g$ generators 
$  \{x_1, \dots, x_{2g} \}$; we identify $\FF$ with the fundamental
group of the surface $\Sigma$ with one point removed, in such a way
that $x_1, \dots, x_{2g} $ correspond to the generators $\alpha_1,
\dots, \alpha_{2g}$ of $H_1(\Sigma,{\bf Z})$ chosen in Section 2.
Then
 $\epsr{K}: \Hom(\FF, K) \to K$ is the evaluation map on the relator
$r = \prod_{j = 1}^g [x_{j}, x_{j+g}]$ 
\beq \label{4.2} 
\epsr{K} (h_1, \dots, h_{2g} ) = \prod_{j = 1}^g [h_{j} , h_{j+g} ]. 
\eeq
The map $\epc: \liek \to K $ is defined by 
\beq \label{4.3} \epc(Y) = \cent \exp (Y),   \eeq
where the generator $\cent$ of the 
centre of $K$ was defined at (\ref{1.p1}) above. 
The diagonal in $K \times K $ is denoted $\bigtriangleup$. 
The space $M(c)$ then has canonical projection maps 
$\proj_1, \proj_2 $ which make the following diagram commute: 
\beq \label{4.4} 
\begin{array}{lcr}
\xc  & \stackrel{\proj_2}{\lrar} & \liek \\
\scriptsize{\proj_1}
\downarrow & \phantom{\stackrel{aaaa}{\lrar} } & \downarrow \scriptsize{e_c} 
 \\
\homfk  & \stackrel{\epsr{K} }{\lrar} & K\\ \end{array} \eeq 
In other words, $\xc$ is the fibre product of $\homfk $ and $\liek$ 
under the maps $\epsr{K}$ and $\epc$. The action of $K$ on $\xc$ 
is given by the adjoint actions on $K$ and $\liek$. 
The space $\xc$ has the following properties (see 
\cite{ext} and \cite{J1}):

\begin{prop} \label{p0}

{\bf
 (a)} The space $\xc$ is smooth near all $(h, \L) \in 
\homfk \times \liek$ for which the linear space 
$z(h) \cap \ker (d \exp)_\L \ne \{0 \} $. Here, $z(h) $ is the Lie
algebra of the stabilizer  $Z(h) $ of $h$.

{\bf (b) } There is a $K$-invariant 2-form $\omega$ on 
$\homfk \times \liek$ whose restriction to $\xc$ is closed and which
defines a nondegenerate bilinear form on the 
Zariski tangent space to $\xc$ at  every $(h, \L)$
in an open dense subset of $\xc$ containing $\xc \cap (K^{2g} \times \{0\})$. 
Thus the form $\omega$ gives rise to a symplectic structure
on this open subset of  $\xc$.

{\bf (c) } With respect to the symplectic structure given by the 
2-form $\omega$, a moment map $\mu: \xc\to \lieks$ for the action of 
$K$ on $\xc$ is given by the restriction to $\xc$ of
$- \proj_2$, where 
$\proj_2: \xc \to \liek$ is the projection map to 
$\liek$  (composed
with the canonical isomorphism $\liek \to \lieks$ given by the invariant 
inner product on $\liek$). 

{\bf (d) } The space $\xc$
 is smooth in a neighbourhood of 
$\mu^{-1} (0).$ 

{\bf (e)} The symplectic quotient $\xred = \xc\cap\mu^{-1}(0)/K$ can be 
naturally identified with $\e_K^{-1}(c)/K = \mnd$.

\end{prop}

\begin{rem} We shall also use $\mu$ to denote the map
$$\mu:K^{2g}\times \liek \to \liek$$
defined by
$$\mu(h,\L) =  -\L,$$
even though it is only its restriction to $M(c)$ which is a moment map in any
obvious sense. That is why we write $\xc\cap\zloc/K$ instead of
$\zloc/K$ in (e) above.
\end{rem}

\begin{rem} \label{r:fibprod}
Using our description (\ref{4.4})
of $\xc$ as a fibre product, 
it is easy to identify the components $F$ 
of the fixed point set of the action of $T$.
We examine the fixed point sets of the action of  $T$ on 
$\homfk$ and
$\liek$ and find 
\beq \label{4.5} 
\begin{array}{lcr}
\xc^T   & \stackrel{\proj_2}{\lrar} & \liet \\
\scriptsize{\proj_1}
\downarrow & \phantom{\stackrel{aaaa}{\lrar} } & \downarrow \scriptsize{e_c} 
 \\
{\rm Hom}(\FF, T)
  & \stackrel{\epsr{K} }{\lrar} & 1 \in T \\ \end{array} \eeq 
(Notice that $\epsr{K}$ sends ${\rm Hom}(\FF, T)$ to $1$ because 
$T$ is abelian.) Thus 
\beq \label{fixed} M(c)^T = {\rm Hom}(\FF, T) \times e_c^{-1}(1) = 
T^{2g} \times \{ \d - \tilde{c}: \onebl \d \in \intlat 
\subset \liet \} \eeq
where $\tilde{c} $ is a fixed element of $\liet $ for which 
$\exp \tilde{c} = c$. (Here, $\intlat$ denotes the 
integer lattice ${\rm Ker}(\exp) \subset \liet$.)
If we ignore the singularities of $M(c)$, this
description also enables us to find a plausible candidate for the equivariant 
Euler class $e_{\fd}$ of the normal bundle of each component 
$T^{2g} \times (\d  - \tilde{c}) $ in $M(c)^T$
(indexed by $\d \in \intlat$).
This should be simply the equivariant
Euler class of the normal bundle to
$T^{2g}$  in $K^{2g}$, implying that 
$e_{\fd}$ is in fact independent of $\d$ and is given by 
\beq \label{33} e_{\fd} (\xvec) = (\prod_{\g} \g)^g = ((-1)^{n_+ }
\nusym(\xvec)^2)^{g}. \eeq
The symplectic volume of the component
$F_\d $ is independent of $\d$ (indeed these components are 
all identified symplectically with $T^{2g}$):
 we denote the volume of $F_\d$ 
by $\int_F e^{\omega}$.
The constant  value taken by the  moment map 
$\mu_T  $ on the component $F = F_\d$ is given by $\tc - \d$.
\end{rem}

We shall need also the following property (proved in \cite{J2}):
\begin{prop} \label{abftil}
The generating classes $a_r$, $b_r^j$ and $f_r$ ($r = 2, \dots, n$, 
$j = 1, \dots , 2g$) extend to classes $\tar$, $\tbrj$ and 
$\tfr$ $ \in H^*_K(\xc)$.  
\end{prop}
Indeed, because of our conventions on the equivariant
differential, the construction of \cite{J2} (which will be
described at the beginning of Section 9) tells us that
the equivariant differential form $\tar \in \Omega^*_K(\xc)$
whose restriction represents the cohomology class $a_r \in H^*(\mnd)$
is $\tau_r(-X)$, where as above $\tau_r \in S^r(\lieks)^K$
$\cong \hk({\rm pt})$ is the invariant polynomial 
which is associated to the $r$th Chern class (see 
\cite{J2}).
Moreover $\tf_2$ is the extension $\bar{\omega}=\omega + \mu$ of the
symplectic form $\omega$ to an equivariantly closed differential
form (see \cite{J2} again).

Finally we shall need to work with the symplectic 
subspace $\emtc =\xc \cap  \mu^{-1}(\liet)$ 
of $\xc$, which is no longer acted on by $K$ but is acted on by $T$. 
The space $\emtc$ has an important periodicity property:

\begin{lemma} \label{l4.3} 
Suppose $\tran $ lies in the integer lattice
$\L^I =  {\rm Ker} (\exp )$ in $\liet$. 
Then there is a homeomorphism $s_\tran: 
K^{2g}\times \liek \to K^{2g} \times \liek$ 
defined by 
$$ s_\tran: (h, \L) \mapsto (h, \L + \tran) $$
which restricts to a homeomorphism $s_\tran: \emtc \to \emtc$.
\end{lemma}
\Proof This is an immediate consequence of the definition of
$\emtc$ and the fact that
$\exp(\L + \L_0) = \exp(\L)\exp(\L_0)$ when $\L$ and $\L_0$ commute.
\hfill $\square$

Let us examine the behaviour of the images in $\hht(\emtc)$ 
of these extensions 
$\tar$, $\tbrj$, $\tfr$ $ \in \hk(\mc)$ 
of the generating classes $a_r, b_r^j, f_r$ (see Proposition \ref{abftil})
under pullback under these
homeomorphisms $s_{\L_0}: \emtc \to \emtc$.
By abuse of language, we shall
refer to these images also as $\tar, \tbrj $ and $\tfr$. 
We noted above that the classes $\tar$ are the images in
$\hk(\mc)$ of
the polynomials $\tau_r(-X) \in \hk = S(\lieks)^K$ (cf. (\ref{9})). Moreover 
(by \cite{J2}, (8.18)) the classes
$\tbrj \in \hk(\mc)$ are of the form 
$\tbrj = \proj_1^* (\tbrj)_1 $ where $  (\tbrj)_1 \in \hk(K^{2g} )$ 
and $\proj_1: \mc \to K^{2g}$ is the projection in
(\ref{4.4}).   It follows that
$$\stran^* \tbrj = \tbrj$$
 and 
$$\stran^* \tar = \tar. $$
Furthermore we see from (8.30) of \cite{J2} that $\tf_2(X)$ is of 
the form 
\beq \label{8.3} \tf_2(X) = \proj_1^* f_2^1 + \langle \mu,X \rangle \eeq
where $f_2^1 \in H^*_K(K^{2g})$ and
 $\mu: \mc \to \liek $ is the moment map 
(which is the restriction to $\mc$ of minus the
projection $K^{2g} \times \liek \to \liek$: see 
Proposition \ref{p0}). 
 It follows from this 
that for any $\tran $ in the integer lattice $\L^I$ of $\liet$ (the kernel
of the exponential map),  
\beq \label{8.4} \stran^* \tf_2(X) = \tf_2(X) -
\langle \tran,X \rangle . \eeq

\renorm
\section{Equivariant Poincar\'{e} duals}

We are aiming to apply nonabelian localization to the extended moduli space
$M(c)$ defined in the previous section.
In order to overcome the problem that $M(c)$ is singular, instead
of working with integrals
 over $M(c)$ of equivariant differential forms, we shall
integrate over $K^{2g}\times \liek$ after first multiplying by a suitable
equivariantly closed differential form
 on $K^{2g}\times \liek$ with support near
$M(c)$ which can be thought of as representing the equivariant
Poincar\'{e} dual to $M(c)$ in $K^{2g}\times \liek$. So we need to construct
such an equivariantly closed differential form.

\begin{rem} In our earlier article \cite{JK2} covering the case when the
bundles have rank $n=2$, we overcame the problem of the singularities
of $M(c)$ in a slightly different way, by perturbing the central constant
$c\in SU(n)$ to a nearby element of the maximal torus $T$. This method can be 
generalized to cover the cases when $n>2$, but it seems a little more
straightforward to
use equivariant Poincar\'{e} duals, so we adopt the latter approach here.
\end{rem}

\begin{rem}
Related 
constructions of equivariantly closed
differential forms representing the Poincar\'e dual to 
a submanifold   appear already in the literature.\footnote{We thank
P. Paradan for pointing out that  the references cited below
contain such constructions.}
In  Kalkman's
paper \cite{Kalkman2} and Mathai-Quillen's paper \cite{MQ}, 
an equivariantly closed  differential form which is rapidly decreasing
away from a submanifold and represents the Poincar\'e dual to the
submanifold is given: such a form is often referred to as the
{\em Thom form}, as the cohomology class it represents is  the
 Thom class of the normal bundle to the submanifold.
 The forms constructed in \cite{Kalkman2} and \cite{MQ} are not
compactly supported: a 
construction of a compactly
supported equivariantly closed form representing the
Poincar\'e dual  of a submanifold is given in section 2.3 of 
\cite{DV}. For completeness, in this section
we  provide a construction of an equivariantly
closed form representing the Poincar\'e dual.
\end{rem}

First we consider  the simpler problem of constructing
an equivariant Poincar\'{e} dual to the origin in a one dimensional 
representation $\chi$ of a circle. If we did not need to find a form
with support near the origin we could represent the equivariant Poincar\'{e}
dual by $\chi$ itself, regarded as an equivariant differential form.
However compact support will be important later, so we need to be a little 
more careful.

\begin{lemma} \label{pd1} Let $T=U(1)$ act on $\CC$ via 
a weight $\chi:T \to U(1)$.
Then we can find an equivariantly closed differential 
form $\a_{\chi}\in\Omega_T^2(\CC)$
on $\CC$ with compact support arbitrarily close to 0, such that
$$\int_{\CC} \eta \a_{\chi}  = \eta|_0 \in H^*_T$$
for all equivariantly closed forms $\eta \in \Omega_T^*(\CC)$. Moreover
$\a_{\chi} \in \chi + D(\Omega^*_T(\CC))$, so that $\a_{\chi}$ represents
the same equivariant cohomology class on $\CC$ as $\chi$.
\end{lemma}
\Proof  
Let $X^{\sharp}$ 
denote the vector field on $\CC$ given by the infinitesimal
action of $X\in \liet$. There is a $T$-invariant closed
differential 1-form on
$\CC-\{0\}$, given in polar coordinates $(r,\theta)$
by $\frac{d\theta}{2\pi}$,  such that 
$\iota_{X^{\sharp}}(\frac{d\theta}{2\pi})$ is identically
equal to $\chi(X)$ for every $X\in\liet$. We can choose 
a smooth $T$-invariant
function $b:\CC \to [0,\infty)$ with support in 
an arbitrarily small neighbourhood
of 0 which is identically equal to 1 on some smaller neighbourhood of 0, and
let
$$\a_{\chi}(X) = \chi(X) + D((1-b)\frac{d\theta}{2\pi}) 
= \chi(X) + d((1-b)\frac{d\theta}{2\pi})
+ (b-1)\chi(X)$$
where $D$ is the equivariant differential defined at (3.2) and 
$d$ is the ordinary differential.  Then $\a_{\chi}$ is
equivariantly closed and
is zero outside the support of $b$. 

Suppose that $\eta\in\Omega_T^*(\CC)$ is equivariantly closed.
We wish to show that
$$\int_{\CC} \eta \a_{\chi}  = \eta|_0.$$
First we shall show that the integral 
$$\int_{\CC} \eta \a_{\chi} $$
is independent of 
the choice of the function $b$. 

If $\rho>0$ is sufficiently small and $R>0$ is sufficiently large,
then $b$ is identically equal to 1 on the disc $D_{\rho}$ centre
0 and radius $\rho$, and  $b$ is identically equal to 0 outside
the disc $D_R$ centre
0 and radius $R$. Then
$$\int_{\CC} \eta \a_{\chi}  = \chi\int_{D_{\rho}}\eta
+ \int_{D_R-D_{\rho}} \eta \a_{\chi} .$$
Now $\eta$ is a polynomial function from $\liet$ to the ordinary
de Rham complex $\Omega^*(\CC)$, so we can write
$$\eta = \eta^{(0)} + \eta^{(1)} + \eta^{(2)}$$
where $\eta^{(j)}$ is a polynomial function from $\liet$ to
$\Omega^j(\CC)$ for $j=0,1,2$. Similarly
$$\alpha_{\chi} = \alpha_{\chi}^{(0)} + \alpha_{\chi}^{(1)} + 
\alpha_{\chi}^{(2)}$$
where $\alpha_{\chi}^{(0)}=b\chi$, $ \alpha_{\chi}^{(1)} = 0$
and $\alpha_{\chi}^{(2)} = d((1-b)\frac{d\theta}{2\pi})$. Since
$D\eta = d\eta -\iota_{X^{\sharp}}\eta$ is zero, we have
$d\eta^{(0)} = \iota_{X^{\sharp}}\eta^{(2)}$. As any 2-form
on $\CC$ is a $C^{\infty}$ function on $\CC$ multiplied by the
nowhere vanishing 2-form given in polar coordinates
by $\frac{rd\theta dr}{2\pi}$, and since
$\iota_{X^{\sharp}}(\frac{rd\theta dr}{2\pi}) = \chi(X) r dr$, it follows
that
$$\chi(X) \eta^{(2)}(X) = \frac{d\theta}{2\pi} d\eta^{(0)}(X)$$
on $\CC-\{0\}$ where $d\theta$ is defined. Hence
$$\int_{D_R-D_{\rho}} \eta \a_{\chi}  = \int_{D_R-D_{\rho}} 
\eta^{(2)}\a^{(0)}_{\chi}  +  \eta^{(0)}\a^{(2)}_{\chi}$$
$$= \int_{D_R-D_{\rho}} b \frac{d\theta}{2\pi} d\eta^{(0)} + \eta^{(0)}
d((1-b)\frac{d\theta}{2\pi})$$
$$= -\int_{D_R-D_{\rho}} d(b \eta^{(0)} \frac{d\theta}{2\pi}) $$
$$ = \int_{\partial D_{\rho}} b \eta^{(0)} \frac{d\theta}{2\pi}
- \int_{\partial D_{R}} b \eta^{(0)} \frac{d\theta}{2\pi}$$
$$ = \int_{\partial D_{\rho}} \eta^{(0)} \frac{d\theta}{2\pi}$$
by Stokes' theorem, since $b$ is identically one on
$\partial D_{\rho}$ and identically zero on $\partial D_R$.
It follows that  
$$\int_{\CC} \eta \a_{\chi}=\chi \int_{D_{\rho}} \eta
+ \int_{\partial D_{\rho}} \eta^{(0)} \frac{d\theta}{2\pi}$$
is independent of the
choice of $b$. 

Now $\rho$ can be taken arbitrarily small, and $\chi \int_{D_{\rho}} \eta
\to 0$ as $\rho \to 0$. 
Moreover by continuity, for fixed $X \in \liet$ and any $\e >0$ we can choose
$\rho$ so that $\eta^{(0)}$ differs from $\eta^{(0)}|_0 = \eta|_0$ by at most
$\e$ on $D_{\rho}$. Then
$$|\int_{\partial D_{\rho}} \eta^{(0)} \frac{d\theta}{2\pi} - 
\eta^{(0)}|_0| = |\int_{\partial D_{\rho}} (\eta^{(0)}-\eta^{(0)}|_0)
 \frac{d\theta}{2\pi}| \leq \epsilon.$$
Thus
$\int_{\CC} \eta \a_{\chi} - \eta|_0 $
tends to zero as $\rho$ tends to 0. Since $\int_{\CC} \eta\a_{\chi}$ 
and $\eta|_0$ 
are independent of $\rho$ we deduce that
$$\int_{\CC}\eta \a_{\chi}  = \eta|_0 $$
as required.
$\square$

\begin{lemma} \label{pd1.5} Let $T$ be a torus acting trivially on $\RR$.
Then we can find an
 equivariantly closed differential form $\a_{0}\in\Omega_T^*(\RR)$
on $\RR$ with compact support arbitrarily close to 0, such that
$$\int_{\RR} \eta \a_0  = \eta|_0 \in H^*_T$$
for all equivariantly closed forms $\eta \in \Omega_T^*(\RR)$. 
\end{lemma}
\Proof  
We have $\Omega_T^*(\RR)=S(\liet^*)\otimes \Omega^*(\RR)$ and
$\eta \in S(\liet^*)\otimes \Omega^0(\RR)$
is equivariantly closed if and only if it is constant on $\RR$, so we can take 
$\a_{0}$ to be the standard volume form on $\RR$ multiplied by any bump
function compactly supported near $0$ with unit integral.
$\square$

\begin{corollary} \label{pd2} Let $T$ be a
 torus acting linearly on $\CC^n$ with
weights $\chi_1,\ldots,\chi_n$ and trivially on $\RR^m$. 
Then we can find an equivariantly closed
differential form $\a \in \Omega_T^{2n}(\CC^n\times\RR^m)$ on 
$\CC^n\times\RR^m$ with compact support
arbitrarily close to 0, such that
$$\int_{\CC^n\times\RR^m}  \eta \a = \eta_0 \in H_T^*$$
for all equivariantly closed forms
  $\eta\in \Omega_T^*(\CC^n\times\RR^m)$. Moreover
if $m=0$ then $\a \in \chi_1\ldots\chi_n + D(\Omega_T^*(\CC^n))$.
\end{corollary}
\Proof The action of $T$ on the copy
  of $\CC$ in $\CC^n$ on which it acts via the
weight $\chi_j$ factors through an action of $T/\ker \chi_j \cong U(1)$ (unless
$\chi_j = 0$ in which case we can replace $\ker\chi_j$ by any subtorus of
$T$ of codimension one). We can construct $\a_{\chi_j}\in \Omega_{U(1)}^*(\CC)$
as in Lemma \ref{pd1} and $m$ copies of $\a_0$ as in Lemma
\ref{pd1.5}, and then define $\a$ to be the wedge
product of the pullbacks of the $\a_{\chi_j}$ and $\a_0$ 
to $\O_T^*(\CC^n\times \RR^m)$ via the projections
of $\CC^n\times\RR^m$ to $\CC$ and $\RR$ and the 
homomorphisms $T\to U(1)$ induced by the weights $\chi_j$.
$\square$

Now we shall relax our assumption that $c$ is a central element of $K$, 
and assume only that $c \in T$. This will be important later when we
apply induction on $n$ (see Remark \ref{induct} below).

\begin{corollary} \label{pd3} Let $T$ be
the maximal torus of $K=SU(n)$ acting on $K$ by 
conjugation. If $c \in T$ then we
can find a $T$-equivariantly closed 
differential form $\a\in\O_T^*(K)$ on $K$ with support
arbitrarily close to $c$ such that
$$\int_K  \eta \a = \eta|_c \in H^*_T$$
for all $T$-equivariantly closed differential forms $\eta\in\O_T^*(K)$.
\end{corollary}
\Proof There is a $T$-equivariant diffeomorphism $\phi$ from a $T$-invariant
neighbourhood $U$ of 0 in the Lie algebra $\liek$ of $K$ to a $T$-invariant
neighbourhood $V$ of $c$ in $K$ given by
$$\phi (X) = c \exp(X).$$
By Corollary \ref{pd2} we can find $\tilde{\a} \in \O_T^*(\liek)$ with
arbitrarily small compact support contained in $U$, such that
$$\int_{\liek} \eta \tilde{\a} = \eta|_0 \in H^*_T$$
for all equivariantly closed forms $\eta\in \O^*_T(\liek)$. Then we can define
$\a $ to be $(\phi^{-1})^*(\tilde{\a})$. $\square$

Note that 
$$M(c)=
\Bigl \{   (h_1, \dots, h_{2g}, \L )  \in K^{2g} \times 
\liek \; : \; \prod_{j = 1}^g h_{2j-1} h_{2j} h_{2j-1}^{-1} h_{2j}^{-1} 
= c \exp (\L) \Bigr \} $$
can be expressed as $\mc = P^{-1}(c)$ where 
$P: K^{2g} \times \liek \to K$ is defined by 
$$
 P \Bigl ( h_1, \dots, h_{2g}, \L \Bigr )  
= \prod_{j = 1}^g h_{2j-1} h_{2j} h_{2j-1}^{-1} h_{2j}^{-1} 
\exp (-\L). $$ 

\begin{prop}
\label{defa} If $T$ is the maximal torus of
$K=SU(n)$ and $c \in T$ then there is a $T$-equivariantly closed differential 
form $\a\in\O^*(K^{2g}\times\liek)$ of degree $n^2 -1$
on $K^{2g}\times \liek$ with support contained in a neighbourhood of $M(c)$
of the form $P^{-1}(V)$ where $V$ is an arbitrarily small neighbourhood of
$c$ in $K$, such that
$$\int_{K^{2g}\times\liek} \eta\a = \int_{M(c)} \eta|_{M(c)} \in H^*_T$$
for any $T$-equivariantly closed form $\eta\in \O^*_T(K^{2g}\times\liek)$
for which the intersection of $P^{-1}(\bar{V})$ with the support of $\eta$ is
compact. 
\end{prop}
\Proof By Corollary \ref{pd3} we can find a $T$-equivariantly 
closed differential
form $\hat{\a} \in \O_T^*(K)$ on $K$ with support in $V$ such that
$$\int_K \eta \hat{\a} = \eta|_c \in H^*_T$$
for all $T$-equivariantly closed forms $\eta\in \O^*_T(K)$. 
Let $\a = P^*(\hat{\a})$;
by the functoriality of the equivariant pushforward map (cf. Section 3
of \cite{ABMM}) this has the properties we want. $\square$

\begin{rem} In fact if $V'$ is any neighbourhood of $c$ in $K$
containing $\bar{V}$ then we have
$$\int_{P^{-1}(V')} \eta\a = \int_{M(c)} \eta|_{M(c)} \in H^*_T$$
for any $T$-equivariantly closed
 form $\eta\in \O^*_T(P^{-1}(V'))$ on $P^{-1}(V')$
for which the intersection of $P^{-1}(\bar{V})$ with the support of $\eta$ is
compact.
\end{rem}

\begin{rem} As we are going to use Proposition \ref{defa} to convert integrals
over $M(c)$ into
integrals over $K^{2g}\times \liek$ 
(or at least over neighbourhoods of $M(c)$ in
$K^{2g}\times \liek$ of the form $P^{-1}(V)$ for 
arbitrarily small neighbourhoods $V$ of $c$
in $K$) we shall need to be able to 
extend $T$-equivariant cohomology classes $\eta$ on $M(c)$
to $T$-equivariant cohomology classes on neighbourhoods of  $M(c)$ in
$K^{2g}\times \liek$ of this form $P^{-1}(V)$. This 
will always be possible by the continuity
properties of cohomology (see e.g. \cite{Dold} VIII 6.18) 
because $\eta$ will always have compact support in $M(c)$; more precisely
we will in fact be converting integrals over $M(c)\cap(K^{2g} \times B)$
for compact subsets $B$ of $\liek$ into integrals over
$P^{-1}(V)\cap(K^{2g} \times B)$.
\end{rem}

Note that the centre $Z_n$ of $K=SU(n)$ is a finite group of
order $n$ which acts trivially on $K^{2g} \times \liek$.

\begin{lemma} \label{l:5.10}
Suppose that $c = \diag(c_1,\ldots,c_n) \in T$
is such that the product of no proper subsequence of $c_1,\ldots,c_n$
is 1. Then the quotient $T/Z_n$ of $T$ by the centre
$Z_n$ of $K=SU(n)$ acts freely on $P^{-1}(V) \cap \mu^{-1}(0)$
for any sufficiently small $T$-invariant neighbourhood $V$ of $c$ in $K$.
\end{lemma}
\Proof Suppose that $T/Z_n$ does not act 
freely on $P^{-1}(V)\cap\mu^{-1}(0)$.
Then 
there exist $t_1,\ldots,t_n\in \CC$, not all equal, such that
$t_1\dots  t_n =1$, and some element $(h,0)=(h_1,\ldots,h_{2g},0)$ 
of $P^{-1}(V)\cap
\mu^{-1}(0)$ fixed by $\diag(t_1,\ldots,t_n)$. 
Then each $h_j$ is block
diagonal with respect to the decomposition of $\{1,\ldots,n\}$ as the union
of $\{i:t_i = t_1\}$ and $\{i:t_i \neq t_1\}$,
which implies that 
$$P(h,0) = \matr{A}{0}{0}{B}$$
where $A$ and $B$ are products of commutators and hence satisfy
$\det A =1 = \det B$. The result follows.
\hfill $\square$

\begin{rem} \label{aa}
It follows from this lemma that  we can extend the definition of the
composition
$$\Phi: H^*_T(M(c)) \to H_T^*(M(c) \cap \mu^{-1}(0)) 
\cong H^*(M(c) \cap \mu^{-1}(0)/T)$$
to
$$\Phi: H^*_T(P^{-1}(V)) \to H_T^*(P^{-1}(V) \cap \mu^{-1}(0)) \cong 
H^*(P^{-1}(V) \cap \mu^{-1}(0)/T).$$
By 1.18 of \cite{GK} (see Remark \ref{arm} above), 
when $T=U(1)$ is a circle then
$\Ph$ is given on the level of forms by
$$\Ph(\eta) = \res_{X=0} \iota_v (\frac{\zeta \eta}{X-d\zeta})$$
where the vector field $v$ is the infinitesimal generator of the $U(1)$ action 
and $\zeta$ is a $U(1)$-invariant differential 1-form on
$P^{-1}(V)\cap\mu^{-1}(0)$ such that $\iota_v(\zeta)=1$. 
(Strictly speaking the residue is an invariant form on 
$P^{-1}(V)\cap\mu^{-1}(0)$ which descends to a form on
$(P^{-1}(V)\cap\mu^{-1}(0))/T$).
Thus when $T=U(1)$ we 
have
$$\int_{M(c)\cap\zloc/T} \Ph(\eta) = \int_{M(c)\cap\zloc} 
\res_{X=0} \frac{\zeta \eta}{X-d\zeta},$$
and it  follows that
 if $\a$ is defined as in Proposition \ref{defa} for $n=2$
and $V'$ is any neighbourhood
of $c$ in $SU(2)$ containing $\bar{V}$ we have
$$ \int_{P^{-1}(V') \cap \mu^{-1}(0)/T} \Phi(\eta\a) 
= \int_{P^{-1}(V')\cap\zloc} 
\res_{X=0} \frac{\zeta \eta\alpha }{X-d\zeta} $$
$$ =  \int_{M(c)\cap\zloc} 
\res_{X=0} \frac{\zeta \eta}{X-d\zeta} = 
 \int_{M(c) \cap 
\mu^{-1}(0)/T} \Phi(\eta)$$
for any $T$-equivariantly closed differential form $\eta 
\in \O^*_T(P^{-1}(V'))$ such that
the intersection of $P^{-1}(\bar{V})$ with 
the support of $\eta$ is compact. Here
we have used the same notation for $\eta$ and its restriction to $M(c)$.

When $n>2$, so that the maximal torus $T$ of $K=SU(n)$
has dimension higher than one, then $\Ph(\eta)$
and $\int_{M(c)\cap\zloc/T} \Ph(\eta)$ are given by similar
formulas involving $n-1$ iterated residues (see \cite{GK}). In 
particular the support of $\Ph(\eta)$
is contained in the image of the support of $\eta$, and
$$ \int_{P^{-1}(V') \cap \mu^{-1}(0)/T} \Phi(\eta\a) 
= \int_{M(c) \cap 
\mu^{-1}(0)/T} \Phi(\eta)$$
for any $T$-equivariantly closed differential form $\eta 
\in \O^*_T(P^{-1}(V'))$ such that
the intersection of $P^{-1}(\bar{V})$ with 
the support of $\eta$ is compact.

\end{rem}

  \renorm
 \section{Nonabelian localization applied to extended moduli spaces}

Na\"{\i}ve application of the residue formula (Theorem 
\ref{t4.1}) to the
extended moduli space $M(c)$, using (\ref{9}) and Remark
\ref{r:fibprod}
 and
ignoring the fact that $M(c)$ is noncompact and has
singularities,
yields
\beq \label{1.7}
\prod_{r = 2}^n a_r^{m_r} \exp ({f_2}) [\mnd] = n C_K
\res 
 \Biggl (  \nusym^2(X) (\int_F e^{\omega})  \sum_{\d \in \intlat} 
\frac{ \prod_{r = 2}^n \tau_r(-X)^{m_r}e^{(\tc - \d) (X)}  }
{((-1)^{n_+} \nusym^{2}(X))^g} \Biggr ) \eeq
where the constant $C_K$ is defined at (3.7).
The main problem with (\ref{1.7}) (related to the noncompactness of $\xc$,
which permits the fixed point set $M(c)^T$ to be the union of infinitely 
many components $F_\d$) is that the sum over $\d$
does not converge for $\xvec \in \liet$. In this section we shall 
instead apply the version of
nonabelian localization due to
Guillemin-Kalkman and Martin (Propositions 3.6 and 3.9) to $M(c)$, using
Remarks \ref{leg} and \ref{arm};
this will lead to a proof 
that (\ref{1.7}) is true if interpreted appropriately (see Remark \ref{naive}).
First we use Proposition \ref{p:sm}.

\begin{lemma} \label{l3} Let $|W|=n!$ be the order of the Weyl group $W$
of $K=SU(n)$, and let $c=\diag(e^{2\pi i d/n},\ldots,e^{2\pi i d/n})$ where
$d$ is coprime to $n$.
If $V$ is a sufficiently small neighbourhood of $c$ in 
$K$ that the quotient $T/Z_n$ of $T$ by the centre $Z_n$ of $K=SU(n)$ acts  
freely on $P^{-1}(V) \cap \zloc$ (see Lemma \ref{l:5.10}),
then for any
$\eta\in\hk(X)$ we have 
$$ \int_{\mnd} \Ph(\eta e^{ \bom} ) = 
 \frac{1}{|W|} \int_{N(c)} 
\Ph(\nusym\eta e^{ \bom} ) =
\frac{1}{|W|}\int_{N(V)} \Ph(\nusym \eta e^{ \bom}\alpha)$$
where 
$$N(c)=\xc \cap \mu^{-1}(0)/T $$
for $\mu:K^{2g} \times \liek \to \liek$ given by minus the projection onto
$\liek$ and 
$$N(V) =  P^{-1}(V)\cap \zloc /T.$$
Also $\alpha$ is a
$T$-equivariantly closed form on $K^{2g}\times \liek$ 
representing the
$T$-equivariant Poincar\'{e} dual to $M(c)$, which is  chosen 
as in Proposition \ref{defa} so that
the  support of $\alpha$ is contained in $P^{-1}(V)$
and has compact intersection with $\zloc$.
\end{lemma}
\Proof Since $\mnd = \xc \cap \zloc/K$, we can first identify
$\int_{\mnd}\Phi(\eta e^{ \bom} ) $ with 
$$ \frac{1}{|W|} \int_{N(c)} 
\Ph(\nusym\eta e^{ \bom} ) $$
via Proposition \ref{p:sm}, whose proof works
in this situation even though $M(c)$ is noncompact and singular,
because $\mu$ is proper and $M(c)$ is nonsingular in a 
neighbourhood of $\mu^{-1}(0)$ (see Remark \ref{leg}). Then
we use Remark \ref{aa}. \hfill $\square$

Next we need to summarize some conventions on the roots
  and weights of $SU(n)$. 
The simple roots $\{e_j: j = 1, \dots, n-1\} $ of 
$SU(n)$ are elements of $\liets$; 
in terms of the standard identification of $\liet$ 
with $\{ (X_1, \dots, X_n ) \in \RR^n: \sum_i X_i = 0 \}$ under which
$(X_1, \dots, X_n ) \in \RR^n$ satisfying $ \sum_i X_i = 0$
corresponds to $X={\rm diag}({2\pi i X_1},\ldots,{2\pi i X_n})\in \liet$,
they are given by
\beq \label{6.1} e_j(X) = X_j - X_{j+1}. \eeq
The dual basis to the basis of simple roots 
(with respect to the 
inner product  $<\cdot,\cdot>$ defined at (\ref{1.02}) above,
which is the usual
Euclidean inner product 
on $\RR^n$) is the set of {\em fundamental weights}
$w_j \in \liets$ given by 
\beq \label{5.fundwts} w_j(X) = X_1 + \dots + X_j. \eeq
If we use this same inner
  product to identify $\liets$ with $\liet$, the simple roots
become identified with a set of generators 
$$\he{j}=(0,\ldots,0,1,-1,0,\ldots,0)$$ 
for the integer lattice $\intlat$ of $\liet$, and the fundamental weights
correspond to  
elements $\hw{j} \in \liet$ given by
$$\hw{j} = (1, \dots, 1, 0, \dots, 0 ) - \frac{j}{n}(1, \dots, 1). $$
In particular we have 
\beq \hat{w}_{n-1} = \frac{1}{n}(1, \dots, 1, - (n-1) ).  \eeq

Since we shall later apply induction on $n$,
it will be convenient to label certain spaces, groups and Lie algebras 
by the associated value of $n$. In particular the space $M(c )$
will sometimes be denoted by $M_n(c)$,
the maximal torus $T$ of $SU(n)$ 
by $\tor$, its Lie algebra $\liet$ 
by $\liet_n$, and the map $\Phi$ by $\Phi_n$.

We define a one dimensional torus 
$\tone \cong S^1$ in 
$SU(n) $  generated by $\he{1}$: it is identified 
with $S^1 $ via 
\beq t \in S^1 \mapsto ( t, t^{-1},1, \dots, 1 ) \in \tone. \eeq 
The (one dimensional) Lie algebra $\lietone$ is spanned
by $\he{1}$. Its orthocomplement in $\liet$ is 
\beq \lietc = \{ (X_1, \dots, X_n) \in \RR^n: X_{1} = X_2, 
\sum_{j = 1}^{n} X_j = 0 \}. \eeq
Define $\torc $ to be the torus given by 
$\exp (\lietc)$:
$$ \torc  = \{ (t_1,t_{1},t_3 \dots, t_{n-1}, t_{n} ) 
\in U(1)^{n} : \, \, (t_1)^2(\prod_{j = 3}^{n} t_j ) = 1 \}; $$
then $\torc$ is isomorphic to the maximal torus of
$SU(n-1)$ (i.e. $\torc \cong
(S^1)^{n-2} $) so this does not conflict with the notation already
adopted.

\begin{rem} The multiplication map \label{f6.4} $\tone \times \torc \to 
\tone \torc = \tor $ is a covering map with fibre 
$\tone \cap \torc = \ZZ_2 = \{(t,t^{-1},1, \dots, 1 ): \;\; t = t^{-1}\}. $
\end{rem}

There is  the following decomposition 
of the ring homomorphism
$ \phil$.
\begin{prop} \label{stages} For any symplectic manifold 
$M$ equipped with a Hamiltonian
action of $T_n$ such that $\torc$ acts locally freely on $\mu_{\torc}^{-1}(0)$,
the symplectic quotient $\mu_{\tor}^{-1}(0)/\tor$ may 
be identified with the symplectic
quotient of $\mu_{\torc}^{-1}(0)/\torc$ by 
the induced Hamiltonian action of $\tone$.
Moreover if in addition $T_n$ acts locally freely on $\mu_{\tor}^{-1}(0)$
then the ring homomorphism $\phil: H^*_\tor(M) \to H^*(
\mu_\tor^{-1}(0)/\tor)$ factors as 
$$ \phil = \phione \circ \phicom$$ where 
$$\phicom: H^*_{\tor}(M) \to 
H^*_{\tor}(\mu_{\torc}^{-1}(0)) \cong 
H^*_{\tone\times T_{n-1}}(\mu_{\torc}^{-1}(0)) \cong 
H^*_\tone (\mu_\torc^{-1}(0)/\torc) $$
and
$$ \phione: H^*_\tone (\mu_\torc^{-1}(0)/\torc ) 
\to H^*\Bigl((\mu_\torc^{-1}(0) \cap
 \mu_\tone^{-1}(0)/\torc\times \tone \Bigr ) 
\cong H^*(\mu_\tor^{-1}(0)/\tor). $$
\end{prop}
\Proof The isomorphisms 
$$ H^*_{\tor}(\mu_{\torc}^{-1}(0)) \cong 
H^*_{\tone\times T_{n-1}}(\mu_{\torc}^{-1}(0)) \cong 
H^*_\tone (\mu_\torc^{-1}(0)/\torc) $$
follow from Remark \ref{f6.4} and the fact that the cohomology with complex
coefficients of the classifying space of a finite group is trivial. 

Since $\mu_{\tor}$ is a $\tor$-invariant map, 
its projection $\mu_{\tone}$ onto $\lietone$
descends to $\mu_{\torc}^{-1}(0)/\torc$ and 
defines a moment map for the induced
$\tone$-action with respect to the induced symplectic 
structure on $\mu_{\torc}^{-1}(0)/\torc$.
The rest then follows from Remark \ref{f6.4} and naturality 
(cf. \cite{GK}, after (2.9)).\hfill $\square$

\begin{rem} \label{induct}
{}From now on, thanks to Lemma \ref{l3}, we shall be working with 
quotients by $T$ and subgroups of $T$, rather than quotients by $K$.
Because of this our arguments will apply to $M(c)$ when $c$ 
belongs to $T$ but is no longer necessarily a 
central element of $K$. This will be important
later, when we apply induction on $n$ using Proposition
\ref{stages}.  The only condition we will need to impose on $c\in T$
is that $c = \diag(c_1,\ldots,c_n)$ where the product of no
proper subsequence of $(c_1,\ldots,c_n)$ is 1; this is certainly
true for our original choice of $c$ when $c_j = e^{2\pi i d/n}$
for all $j$ with $d$ coprime to $n$. 

So for any $c\in T$, let us 
define
$$M(c) = M_n(c) = P^{-1}(c) = 
\Bigl \{   (h_1, \dots, h_{2g}, \L )  \in K^{2g} \times 
\liek \; : \; \prod_{j = 1}^g h_{2j-1} h_{2j} h_{2j-1}^{-1} h_{2j}^{-1} 
= c \exp (\L) \Bigr \} $$
where 
$P: K^{2g} \times \liek \to K$ is defined by 
$$  P \Bigl ( h_1, \dots, h_{2g}, \L \Bigr )  
= \prod_{j = 1}^g h_{2j-1} h_{2j} h_{2j-1}^{-1} h_{2j}^{-1} 
\exp (-\L). $$ 
Let us also define 
\beq  \nl(c)   = M(c ) \cap \mu^{-1}(0) / \tor 
\eeq
and
\beq  \nl(V)   = P^{-1}(V) \cap \mu^{-1}(0) / \tor \eeq
where $V$ is a small $T$-invariant neighbourhood of $c$ in $K$.
\end{rem}

\begin{prop} \label{p7.1} \label{f6.6} Suppose $c = \diag (c_1, 
\dots, c_n)\in T$  is 
such that the product of no proper subset of $(c_1,\ldots,c_n)$
is 1.  Then the group $\torc/\ZZ_n$, 
where $\ZZ_n$ 
consists of the
identity matrix multiplied by $n$th roots of unity, acts freely 
on $P^{-1}(V) \cap \mu_{\torc}^{-1} (0)$ for any sufficiently
small $T$-invariant neighbourhood $V$ of $c$ in $K$. 
Hence the quotient $P^{-1}(V) \cap \mu_{\torc}^{-1} (0)/\torc$
is smooth.
\end{prop}
\Proof  
The conjugation action of $(t_1, \dots, t_n) \in 
U(1)^n$ on the space of $n \times n$ matrices sends
$$ (A_{ij}) \mapsto (t_i t_j^{-1} A_{ij}). $$
Clearly $\ZZ_n$ acts trivially. Let us assume that 
$(h, \L) \in M(c) \cap \mu_{\torc}^{-1}(0) $ is fixed by the 
action of some element of $\torc$ which is not in $\ZZ_n$. After
rearranging the coordinates $X_3, \ldots, X_n$ if necessary,
we may assume that there is some $k$ between $3$ and $n$ such
that this element of $\torc$ is of the form
$(t_1, t_1, t_3, \ldots, t_{n-1}, t_{n})$ where $t_i = t_{1}$
if and only if $i\leq k$. Then each $h_j $ is block diagonal of the form
$$\matr{h_j^1}{0}{0}{h_j^2} 
$$
where $h_j^1$ is a $k\times k$ matrix and $h_j^2$ is $(n-k)\times(n-k)$.
As the determinant of any commutator is one, it follows that 
$\prod_{j=1}^n [h_{2j-1}, h_{2j}] $ is block diagonal of the form
\beq \label{7.0} \matr{A}{0}{0}{B} 
\eeq
where $\det A = \det B = 1$. But 
$\L$ is also block diagonal of the same form
$$\matr{\L_1}{0}{0}{\L_2} ,$$
and since $(h, \L) \in \mu_{\torc}^{-1}(0) $
 the diagonal entries of $\L$ are
$(2\pi i\l, -2\pi i\l,0, \dots, 0) $ for some $\l \in \RR$. 
Thus as $k\geq 3$ both $\L_1$
and $\L_2$ have trace $0$, so $\det\exp\L_1 = 1 = \det\exp\L_2$.
Since $(h, \L) \in M(c)$ it follows that the
matrix $A$ must equal
$$  \diag (c_{1}, \ldots, c_{k})\exp\L_1, $$ 
and hence
$$c_1\ldots c_k = \det A =1.$$ 
This contradiction to the hypotheses on $c$ shows that
$\torc/\ZZ_n$ acts freely on $M(c)\cap \mu_{\torc}^{-1}(0)$, and the same 
argument shows that $\torc/\ZZ_n$ acts freely
on $P^{-1}(V) \cap \mu_{\torc}^{-1}(\hat{V})$ for any
sufficiently small $T$-invariant neighbourhood $V$ of $c$ in $K$
and any sufficiently small neighbourhood $\hat{V}$ of 0 in $\liet_{n-1}$. 
The result follows. \hfill $\square$

\begin{definition}\label{YY}
Let us introduce coordinates
 $$Y_k = e_k(X) = \inpr{\he{k},X} $$ on 
$\liet$, corresponding to  the simple roots
$e_k \in \liets$. 
\end{definition}

We are now in a position to exploit Proposition \ref{p:gkm}
and Remark \ref{arm}, by using the
 translation map $s_{\tran} $ defined by Lemma
\ref{l4.3}, where $\tran = \he{1}$ lies in the integer
lattice  $\intlat$ and so satisfies $\exp(\L_0)=1$. 

\begin{lemma} \label{old5.17} \label{l6.4} Suppose $c = \diag (c_1, 
\dots, c_n)\in T$  is 
such that the product of no proper subset of $(c_1,\ldots,c_n)$
is 1. Suppose also that $\eta$ is a polynomial in 
the $\tar$ and $\tbr$, so that $s_{\he{1}}^*\eta  = \eta$. 
If $V$ is a sufficiently small $T$-invariant neighbourhood of $c$
in $K$ 
so that $P^{-1}(V) \cap \mu^{-1}(\lietone)/\torc$ is smooth
(see Proposition \ref{f6.6}),
and if $\nl(V) = P^{-1}(V) \cap \mu^{-1}(0)/\tor$ as before, then
$$\int_{\nl(V)} \phil( \eta e^{ \bom} e^{ - 
Y_1 } \alpha ) 
= \int_{P^{-1}(V) \cap \mu^{-1}(- \he{1} )/\tor}
\phil (\eta e^{ \bom}\alpha )  $$
$$ = \int_{\nl(V)} \phil (\eta e^{ \bom}\alpha ) -
n_0 \sum_{F\in\calf: -||\he{1}||^2  < \inpr{\he{1},\mu(F)} < 0}
\res_{Y_1 = 0 } \int_F 
\frac{ \phicom  ( 
 \eta e^{ \bom}\alpha)} {e_{{F}} } 
 $$
where $\calf$ is the set of components of the fixed point set of the 
action of $\tone$ on $P^{-1}(V)\cap \mu^{-1}(\lietone)/\torc$,
and $e_F$ denotes the $\tone$-equivariant Euler
class of the normal to $F$ in
$P^{-1}(V)\cap \mu^{-1}(\lietone)/\torc$
for any $F\in\calf$, 
while  $n_0$  is the order of the subgroup
of $\tone/\tone\cap\torc$ that acts trivially on $P^{-1}(V) \cap
\mu^{-1}(\lietone)/\torc$. 
 Also $\alpha$ is the $T$-equivariantly 
closed differential form on $K^{2g}\times \liek$ given by Proposition
\ref{defa} which
represents the equivariant Poincar\'{e} dual of $M(c)$, 
chosen so that the support of
$\alpha$ is contained in $P^{-1}(V)$.
\end{lemma}
\Proof  Since $\mu^{-1}(\lietone)=K^{2g}\times\lietone$ is contained
in $\mu_{\torc}^{-1}(0)$, it follows from 
Proposition \ref{f6.6} that if $V$ is a sufficiently
small $T$-invariant neighbourhood of $c$ in $K$, then $\torc/\ZZ_n$ acts
freely on $P^{-1}(V)\cap \mu^{-1}(\lietone)$ and so the quotient
$P^{-1}(V)\cap \mu^{-1}(\lietone)/\torc$ is smooth.

Since the 
restriction of $\mu_{\tone}$ to 
$\mu^{-1}(\lietone)$ is proper, and the support of $\alpha$ is contained
in $P^{-1}(V)$, by Remark \ref{arm} Guillemin and Kalkman's proof of 
Proposition \ref{p:gkm} can be 
applied to the $\tone$-invariant function induced by $\mu_{\tone}$ 
on the smooth manifold
$P^{-1}(V)\cap \mu^{-1}(\lietone)/\torc$
and the $\tone$-equivariant form induced by $\eta e^{\bom}\a$. 
In fact since $\tone\cap\torc \cong \ZZ_2$ acts trivially we can
work  with the action of $\tone/\tone\cap\torc$ instead of the action of 
$\tone$ (the Lie algebra and moment map are 
of course the same). This fits better with the choice
of coordinates $Y_k$ defined by the simple roots $\hat{e}_k$ because
the simple root $\hat{e}_{1}$ takes $(t,t^{-1},1,\ldots,1)\in\tone$
to $t^2$ and thus induces an isomorphism from $\tone/\tone\cap\torc$
to $S^1$. 
By combining this with Proposition \ref{stages} 
we get
$$ \int_{P^{-1}(V) \cap \mu^{-1} (0 )/\tor
}  \phil ( \eta e^{ \bom} \a) 
- \int_{P^{-1}(V) \cap \mu^{-1} (- \he{1} )/\tor }
\phil ( \eta e^{ \bom}\a)    $$  
$$ ~~~= n_0  
\res_{Y_1 = 0} \sum_{F\in\calf: -||\he{1}||^2  < \inpr{\he{1},\mu(F)} < 0}
\int_{{F}} 
\frac{\phicom ( \eta e^{ \bom}\a)  }{e_{{F}} }. $$
Now note that 
the restriction of $P:K^{2g}\times\liek \to K$
to $\mu^{-1}(\liet) = K^{2g}\times \liet$ is invariant under the translation
$s_{\L_0}$ for $\L_0 \in \L^I$. 
Therefore by construction the restriction of $\alpha$
to $\mu^{-1}(\liet)$ is also invariant under this translation. Thus
by  (\ref{8.4}) and Definition \ref{YY}
$$\int_{P^{-1}(V) \cap \mu^{-1} (-\he{1})/\tor} 
 \phil (\eta e^{ \bom} \alpha ) 
= \int_{\nl(V)} \phil \Bigl ( s_{\he{1}}^* (\eta e^{ \bom}\alpha ) 
\Bigr )  = \int_{\nl(V)}  \phil (\eta e^{ \bom} 
e^{ - Y_1}\alpha ). $$
The result follows. \hfill $\square$

\begin{rem}
It will follow from the proof of Proposition \ref{p7.2} below that $n_0=1$
here (see Remark \ref{n0=2}).
\end{rem}

\renorm
\section{Fixed point sets of the circle action}

In this section we shall consider the components $F\in\calf$ of the fixed 
point set of the action of $\tone$ on the quotient 
$P^{-1}(V) \cap\mu^{-1}(\lietone)/\torc$ (which appeared in Lemma 
\ref{old5.17}). Since $P^{-1}(c)=M(c)$ and
$V$ is an arbitrarily small $T$-invariant 
neighbourhood of $c$ in $K$, we may assume that
every $F\in\calf$ contains a component of the
fixed point set of the action of $\tone$ on $M(c) 
\cap\mu^{-1}(\lietone)/\torc$, and each of these components is contained
in a unique $F\in\calf$. So we shall start by analysing the
components of the
fixed point set of the action of $\tone$ on $M(c) 
\cap\mu^{-1}(\lietone)/\torc$.   
We shall find that they can be described
inductively in terms of products of spaces of the form $N(c)$
(see Remark \ref{induct}) for smaller values
of $n$. This will enable us to use induction in the next two sections 
to express
the intersection pairings $\int_{\mnd} \Ph(\eta e^{\bom})$
on the moduli spaces $\mnd$ as iterated residues (see  Theorem
\ref{mainab}
 and 
Theorem \ref{t9.6}).

\begin{prop} \label{p7.2} Suppose that $c=\diag(c_1,\ldots,c_n)\in SU(n)$
is such that the product of no proper subsequence of $(c_1,\ldots,c_n)$
is 1.
Then the components of the fixed  \label{fil}
point set of the action of $\tone$ on the quotient $(M(c)
\cap\mu^{-1}(\lietone))/\torc$
may be described as follows. For any subset  $I$ of $\{ 3, \dots, n\}$ 
let $I_1=I\cup\{1\}$ 
and let $I_2= \{ 1, \dots, n\} - I_1$.
Let $H_I$ be the subgroup of $SU(n)$ given by
$$H_I = \{ (a_{ij}) \in SU(n): 
a_{ij} = 0 ~\mbox{if  $(i,j) \in (I_1\times I_2)\cup (I_2\times I_1)$ } \}. $$
Suppose that  $\l \in \RR$ is a solution of 
$$ e^{-2\pi i \l} = c_{i_1} \dots c_{i_{r}} = \prod_{j\in I_1} c_j$$
where $r$ is the number of elements of $I_1=\{i_1,\ldots,i_r\}$, so that
$$ e^{2\pi i \l} =  \prod_{j\in I_2} c_j.$$
Then we have a component of the fixed point set given by
$\fisml = \tilde{F}_{I,\l}/\torc$ where 
$$\tilde{F}_{I,\l} = M(c) \cap (H_I^{2g} \times \{\l \hat{e}_1 \}),$$
and every component is of this form for some 
subset $I$ of $\{3,\ldots, n\}$
and solution $\l$ to the equation above. 
\end{prop}
\Proof  Suppose the $\tone$ orbit of a 
point $ (h_1, \dots, h_{2g}, \L) 
\in SU(n)^{2g} \times \lietone$  is contained in its
orbit under $\torc$. A general element of the $\tone$ orbit 
of an $n\times n$ matrix  $ A =  (a_{ij})$ under conjugation looks like 
$$
\left \lbrack \begin{array}{lcccr}
a_{11} & t^2a_{12} & ta_{13} & \dots &  t a_{1n} \\
t^{-2}a_{21} & a_{22} & t^{-1}a_{23} & \dots & t^{-1} a_{2n} \\
t^{-1}a_{31} & ta_{32} & a_{33} & \dots & a_{3n} \\
\vdots & \vdots & \vdots & \dots           & \vdots \\
t^{-1} a_{n1}& ta_{n2} & a_{n3} & \dots  & a_{nn} 
\end{array}\right \rbrack $$
while a general element of the $\torc$ orbit of $A$ looks like 
$$ \left \lbrack \begin{array}{lcccr}
a_{11} & a_{12} & t_1 t_3^{-1}a_{13} & \dots & t_1 t_n^{-1} a_{1n} \\
a_{21} & a_{22} & t_1 t_3^{-1}a_{23} & \dots & t_1 t_n^{-1} a_{2n} \\
t_3 t_1^{-1} a_{31} & t_3 t_1^{-1}a_{32} & a_{33} & \dots & t_3 t_n^{-1} 
a_{3n}    \\
\vdots & \vdots & \vdots & \dots           & \vdots \\
t_n t_1^{-1} a_{n1}&t_n t_1^{-1} a_{n2} & t_n t_3^{-1} a_{n3} &\dots 
& a_{nn}
\end{array}\right \rbrack . $$
For each $t$ there exist $t_1,t_3,\ldots,t_n$ such that these two matrices
are equal when $A$ is any of $h_1,\ldots, h_{2g}$ and $\L$. 
Choose $t\neq t^{-1}$ and
let $I$ denote the set of $j$ in $\{ 3, \dots, n\}$ for which 
$t_1 t_j^{-1}
= t$. Similarly, define $J$ to be the set 
of $j$ in $\{ 3, \dots, n\}$ for which
$t_1 t_j^{-1}
= t^{-1}$, and let $K = \{ 3, \dots, n\} - I - J.$ 
Reordering the coordinates one finds that 
all the $h_j$ and $\L$ are  block diagonal
where the blocks correspond to $ I \cup \{1\}$, 
$ J \cup \{2\} $ and $K$. 
 Conversely, if all the $h_j$ 
are block diagonal of this form and $\L\in\lietone$,
then the $\tone$ orbit of $(h_1,\ldots,h_{2g},\L)$  is contained
in its $T_{n-1}$ orbit since given any $t\in U(1)$ we can find $(t_1, 
t_1, t_3, \ldots,t_n)$ in $T_{n-1}$ satisfying
$t_1 t_j^{-1}
= t$ if $j\in I$ and
$t_1 t_j^{-1}
= t^{-1}$ if $j \in J.$ 

We next prove that $K$ is empty. Suppose otherwise; then 
as the determinant of any commutator is one, 
$\det \prod_{j=1}^g [h_{2j-1}^{[K]}, h_{2j}^{[K]}] =1 $ 
(where the superscript $[K]$ denotes the block
of the matrix  corresponding to $K$). Thus the $K$ block in 
$c $ also has determinant $1$. This is impossible by
the hypothesis on $c$.

Suppose now that $(h_1, \dots, h_{2g}, \L) $ $ \in M(c)\cap \mu^{-
1}(\lietone)$
lies in $H_I^{2g} \times \lietone.$ 
Then 
$$\L = \l \he{1} = 2 \pi i {\rm diag}(\l,-\l,0,\ldots,0)$$
for some $\l \in \RR$, so the blocks $\L^{[I_1]}$ and $\L^{[I_2]}$ of $\L$
corresponding to $I_1 = I \cup \{1\}$ and 
$I_2 = J \cup \{2\}$ 
satisfy $\det\exp\L^{[I_1]} = e^{2 \pi i \l}$ 
and $\det\exp\L^{[I_2]}= e^{-2 \pi i \l}$. 
But
$$\det (\prod_{j=1}^g [h_{2j-1}^{[I_1]}, h_{2j}^{[I_1]}])=1=
\det (\prod_{j=1}^g [h_{2j-1}^{[I_2]}, h_{2j}^{[I_2]}])  $$
because the determinant of any commutator is one.
It therefore follows from the definition of $M(c)$ 
 that 
$$ e^{-2\pi i \l} = \prod_{j\in I_1} c_j.$$
This is enough to complete the proof. $\square$

\begin{rem} \label{n0=2} The proof of this proposition shows
that the elements of $\tone$ which act trivially on the quotient $(M(c)
\cap\mu^{-1}(\lietone))/\torc$ are
 precisely those represented by $t$ satisfying
$t=t^{-1}$, i.e. $t=\pm 1$, or equivalently those in 
$\tone\cap \torc$. Thus the size $n_0$ of the subgroup
of $\tone/\tone\cap\torc$ acting trivially on the quotient $M(c)
\cap\mu^{-1}(\lietone))/\torc$ is 1 (cf. Lemma \ref{old5.17}).
\end{rem}

\begin{prop} \label{p7.3}   Suppose that $c=\diag(c_1,\ldots,c_n)\in SU(n)$
is such that the product of no proper subsequence of $(c_1,\ldots,c_n)$
is 1. Suppose that $I$ is a subset of $\{3,\ldots,n\}$
with $r-1$ elements where $1\leq r\leq n-1$, and let 
$I_1 = I \cup \{1\} = \{i_1,\ldots,i_r\}$
and $I_2 =  \{1,\ldots,n \}- I_1 = \{i_{r+1},\ldots, i_n\}$. 
Suppose also that  $\l \in \RR$ is a solution of 
$$ e^{-2\pi i \l} = \prod_{j\in I_1} c_j$$
so that
$ e^{2\pi i \l} =  \prod_{j\in I_2} c_j.$
Let 
$$ c(I_1,\lambda) = \diag(c^{I,\lambda}_{i_1},\ldots,
c^{I,\lambda}_{i_r})$$
and
$$ c(I_2, -\lambda) = 
\diag(c^{I,\lambda}_{i_{r+1}}, \ldots, c^{I,\lambda}_{i_n})$$
where $ c^{I,\lambda}_j = c_j$
if $j \geq 3$, while
$ c^{I,\lambda}_{1} = c_{1} e^{2\pi i \l}$
and
$c^{I,\lambda}_2 = c_2 e^{-2\pi i \l}.$
Let $\fisml$ be defined as in Proposition \ref{p7.2}. 
Then there is a finite to one
(in fact $(r(n-r))^{2g}$ to one) surjective smooth map 
$$ \Psi_{I,\lambda}: (S^1)^{2g} \times N_r (c(I_1,\lambda)) 
\times N_{n-r} (c(I_2,-\lambda)) 
\to \fisml. $$ 
\end{prop} 
\Proof We define a homomorphism 
$$\r_I: S^1 \times SU(r) \times SU(n-r) \to 
H_I \subset SU(n)$$
given by 
$$\r_I : (s,A,B)  \mapsto \matr{s^{n-r} A}{0}{0}{s^{-r} B} $$ 
with respect to the decomposition of $\{1,\ldots,n\}$ as $I_1\cup I_2$.
Note that $\r_I$ restricts to an $r(n-r)$ to one surjective
homomorphism
$$\r_I:S^1 \times T_r \times T_{n-r} \to \tor.$$
               
If $\epsr{r}: SU(r)^{2g} \to SU(r)$ is defined by 
$  \epsr{r} (h_1, \dots, h_{2g}) =
\prod_{j=1}^g [h_{2j-1}, h_{2j}] $ then 
\beq \label{7.4}  \epsr{n}\Bigl (\r_I(s_1, A_1, B_1), \dots, 
\r_I(s_{2g}, A_{2g}, B_{2g}) \Bigr ) 
= \matr{\epsr{r} (A_1, \dots, A_{2g}) } {0}{0}{ \epsr{n-r}(B_1, \dots, 
B_{2g} ) }. \eeq
Let us define a map 
$$ \Psi_{I,\lambda}: (S^1)^{2g} \times N_r (c(I_1,\lambda)) 
\times N_{n-r} (c(I_2,-\lambda)) 
\to \fisml$$
as the quotient of 
$$\tpsi_{I,\lambda}: (S^1)^{2g} \times \Bigl (\mu_{SU(r)}^{-1}(0) \cap
M_r (c(I_1,\lambda))   \Bigr )
\times \Bigl (\mu_{SU(n-r)}^{-1}(0)  \cap 
M_{n-r} (c(I_2,-\lambda)) \Bigr )
\to \tilde{F}_{I,\l} $$ 
defined by $$\tpsi_{I,\lambda} \Bigl ( (s_1, \dots, s_{2g}), 
(h_1^{[I_1]}, \dots, h_{2g}^{[I_1]}, 0), 
(h_1^{[I_2]}, \dots, h_{2g}^{[I_2]}, 0) \Bigr ) $$ $$
= \Bigl ( \r_I (s_1, h_{1}^{[I_1]}, h_{1}^{[I_2]}),\r_I (s_{2}, 
h_{2}^{[I_1]}, h_{2}^{[I_2]} ),
\dots, \r_I (s_{2g},h_{2g}^{[I_1]}, h_{2g}^{[I_2]} ), 
2\pi i \diag(\l, -\l,0,\ldots, 0) \Bigr ). $$
Here, $\tilde{F}_{I,\l}$ was defined in Proposition \ref{p7.2}. 

We must check that the image of $\tpsi_{I,\lambda}$ is 
contained in $\tilde{F}_{I,\l}$.
We have 
$$\epsr{r} (h_1^{[I_1]}, \dots, h_{2g}^{[I_1]})
= \diag(c^{I,\lambda}_{i_1},\ldots,
c^{I,\lambda}_{i_r}) = 
c(I_1,\lambda)$$
and
$$ \epsr{n-r} (h_1^{[I_2]}, \dots, h_{2g}^{[I_2]})
=\diag(c^{I,\lambda}_{i_{r+1}}, \ldots, c^{I,\lambda}_{i_n})
= c(I_2,-\lambda).  $$
In order to show that
$\tpsi_{I,\lambda} \Bigl ( (s_1, \dots, s_{2g}), 
(h_1^{[I_1]}, \dots, h_{2g}^{[I_1]}, 0), 
(h_1^{[I_2]}, \dots, h_{2g}^{[I_2]}, 0) \Bigr ) $ lies 
in $\tilde{F}_{I,\l}$
we need to check that if $\L = 2\pi i \diag(\l,-\l,0,\ldots,0)$ then
$ c \exp ( \L)$ is block diagonal of the form
$$ \matr{c(I_1,\lambda) }{0}{0}
{c(I_2,-\lambda)} $$
with respect to the decomposition of $\{1,\ldots,n\}$ as $I_1\cup I_2$.
This follows by the choice of $c^{I,\lambda}_1,\ldots,c^{I,\lambda}_n$.

We must also check that the map $\Psi_{I,\lambda}$ is well defined on the 
quotient by the action of $T_{r}\times T_{n-r}$: in other words we must 
check that for any  $t=(t_1, \dots,t_n)$ $\in U(1)^n$ 
satisfying $t_{i_1}, \dots, t_{i_r} = t_{i_{r+1}}, \dots, t_{i_n} = 1 $ so that
$t_{I_1}=(t_{i_1} \dots t_{i_r})\in T_r$ 
and $t_{I_2} =( t_{i_{r+1}} \dots t_{i_n} ) \in T_{n-r}$, we have
$$\tpsi_{I,\lambda} \Bigl ( (s_1, \dots, s_{2g}), 
(t_{I_1}h_1^{[I_1]}(t_{I_1})^{-1}, \dots, t_{I_1}h_{2g}^{[I_1]}
(t_{I_1})^{-1}, 0), 
(t_{I_2}h_1^{[I_2]}(t_{I_2})^{-1}, \dots, t_{I_2}h_{2g}^{[I_2]}
(t_{I_2})^{-1}, 0) \Bigr ) $$ 
\beq \label{7.009} 
 = \tilde{t} \tpsi_{I,\lambda} \Bigl ( (s_1, \dots, s_{2g}),
(h_1^{[I_1]}, \dots, h_{2g}^{[I_1]}, 0),
(h_1^{[I_2]}, \dots, h_{2g}^{[I_2]}, 0) \Bigr )  \eeq
for some 
$\tilde{t} = (\tilde{t}_1, \dots, \tilde{t}_n) $ satisfying
$\prod_{j = 1}^n \tilde{t}_j = 1$ and 
$\tilde{t}_{1} = \tilde{t}_2$. For any 
$s \in U(1)$, we may conjugate  all the $h_j^{[I_1]}$ by
$s^{n-r}$ and all the $h_j^{[I_2]}$ by $s^{-r} $ without
changing the image under $\tpsi_{I,\lambda} $;
choosing $s$ so that 
$s^{n-r} \tilde{t}_{1} = s^{-r} \tilde{t}_2$ we find that the equation 
(\ref{7.009}) is satisfied for $\tilde{t}_j = t_j s^{n-r} $ 
(when $j \in I_1$) and $\tilde{t}_j = t_{j} s^{-r} $ (when 
$j \in I_2$). 

To show that $\Psi_{I,\lambda}$ is finite-to-one and surjective, suppose that
$(h_1,\ldots,h_{2g},\L) \in \tilde{F}_{I,\l} $; we must check that
a finite (and nonzero) number of $T_r \times T_{n-r}$ orbits in
$$(S^1)^{2g} \times \Bigl (\mu_{SU(r)}^{-1}(0) \cap
M_r (c(I_1,\lambda))   \Bigr )
\times \Bigl (\mu_{SU(n-r)}^{-1}(0)  \cap 
M_{n-r} (c(I_2,-\lambda)) \Bigr )$$
map into the $\torc$ orbit of $(h_1,\ldots,h_{2g},\L)$. Now by the definition
of $\tilde{F}_{I,\l}$ we have
$\L = 2 \pi i \diag (\l,-\l,0,\ldots,0)$ 
and each $h_j$ is block diagonal of the form
$$\matr{\ha{j}}{0}{0}{\hb{j}} $$
with respect to the decomposition of $\{1,\ldots,n\}$ as $I_1\cup I_2$. So
$$\tpsi_{I,\lambda} \Bigl ( (s_1, \dots, s_{2g}), 
(H_1^{[I_1]}, \dots, H_{2g}^{[I_1]}, 0), 
(H_1^{[I_2]}, \dots, H_{2g}^{[I_2]}, 0) \Bigr ) $$
belongs to the $\torc$ orbit of  $(h_1,\ldots,h_{2g},\L)$ if 
 and only if there is some
$\tilde{t} = (\tilde{t}_1, \dots, \tilde{t}_n) $ satisfying
$\prod_{j = 1}^n \tilde{t}_j = 1$ and 
$\tilde{t}_{1} = \tilde{t}_2$ such that
$$s_j^{n-r} H_j^{[I_1]} = \tilde{t}_{I_1} h_j^{[I_1]}(\tilde{t}_{I_1})^{-1}$$
and
$$s_j^{-r} H_j^{[I_2]} = \tilde{t}_{I_2} h_j^{[I_2]}(\tilde{t}
_{I_2})^{-1}$$
where $\tilde{t}_{I_1}=(\tilde{t}_{i_1}, \dots ,\tilde{t}_{i_r})\in T_r$ 
and $\tilde{t}_{I_2} =( \tilde{t}_{i_{r+1}}, \dots, \tilde{t}_{i_n} ) 
\in T_{n-r}$.
Since $\det  H_j^{[I_1]} =1 = \det H_j^{[I_2]} $ and
$\det h_j^{[I_1]}  \det h_j^{[I_2]} =1$, by the argument of the 
previous paragraph
this happens if and only if $(s_j)^{r(n-r)} = \det h_j^{[I_1]}$ and
$s_j^{n-r} H_j^{[I_1]}$ is conjugate to $ h_j^{[I_1]}$
and
$s_j^{-r} H_j^{[I_2]}$ is 
conjugate to $h_j^{[I_2]}$. Thus $\Psi_{I,\lambda}$ is surjective
and $(r(n-r))^{2g}$ to one.
$\square$ 

\begin{rem}  Note that by the definition of $c_j^{I,\l}$
(see Proposition \ref{p7.3}) no proper subsequence of
$(c^{I,\lambda}_{i_1},\ldots,c^{I,\lambda}_{i_r})$ or 
$(c^{I,\lambda}_{i_{r+1}},\ldots,c^{I,\lambda}_{i_n})$ 
has product equal to 1, because the same 
is true of $(c_1,\ldots,c_n)$.
\end{rem}

\begin{rem} \label{big} 
It follows from Proposition \ref{p7.3} that if 
$\Psi_{I,\lambda}$ is orientation preserving (and we shall see
below in Remark \ref{indbasrem}
 that $\Psi_{I,\lambda}$ takes a natural symplectic
orientation on $(S^1)^{2g} \times N_r(c(I_1,\lambda))
 \times N_{n-r}(c(I_2,-\lambda))$ to the symplectic
orientation induced by $\omega$ on $F_{I,\l}$) then
\beq \label{snow} \int_{F_{I,\l}} \phicom (\eta e^{ \bom}) = (r(n-r))^{-2g}
\int_{(S^1)^{2g} \times N_r(c(I_1,\lambda))
 \times N_{n-r}(c(I_2,-\lambda))
} \Psi_{I,\lambda}^*  \phicom (\eta e^{ \bom})  \eeq
where both sides are elements of $H^*_{\tone}$. To be more precise we should
replace $\eta e^{ \bom}$ on each side of this equation by its restriction to 
$\tilde{F}_{I,\l}$, and as in Proposition \ref{stages} and Lemma \ref{l6.4}
we use the double cover $\tone \times T_{n-1} \to \tor$ to define
$$\phicom:H^*_{\tor}(\tilde{F}_{I,\l}) \to H^*_{\tone}(\fisml)
\cong H^*_{\tone} \otimes H^*(\fisml).$$
Recall from the proof of the last proposition that the homomorphism 
$$\r_I: S^1 \times SU(r) \times SU(n-r) \to 
 H_I \subset SU(n)$$
given by 
$$\r_I : (s,A,B)  \mapsto \matr{s^{n-r} A}{0}{0}{s^{-r} B} $$ 
with respect to the decomposition of $\{1,\ldots,n\}$ as $I_1\cup I_2$
restricts to an $r(n-r)$ to one surjective
homomorphism
$$\r_I:S^1 \times T_r \times T_{n-r} \to \tor.$$
It is easy to check that the inclusions of $\tone$, $T_r$ and $T_{n-r}$
in $T_n$ induce an isomorphism
$$\hat{\r}_r:\tone \times T_r \times T_{n-r} \to \tor$$
such that $\r_I$ and $\hat{\r}_r$ have the same restriction to
$T_r \times T_{n-r}$. The composition of this restriction
with the natural surjection from $\tor$ to $\tor/\tone\cong T_{n-1}
/(\tone \cap T_{n-1})$ gives an isomorphism
$$T_r \times T_{n-r} \to \tor/\tone\cong T_{n-1}
/(\tone \cap T_{n-1}).$$
Moreover the composition of $\r_I$ with the inverse of
$\hat{\r}_r$ defines a finite (in fact $r(n-r)$ to one) cover
$$\nu_I:S^1 \times T_r \times T_{n-r} \to
\tone \times T_r \times T_{n-r}$$ 
which restricts to the identity on $T_r \times T_{n-r}$
and induces a finite cover $\nu_I:S^1 \to \tone$
and isomorphisms on Lie algebras and equivariant
cohomology.

The argument in the proof of the last proposition to 
show that the map $\Psi_{I,\lambda}$ is well defined
on the quotient by the action of $T_r \times T_{n-r}$ may 
be rephrased as the statement that the map
$$\tpsi_{I,\lambda}: (S^1)^{2g} \times \Bigl (\mu_{SU(r)}^{-1}(0) \cap
M_r (c(I_1,\lambda))   \Bigr )
\times \Bigl (\mu_{SU(n-r)}^{-1}(0)  \cap 
M_{n-r} (c(I_2,-\lambda)) \Bigr )
\to \tilde{F}_{I,\l} $$ 
defined in the proof of Proposition \ref{p7.3} satisfies
$$\tpsi_{I,\lambda} \Bigl (t ((s_1, \dots, s_{2g}), 
(h_1^{[I_1]}, \dots, h_{2g}^{[I_1]}, 0), 
(h_1^{[I_2]}, \dots, h_{2g}^{[I_2]}, 0)) \Bigr )  $$
$$ = \r_I(t) \tpsi_{I,\lambda} \Bigl ( (s_1, \dots, s_{2g}),
(h_1^{[I_1]}, \dots, h_{2g}^{[I_1]}, 0),
(h_1^{[I_2]}, \dots, h_{2g}^{[I_2]}, 0) \Bigr ) $$
for all $t \in S^1 \times T_r \times T_{n-r}$, where $S^1$
acts trivially. Thus $\tpsi_{I,\lambda}$ and $\r_I$ induce
$\tpsi_{I,\lambda}^*$ from
$H^*_{\tor}(\tilde{F}_{I,\l}) $
to 
$$H^*_{S^1} \otimes H^*((S^1)^{2g}) \otimes 
H^*_{T_r} (\mu_{SU(r)}^{-1}(0) \cap
M_r (c(I_1,\lambda) ))
\otimes H^*_{T_{n-r}} (\mu_{SU(n-r)}^{-1}(0)  \cap 
M_{n-r} (c(I_2,-\lambda)))$$
and
$$ \nu_I^* \int_{(S^1)^{2g} \times N_r(c(I_1,\lambda))
 \times N_{n-r}(c(I_2,-\lambda))   
} \Psi_{I,\lambda}^*  \phicom (\eta e^{ \bom})  $$        
$$ \label{7.51} =  
 \int_{(S^1)^{2g} \times N_r(c(I_1,\lambda))
 \times N_{n-r}(c(I_2,-\lambda))
} (1\otimes\Phi_r \otimes \Phi_{n-r}) 
\tpsi_{I,\lambda}^*(\eta e^{\bom}), $$
where both sides are elements of $H^*_{S^1}$. Hence
by (\ref{snow}),
if $Y_1$ is the coordinate on $\lietone$ given by the
restriction of $Y_1=X_1-X_2$ on $\liet$ and $Y_1^I$ is the
coordinate on the Lie algebra of $S^1$ obtained from $Y_1$
via the isomorphism on Lie algebras induced by $\nu_I:S^1 \to \tone$, 
then
$$ \res_{Y_1 =0} \int_{F_{I,\l}} \phicom (\eta e^{ \bom}) = $$
$$(r(n-r))^{-2g}
\res_{Y_1^I=0}  \int_{(S^1)^{2g} \times N_r(c(I_1,\lambda))
 \times N_{n-r}(c(I_2,-\lambda))
} (1\otimes\Phi_r \otimes \Phi_{n-r}) 
\tpsi_{I,\lambda}^*(\eta e^{\bom}). $$

Since $S^1$ acts trivially, the residue operation
$\res_{Y^I_1=0}:H^*_{S^1} \to \CC$
can be extended to map
$$H^*_{S^1} \otimes H^*((S^1)^{2g}) \otimes 
H^*_{T_r} \Biggl( \mu_{SU(r)}^{-1}(0) \cap
M_r (c(I_1,\lambda) )\Biggr )
\otimes H^*_{T_{n-r}} \Biggl( \mu_{SU(n-r)}^{-1}(0)  \cap 
M_{n-r} (c(I_2,-\lambda))\Biggr) $$
to
$$H^*((S^1)^{2g}) \otimes 
H^*_{T_r} \Biggl( \mu_{SU(r)}^{-1}(0) \cap
M_r (c(I_1,\lambda) )\Biggr)
\otimes H^*_{T_{n-r}} \Biggl( \mu_{SU(n-r)}^{-1}(0)  \cap 
M_{n-r} (c(I_2,-\lambda))\Biggr) $$
so that it commutes with $\Ph_r$ and $\Ph_{n-r}$ and with
integration over $N_r(c(I_1,\lambda))$ and integration over
$ N_{n-r}(c(I_2,-\lambda))$. In particular by expressing
integrals over products as iterated integrals we obtain
$$\res_{Y_1 =0} \int_{F_{I,\l}} \phicom (\eta e^{ \bom})$$
$$ = (r(n-r))^{-2g}
\int_{N_r(c(I_1,\lambda)) } \Phi_r (\res_{Y_1^I=0} 
\int_{N_{n-r}(c(I_2,-\lambda))} \Phi_{n-r}( \int_{(S^1)^{2g} }
\tpsi_{I,\lambda}^*(\eta e^{\bom})))$$                     
$$=(r(n-r))^{-2g}
\int_{N_{n-r}(c(I_2,-\lambda)) } \Phi_{n-r} (\res_{Y_1^I=0} 
\int_{N_r(c(I_1,\lambda))} \Phi_{r}( \int_{(S^1)^{2g} }
\tpsi_{I,\lambda}^*(\eta e^{\bom}))).$$                     
This will be important when we apply induction later.
\end{rem}

Recall from Lemma \ref{old5.17} that $\calf$ is the set of components
of the fixed point set
 of the action of $\tone$ on the quotient 
$P^{-1}(V) \cap\mu^{-1}(\lietone)/\torc$, where
$V$ is a sufficiently small $T$-invariant neighbourhood of $c$ in $K$.
Every $F\in\calf$ contains a component $\fisml$ of the
fixed point set of the action of $\tone$ on $M(c) 
\cap\mu^{-1}(\lietone)/\torc$, and each $\fisml$ is contained
in a unique $F\in\calf$ (see Proposition \ref{p7.2}
for the definition of $\fisml$). 
For each $I$ and $\l$ we now need to understand the
normal bundle in $P^{-1}(V) \cap\mu^{-1}(\lietone)/\torc$
to the component $F\in \calf$ of the fixed point set which
contains $\fisml$.
First, we observe that there is  the following decomposition:

\begin{rem} \label{p5.11} \label{feb} Let $I$ be a subset of $\{3,\ldots,n\}$
with $r-1$ elements where $1\leq r\leq n-1$, 
let $I_1=I\cup\{1\}$ 
and let $I_2= \{ 1, \dots, n\} - I_2$. Then
$$\nusym_n(X) = \nusym^{[I_1]}_{r}(X) \nusym^{[I_2]}_{n-r}(X) \tau_I(X)  $$ 
where 
$$\nusym^{[I_1]}_{r}(X) 
 = \prod_{1\leq j<k\leq r } (X_{i_j} - X_{i_k}) $$
is the product of the positive roots of $SU(r)$ embedded in $SU(n)$
via the inclusion of $I_1$ in $\{1,\ldots,n\}$,
$$\nusym^{[I_2]}_{n-r}(X) 
= \prod_{r+1\leq j < k \leq n } (X_{i_j} - X_{i_k}) $$
is the product of the positive roots of $SU(n-r)$ embedded in $SU(n)$
via the inclusion of $I_2$ in $\{1,\ldots,n\}$,
and $$\tau_I(X) = 
\pm \prod_{1\leq j \leq r < k \leq n} ( X_{i_j} - X_{i_k}), $$
where the sign is $+$ or $-$ depending on whether the permutation
$$ \left \lbrack \begin{array}{lcccr}
1&2&\dots & n \\
i_1 & i_2 & \dots &i_n
\end{array} \right \rbrack $$
is even or odd. Note also that
$$(-1)^{r(n-r)} (\tau_I(X))^2 = \prod_{(i,j)\in I_1\times 
I_2\cup I_2\times I_1}
( X_i - X_j). $$
\end{rem}

Now we can find the $\tone$-equivariant Euler class of
the normal bundle in $P^{-1}(V) \cap\mu^{-1}(\lietone)/\torc$
to the component $F\in \calf$ of the fixed point set which
contains $\fisml$.

\begin{lemma} \label{p5.13} Let $I$ be a subset of $\{3,\ldots,n\}$
with $r-1$ elements where $1\leq r\leq n-1$, and let
$\l \in \RR$ be a solution of the equation
$$ e^{-2\pi i \l} =  \prod_{j\in I_1} c_j.$$
Then the $\tone$-equivariant Euler class of the normal bundle 
in $P^{-1}(V) \cap\mu^{-1}(\lietone)/\torc$
to the component $F\in \calf$ of the fixed point set 
 of the action of $\tone$ on  
$P^{-1}(V) \cap\mu^{-1}(\lietone)/\torc$ which
contains $\fisml$  is given by 
$ e_F = $ $(-1)^{r(n-r)g} \Ph_{n-1}(\tau_I^{2g})$.
\end{lemma}

\Proof The proof of Proposition \ref{p7.2} shows that
the component $F\in \calf$ of the fixed point set 
 of the action of $\tone$ on  
$P^{-1}(V) \cap\mu^{-1}(\lietone)/\torc$ which
contains $\fisml$  is
$$F = P^{-1}(V) \cap (H_I^{2g} \times \lietone)/\torc,$$
whereas $\mu^{-1}(\lietone)=K^{2g}\times\lietone$.
The  $T$-equivariant
Chern roots of the normal
bundle to $H_I^{2g}$ in $K^{2g}$ are $X_i-X_j$ for $(i,j) \in
I_1\times I_2 \cup I_2 \times I_1$ with multiplicity
$g$. The result follows by Remark \ref{p5.11}.
 \hfill $\square$

\begin{lemma} \label{l5.21}
Let $I$ be a subset of $\{3,\ldots,n\}$
with $r-1$ elements where $1\leq r\leq n-1$, 
let $I_1=I\cup \{1\}$, let $I_2 = \{1,\ldots,n\} - I_1$ and let
$\l \in \RR$ be a solution of the equation
$$ e^{-2\pi i \l} =  \prod_{j\in I_1} c_j.$$
Let $F$ be the component of the fixed point set 
 of the action of $\tone$ on  
$P^{-1}(V) \cap\mu^{-1}(\lietone)/\torc$ which
contains $\fisml$, where $\fisml$ is as defined in Proposition
\ref{p7.2}.
 We then have 
$$  \int_F \frac{\phicom ( \eta e^{ \bom} \a)} {e_F } =  
(-1)^{r(n-r)(g-1)}
\int_{F_{I,\lambda}} \phicom (\frac{ \eta e^{ \bom} }{\tau_I^{2g-1}}) $$
$$= (-1)^{r(n-r)(g-1)}
\int_{F_{I,\lambda}} \phicom (\frac{\nusym^{[I_1]}_{r}(X) 
\nusym^{[I_2]}_{n-r}(X) \eta e^{ \bom} }
{\nusym_n(X)\tau_I^{2g-2}})   $$
where $\a$ is the $\tor$-equivariant differential form on 
$K^{2g}\times\liek$ given by Proposition \ref{defa} which is supported
near $M(c)$ and represents the equivariant Poincar\'{e} dual of 
$M(c)$ in $K^{2g}\times\liek$.     \label{bag}
\end{lemma}
\Proof The $\tor$-equivariant differential form $\a$ on 
$K^{2g}\times\liek$ which 
represents the equivariant Poincar\'{e} dual of 
$M(c)=P^{-1}(c)$ in $K^{2g}\times\liek$ was defined in Proposition \ref{defa} 
as a pullback via the map $P:K^{2g}\times\liek \to K$. By using the 
restriction $P:H_I^{2g}\times\hat{\liet}_1 \to H_I$ we can similarly
define a $\tor$-equivariant differential form $\a_I$ on 
$H_I^{2g}\times\hat{\liet}_1$ which 
represents the equivariant Poincar\'{e} dual of 
$M(c)\cap(H_I^{2g}\times\hat{\liet}_1)$ in $H_I^{2g}\times\hat{\liet}_1$.
The restriction of $\Phi_{n-1}(\a_I)$ then represents the Poincar\'{e}
dual to $F_{I,\l}$ in $F$, provided suitable orientations are chosen.
Note that $\{1,\dots,1\}\times \liek$ is transverse to both
$M(c)=P^{-1}(c)$ and $\mu^{-1}(0)=K^{2g}\times\{0\}$ in
$K^{2g}\times \liek$, and that if $\L\in\liek$ then
$$\mu(1,\dots,1,\L) = -\L$$
while
$$P(1,\dots,1,\L) = \exp(-\L).$$
{}From the orientation conventions of Remark \ref{orient}
 it follows that the
normal to $P^{-1}(H_I)$ in $K^{2g}\times \liek$ is $\tor$-equivariantly
isomorphic to the kernel of the restriction map $\lieks \to {\bf h}_I^*$.
Thus the restriction of $(-1)^{r(n-r)}\tau_I \a_I$ to 
$H_I^{2g}\times\hat{\liet}_1$ has compact support near $M(c)$ and locally
represents the equivariant Poincar\'{e} dual to $M(c)$ in
$K^{2g}\times \liek$, so we can substitute it for $\a$ on
$H_I^{2g}\times\hat{\liet}_1$ and we can substitute
$(-1)^{r(n-r)}\Phi_{n-1}(\tau_I \a_I)$ for $\Phi_{n-1}(\a)$ on $F$.

We have that $e_F = (-1)^{r(n-r)g} \phicom(\tau_I^{2g}) $ by the 
last lemma. 
We therefore get 
$$ \int_{F } \frac{\phicom (
\eta e^{ \bom}\compform ) } {e_F}  = (-1)^{r(n-r)(g-1)}
\int_{F } \phicom (\frac{ 
\eta e^{ \bom}\a_I } {\tau_I^{2g-1}}) $$
$$ = (-1)^{r(n-r)(g-1)} 
\int_{F_{I,\l}} \phicom (\frac{ \eta e^{ \bom} }{\tau_I^{2g-1}}),$$
and Remark \ref{feb} completes the proof.
\hfill $\square$ 

\begin{rem} \label{fi}
The condition for $F\in \calf$ to appear in the sum in the
statement of Lemma \ref{l6.4} was that 
\beq \label{5.54}  -|\!|\he{1}|\!|^2 < \inpr{\he{1},\mu(F)} <0.  \eeq
Let $I$ be a subset of $\{3,\ldots,n\}$
with $r-1$ elements where $1\leq r\leq n-1$, and let
$\l \in \RR$ be a solution of the equation
$$ e^{-2\pi i \l} =  \prod_{j\in I_1} c_j.$$
If $F\in \calf$ is the component of the fixed point set 
 of the action of $\tone$ on  
$P^{-1}(V) \cap\mu^{-1}(\lietone)/\torc$ which
contains $\fisml$, then
$$\mu_{\tone}(F) = \mu_{\tone} (\fisml) =  -\l \hat{e}_1.$$
We thus find that for each $I$ there is precisely one
solution $\l\in\RR$ to the equation
$$ e^{-2\pi i \l} =  \prod_{j\in I_1} c_j$$
such that the component $F$ of the fixed point set 
 of the action of $\tone$ on  
$P^{-1}(V) \cap\mu^{-1}(\lietone)/\torc$ which
contains $\fisml$ contributes  to the sum in Lemma \ref{l6.4}.
This solution is
$\l = \d_I$ where $\d_I$ is the non-integer part of
$$\frac{i}{2\pi } \log \prod_{j\in I_1} c_j,$$
and so we have
$$\mu_{\tone}(F) = - \d_I \hat{e}_1.$$
(Note that since $\prod_{j\in I_1} c_j$ has modulus 1 but is not equal
to 1, the non-integer part of
$\frac{i}{2\pi} \log \prod_{j\in I_1} c_j$ is well defined as an element
of the open interval $(0,1)$ in $\RR$.) 
We therefore define
$$F_I = F_{I,\d_I},$$                               \label{bog}
and also $\Psi_I=\Psi_{I,\d_I}$
and
$\tpsi_I=\tpsi_{I,\d_I}.$
\end{rem}

We can now deduce the following result.

\begin{prop}                        \label{crucial} 
If $\eta(X)$ is a polynomial in 
the $\tar$ and $\tbr$, so that $s_{\he{1}}^*\eta  = \eta$, then 
$$ \int_{\nl(c)} \phil( \eta e^{\bom}  ) -
\int_{\nl(c)}  \phil (\eta e^{ \bom}e^{- Y_1 }  ) 
 = \int_{\nl(c)} \phil \Bigl ((1- e^{- Y_1}) 
\eta e^{ \bom}  \Bigr  ) $$
$$=  \sum_{1\leq r\leq n-1} 
\sum_{I\subseteq \{3,\ldots,n\},|I|=r-1} (-1)^{r(n-r)(g-1)}
\res_{Y_1 = 0 } \int_{F_I}
\phicom (\frac{
 \eta e^{ \bom} }{\tau_I^{2g-1}}). $$
\end{prop}
\Proof Recall that the coordinates $Y_k = e_k(X) 
= < \hat{e}_k,X>$ were introduced
in Definition \ref{YY}. The result then follows immediately
from Lemma \ref{l6.4}, Lemma \ref{bag} and Remark \ref{bog} above, 
together with Lemma
\ref{l3} and Remark \ref{n0=2}. \hfill $\square$ 

\begin{rem} \label{formal} This proposition is also true for formal
equivariant cohomology classes $\eta = \sum_{j=0}^{\infty} \eta_j$
with $\eta_j \in H^j_K(M(c))$, because all but finitely many
$\eta_j$ contribute zero to both sides of the equations.
\end{rem}

\begin{corollary} 
\label{c6.5} Suppose $\eta$ is a polynomial in 
the $\tar$ and $\tbr$, so that $s_{\he{1}}^*\eta  = \eta$. 
Then
$$ \int_{\nl(c)} \phil (\nusym_n
\eta e^{ \bom} ) =  \sum_{1\leq r\leq n-1} 
\sum_{I\subseteq \{3,\ldots,n\},|I|=r-1} (-1)^{r(n-r)(g-1)}
\res_{Y_1 = 0 } \int_{F_I}
\phicom (\frac{\nusym^{[I_1]}_r \nusym^{[I_2]}_{n-r}
 \eta e^{ \bom} }{\tau_I^{2g-2}( 1- e^{- Y_1} ) }). $$
\end{corollary}
\Proof This follows by applying Remark \ref{p5.11}
and Proposition \ref{crucial}
with $\eta$ replaced by the formal equivariant cohomology class
$\eta \nusym_n / (1 - e^{-Y_1}) $. This is valid by Remark
\ref{formal} because $Y_1$ divides $\nusym_n(X)$ and so
$ \nusym_n / (1 - e^{-Y_1}) $ can be expressed as a
power series in $Y_1$ whose coefficients are polynomials in
the other coordinates $Y_2,\ldots, Y_{n-1}$. \hfill $\square$

\bigskip

\begin{rem} \label{indbasrem}
Recall from the proof of Proposition \ref{p7.3} that
$$\tpsi_{I}: (S^1)^{2g} \times \Bigl (\mu_{SU(r)}^{-1}(0) \cap
M_r (c(I_1,\d_I))   \Bigr )
\times \Bigl (\mu_{SU(n-r)}^{-1}(0)  \cap 
M_{n-r} (c(I_2,-\d_I)) \Bigr )
\to \tilde{F}_{I,\d_I} $$ 
is defined for $\d_I$ as in Remark \ref{fi} by 
$$\tpsi_{I} \Bigl ( (s_1, \dots, s_{2g}), 
(h_1^{[I_1]}, \dots, h_{2g}^{[I_1]}, 0), 
(h_1^{[I_2]}, \dots, h_{2g}^{[I_2]}, 0) \Bigr ) $$ $$
= \Bigl ( (\r_I (s_1, \ha{1}, \hb{1}),\r_I (s_{2}, \ha{2}, \hb{2}),
\dots, \r_I (s_{2g}, \ha{2g}, \hb{2g}), 
2\pi i \diag(\d_I,-\d_I,0,\ldots, 0) \Bigr ) $$
using the map
$$\r_I: S^1 \times SU(r) \times SU(n-r) \to 
S(U(r) \times U(n-r)) \subset SU(n)$$
given by 
$$\r_I : (s,A,B)  \mapsto \matr{s^{n-r} A}{0}{0}{s^{-r} B} $$ 
with respect to the decomposition of $\{1,\ldots,n\}$ as $I_1\cup I_2$,
which restricts to an $r(n-r)$ to one surjective
homomorphism
$$\r_I:S^1 \times T_r \times T_{n-r} \to \tor.$$
Since $\bom = \omega + \mu$ is constructed using the inner product $<,>$ 
defined at (2.2) on the Lie algebra $\liek$ of $K=SU(n)$, and since
$\r_I$ embeds the Lie algebras of $S^1$, $SU(r)$ and $SU(n-r)$ as
mutually orthogonal subspaces of $\liek$, we have
$$\tpsi_I^*(\bom) = \bom_r + \bom_{n-r} + \Omega - \d_I \hat{e}_{1}$$ 
for some $\Omega\in H^2((S^1)^{2g})$, where $\bom_r$
and $\bom_{n-r}$ are defined like $\bom$ but with $n$ replaced by
$r$ and $n-r$. Thus we have
$$\tpsi_{I}^*(\frac{\nusym^{[I_1]}_r \nusym^{[I_2]}_{n-r}
 \eta e^{ \bom} }{\tau_I^{2g-2}(1 - e^{- Y_1}) })
=e^{ \bom_r + \bom_{n-r} + \Omega - \d_I Y_1} 
\tpsi_{I}^*(\frac{\nusym^{[I_1]}_r \nusym^{[I_2]}_{n-r}
 \eta }{\tau_I^{2g-2}(1 - e^{- Y_1}) }).$$
\end{rem}

Since $\tpsi_I^*(\nusym^{[I_1]}_r) = \nusym_r$ and
$\tpsi_I^*(\nusym^{[I_2]}_{n-r}) = \nusym_{n-r}$, 
we can combine this with Corollary \ref{c6.5} and Remark
\ref{big} to obtain the result on which is 
based the inductive proof of Witten's formulas in the next section.

\begin{prop}                        \label{indbas} 
If $c\in T$ satisfies the conditions of Remark
\ref{induct}, and if $\eta(X)$ is a polynomial in 
the $\tar$ and $\tbr$ so that $s_{\he{1}}^*\eta  = \eta$, then 
$$ \int_{\nl(c)} \phil (\nusym_n
\eta e^{ \bom} ) =  \sum_{1\leq r\leq n-1} 
\sum_{I\subseteq \{3,\ldots,n\},|I|=r-1} (-1)^{r(n-r)(g-1)}
\res_{Y_1 = 0 } \int_{F_I}
\phicom (\frac{\nusym^{[I_1]}_r \nusym^{[I_2]}_{n-r}
 \eta e^{ \bom} }{\tau_I^{2g-2}(1- e^{- Y_1} ) })$$
where
$$\res_{Y_1 = 0 } \int_{F_I}
\phicom (\frac{\nusym^{[I_1]}_r \nusym^{[I_2]}_{n-r}
 \eta e^{ \bom} }{\tau_I^{2g-2}(1- e^{- Y_1}) })$$
is equal to $(r(n-r))^{-2g}$ times
$$ 
\int_{N_r(c(I_1,\d_I)) } \Phi_r (\nusym_r e^{\bom_r} \res_{Y_1^I=0} 
\int_{N_{n-r}(c(I_2,-\d_I))} \Phi_{n-r}( 
\nusym_{n-r} e^{\bom_{n-r}} \int_{(S^1)^{2g} } e^{\Omega}
\tpsi_{I}^*(\frac{
 \eta e^{- \d_I Y_1} }{\tau_I^{2g-2}(1- e^{- Y_1} ) })))$$                     
and also to $(r(n-r))^{-2g}$ times
$$\int_{N_{n-r}(c(I_2,-\d_I)) } \Phi_{n-r} (\nusym_{n-r} e^{\bom_{n-r}}
\res_{Y_1^I=0} 
\int_{N_r(c(I_1,\d_I))} \Phi_{r}(\nusym_r e^{\bom_r} \int_{(S^1)^{2g} }
e^{\Omega}
\tpsi_{I}^*(\frac{
 \eta e^{-\d_I Y_1} }{\tau_I^{2g-2}(1- e^{- Y_1}) }))).$$  
Here $c(I_1,\d_I))$ and $c(I_2,-\d_I)$ are defined as in Proposition
\ref{p7.3} with $\d_I$ as in Remark \ref{bog} and $\bom_r$, $\bom_{n-r}$
and $\Omega$ as in Remark \ref{indbasrem}.
\end{prop}

\begin{rem}  \label{c} For any $\g \in \tor$ a unique $\tilde{\g}
\in \liet_n$ can be chosen so that $\exp \tilde{\g}=\g$
and $\tilde{\g}$ belongs to the fundamental domain defined
by the simple roots for the translation action on
$\liet_n$ of the integer lattice $\L^I$ (i.e. $\tilde{\g} = \g_1\hat{e}_1
+ \ldots + \g_{n-1}\hat{e}_{n-1}$ with $0\leq \g_j <1$ for $1\leq j\leq n-1$).
Suppose that $\tilde{c}(I_1,\d_I)\in\liet_r$ and 
$\tilde{c}(I_2,-\d_I)\in\liet_{n-r}$ are chosen in this way in the fundamental 
domains defined by the simple roots for the translation actions on
$\liet_r$ and $\liet_{n-r}$ of their integer lattices, satisfying
$$\exp \tilde{c}(I_1,\d_I) = c(I_1,\d_I)= \diag(c^{I,\d_I}_{i_1},\ldots,
c^{I,\d_I}_{i_r})$$
and
$$\exp \tilde{c}(I_2,-\d_I) = c(I_2,-\d_I) =
\diag(c^{I,\d_I}_{i_{r+1}}, \ldots, c^{I,\d_I}_{i_n}),$$
where  (as in Proposition \ref{p7.3} and Remark \ref{fi}) we define
$\d_I$ to be the non-integer part of $\frac{i}{2\pi}\log \prod_{j\in I_1}
c_j$ and let
$ c^{I,\d_I}_j = c_j$
if $j\geq 3$, and
$ c^{I,\d_I}_{1} = c_{1} e^{2\pi i \d_I}$
and
$c^{I,\d_I}_2 = c_2 e^{-2\pi i \d_I}.$

In the proof of the main theorem (Theorem \ref{mainab}) of the next section
we shall need to consider the elements
$w^1_I$ and $w^2_I$ of the subgroup $S_{n-1}$ of
the Weyl group  $W\cong S_n$ of $SU(n)$
given by the permutations
$$ \left \lbrack \begin{array}{lcccr}
1&2&\dots & n \\
i_1 & i_2 & \dots &i_n
\end{array} \right \rbrack $$
and
$$ \left \lbrack \begin{array}{lcccccr}
1&2&\dots & n-r & n-r+1 & \dots & n \\
i_{r+1} & i_{r+2} & \dots &i_n &i_1 & \dots & i_r
\end{array} \right \rbrack, $$
in the cases when $i_r=1$ and $i_{r+1}=2$ and $i_n=n$ and when $i_1=1$
and $i_n=2$ and $i_r=n$ respectively. We will use the fact that
if $i_r=1$ and $i_{r+1}= 2$ and $i_n=n$ then
$$w^1_I(\tc)  = 
\matr{\tilde{c}(I_1,\d_I)}{0}{0}{\tilde{c}(I_2,-\d_I)}
+ (1 -\d_I) \hat{e}_{1}, $$
where the block diagonal form is taken
with respect to the decomposition of $\{1,\ldots,n\}$ as $
\{1,\ldots,r\}\cup \{r+1,\ldots,n\} $.
To see why this is the case, note that
$$w^1_I(\tc)(X) = \g_1 (X_{i_1} - X_{i_2}) + \dots \g_{n-1}
 (X_{i_{n-1}} - X_{i_n})$$
where 
$\g_k$ is the non-integer part of $\frac{1}{2\pi i} {\rm log} 
\prod_{j\leq k} c_{i_j},$         
so that $\g_r = 1 - \d_I$ and if $k<r$
then $\g_k$ is the non-integer part of $\frac{1}{2\pi i} {\rm log} 
\prod_{j\leq k} c_{i_j}^{I,\d_I}$         
whereas if $k>r$ then $\g_k$ 
is the non-integer part of 
$$-\d_I + \frac{1}{2\pi i} {\rm log} 
\prod_{r < j\leq k} c_{i_j} = \frac{1}{2\pi i} {\rm log} 
\prod_{r < j\leq k} c_{i_j}^{I,\d_I}.$$
Similarly if $i_1 =1$ and $i_n =2$ and $i_r = n$ then
$$w^2_I(\tc)  = 
\matr{\tilde{c}(I_2,-\d_I)}{0}{0}{\tilde{c}(I_1,\d_I)} - \d_I \hat{e}_{1}$$
where the block diagonal form is taken 
with respect to the decomposition of $\{1,\ldots,n\}$ as $
\{1,\ldots,n-r\}\cup \{n-r+1,\ldots,n\} $.
\end{rem}

\renorm
\section{Proof of the iterated residue formula}

In this section we shall use induction to prove Witten's formulas in the 
formulation
given in Section 2 (see Proposition \ref{p:sz}) 
involving iterated residues, for 
pairings  of the form  
\beq \prod_{r=2}^n a_r^{m_r}
 \prod_{k_r=1}^{2g} (b_r^{k_r})^{p_{r,k_r}}\exp (f_2) 
[\mnd] \eeq
for nonnegative integers $m_r$ and $p_{r,k}$. The induction is based on
Proposition \ref{indbas}.
In the next section we shall extend the proof to give formulas for all 
pairings, and in the following section we shall show that these formulas are 
equivalent to those of Witten.

We are aiming to prove

\begin{theorem} \label{mainab} Let $c=\diag\, (e^{2\pi i d/n},\ldots,
e^{2\pi i d/n})$ where $d \in\{1,\ldots,n-1\}$ is coprime to $n$,
and suppose that $\eta\in H^*_{SU(n)}(M_n(c))$
is a polynomial $Q(\tilde{a}_2,\ldots,\tilde{a}_n,
\tilde{b}_2^1,\ldots,\tilde{b}_n^{2g})$ in the 
equivariant cohomology classes $\tilde{a}_r$
and $\tilde{b}_r^j$ for $2\leq r\leq n$ and $1\leq j\leq 2g$.
Then the pairing
$Q(a_2,\ldots,a_n,b_2^1,\ldots,b_n^{2g})\exp (f_2) [\mnd]$
is given by
$$\int_{\mnd} \Phi (\eta e^{\bom} ) 
 =  \frac{(-1)^{n_+(g-1)}}{n!} \res_{Y_{1} = 0} \dots \res_{Y_{n-1} = 0} 
\Biggl ( \frac{\sum_{w \in W_{n-1}} e^{ 
\inpr{ \tildarg{w \tc},X}  } \int_{\tor^{2g}} \eta
e^{ \omega} } {\nusym_n^{2g-2} \prod_{1\leq j \leq n-1}
( \exp (Y_j)-1 ) } \Biggr ), $$           
where $n_+ = \frac{1}{2} n(n-1)$ is the number of positive
roots of $K=SU(n)$ and $X\in\tor$ has 
coordinates $Y_1=X_1-X_2,\ldots,Y_{n-1}=X_{n-1}-X_n$ 
defined by
the simple roots, while $W_{n-1} \cong S_{n-1}$ is the Weyl
group of $SU(n-1)$ embedded in $SU(n)$ in the standard way
using the first $n-1$ coordinates. The element
$\tilde{c}$ was defined in Remark \ref{r2.1}: it is the
unique element of $\liet_n$  which 
satisfies $e^{2\pi i\tilde{c}} = c$ and belongs to the
fundamental domain defined by the simple roots for the translation
action on $\liet_n$ of the integer lattice $\L^I$.
Also,  the notation $\bracearg{\gamma}$ (introduced 
in Definition \ref{bracedef}) means the unique element which is in the
fundamental domain defined by the simple roots for the
translation action  on $\liet_n$ of the integer lattice and for
 which $\bracearg{\gamma}$ is equal to $\gamma$ plus some element of the
integer lattice.
\end{theorem}   

\begin{rem} \label{interp} Here the integral
$$ \int_{\tor^{2g}} \eta
e^{ \omega}  $$
is to be interpreted as the integral of the restriction of
$\eta e^{ \omega} $ over a
connected component
$$\tor^{2g} \times \{ \l\}$$
(for some $\l \in \liet_n$ satisfying
$c\exp \l = 1$) of the
fixed point set of the action of $\tor$
on $M_{n}(c)$. It
does not matter which component we choose here, because $\eta$ and
$\omega$ are invariant under the
translation maps $s_{\L_0}$ defined at Lemma \ref{l4.3}
for  $\L_0$ in the integer lattice
of $\liet_n$. 
\end{rem}

\begin{rem} \label{rem8.3} (a) We can 
substitute $-X$ for $X$ in Theorem \ref{mainab} to get
$$\int_{\mnd} \Phi (\eta e^{\bom} ) 
 =  \frac{(-1)^{n_+(g-1)}}{n!} 
\res_{Y_{1} = 0} \dots \res_{Y_{n-1} = 0} 
\Biggl ( \frac{\sum_{w \in W_{n-1}} e^{
- \inpr{\tildarg{w\tc},X}  } \int_{\tor^{2g}} \eta(-X)
e^{ \omega} } {\nusym_n^{2g-2} \prod_{1\leq j \leq n-1}
(1 - \exp (-Y_j)) } \Biggr ). $$           

\noindent (b) When $\eta$ is a polynomial in $a_2,\dots,a_n$ then
$$ \int_{\tor^{2g}} \eta
e^{ \omega}  = \eta \int_{\tor^{2g}} 
e^{ \omega} = n^g \eta $$
(see Lemma \ref{l9.3b} below). Since $\tilde{a}_r$ is represented
by the polynomial $\tau_r(-X)$ for $2\leq r\leq n$ (see
Proposition \ref{abftil} or Section 9 below), this means that, by (a) above,
 Theorem \ref{mainab} combined with Proposition
\ref{p:sz} gives us Witten's formula (2.4).

\noindent (c) We can also replace the symplectic form $\omega$
by any nonzero scalar multiple $\e \omega$. Then the moment map $\mu$
is multiplied by the same scalar $\e$, and the proof of
Theorem \ref{mainab} yields
$$\int_{\mnd} \Phi (\eta e^{\e \bom} ) 
 =  \frac{(-1)^{n_+(g-1)}}{n!} \res_{Y_{1} = 0} \dots \res_{Y_{n-1} = 0} 
\Biggl ( \frac{\sum_{w \in W_{n-1}} e^{
\inpr{\e \tildarg{w\tc},X}  } \int_{\tor^{2g}} \eta
e^{\e \omega} } {\nusym_n^{2g-2} \prod_{1\leq j \leq n-1}
( \exp (\e Y_j)-1) } \Biggr ). $$   
If the degree of $\eta$ is equal to the dimension of $\mnd$
then the left hand side of this equation is equal to
$$\int_{\mnd} \Phi (\eta)$$
and hence is independent of $\e$. Thus in this case we can take
any nonzero value of $\e$ on the right hand side, or let
$\e$ tend to zero, to give alternative formulas for
$\int_{\mnd} \Phi (\eta)$.
\end{rem}  

Recall from Lemma \ref{l3} that
\beq\label{oldl5.10} \int_{\mnd} \Phi (\eta e^{\bom} ) = \frac{1}{n!}
\int_{\nl(c)} \phil (\nusym_n \eta e^{ \bom} ) .\eeq
Proposition \ref{indbas} tells us that 
$\int_{\nl(c)} \phil (\nusym_n \eta e^{ \bom} )$ can be expressed in terms
of iterated integrals of the same form for smaller values of $n$,
but with $c$ no longer central in $K=SU(n)$.
We shall therefore obtain Theorem \ref{mainab} from the following
result involving values of $c$ which are not central (cf. Remark
\ref{induct}), which will be proved by induction on $n$.

\begin{prop} \label{beg} Let $c=\diag(c_1,\ldots,c_n) \in \tor$ be such that
the product of no proper subset of $c_1,\ldots,c_n$ is 1. 
If $\eta(X)$ is a polynomial in 
the $\tar$ and $\tbr$, so that $s_{\he{l}}^*\eta  = \eta$, then 
$$\int_{\nl(c)} \phil (\nusym_n \eta e^{ \bom} ) = (-1)^{n_+(g-1)}
\res_{Y_{1} = 0} \dots \res_{Y_{n-1} = 0} 
\Bigl ( \frac{\sum_{w \in W_{n-1}} e^{
\inpr{\tildarg{w\tc} ,X}  } \int_{\tor^{2g}} \eta
e^{ \omega} } {\nusym_n^{2g-2} \prod_{1\leq j \leq n-1}
( \exp (Y_j)-1 ) } \Bigr ), $$           
where $W_{n-1} \cong S_{n-1}$ is the Weyl group of $SU(n-1)$,
embedded in $SU(n)$ in the standard way using the first
$n-1$ coordinates, and $\tc = (\tc_1,\dots,\tc_n) \in \liet_n$ satisfies 
$e^{2\pi i \tc}=c$ and belongs to the fundamental domain defined
by the simple roots for the translation action on
$\liet_n$ of the integer lattice $\L^I$.
\end{prop}

\noindent{\bf Proof of Theorem \ref{mainab} from Proposition \ref{beg}:}
Note that when $c = \diag(e^{2\pi i d/n},\ldots, e^{2\pi i d/n})$ 
we had introduced an element
$\tilde{c} $ $\in$ 
$\liet_n$ (see Remark \ref{r2.1}) which 
satisfies $e^{2\pi i\tilde{c}} = c$ and 
 belongs to the fundamental domain defined
by the simple roots for the translation action on
$\tor$ of the integer lattice $\L^I$.
Thus 
Theorem \ref{mainab} follows immediately from
(8.2) and Proposition \ref{beg}. \hfill $\square$

\bigskip

\noindent {\bf Proof of Proposition \ref{beg}:} 
The proof is by induction on $n$. 
When $n=1$ then both $SU(n)$ and the torus $\tor$ are trivial, 
$\nusym_n = 1$ and both $M_n(c)$ and
$N_{n}(c)$ are single
points. Thus in this case Proposition \ref{beg} reduces to the tautology
$ \eta =  \eta$
for any $\eta \in H^*_{SU(1)}(M_1(c))$.

Now let us assume that $n>1$ and that the result is true for all smaller 
values of $n$. 
By Proposition \ref{indbas} we have
$$ \int_{\nl(c)} \phil (\nusym_n
\eta e^{ \bom} ) =  \sum_{1\leq r\leq n-1} 
\sum_{I\subseteq \{3,\ldots,n\},|I|=r-1} (-1)^{r(n-r)(g-1)}
\res_{Y_1 = 0 } \int_{F_I}
\phicom (\frac{\nusym^{[I_1]}_r \nusym^{[I_2]}_{n-r}
 \eta e^{ \bom} }{\tau_I^{2g-2}(1- e^{- Y_1}) })$$
where
$$\res_{Y_1 = 0 } \int_{F_I}
\phicom (\frac{\nusym^{[I_1]}_r \nusym^{[I_2]}_{n-r}
 \eta e^{ \bom} }{\tau_I^{2g-2}(1- e^{ -Y_1}) })$$
is equal to
$ (r(n-r))^{-2g}$ times the iterated integral
$$\int_{N_r(c(I_1,\d_I)) } \Phi_r (\nusym_r e^{\bom_r} \res_{Y_1^I=0} 
\int_{N_{n-r}(c(I_2,-\d_I))} \Phi_{n-r}(\nusym_{n-r}
e^{\bom_{n-r}} \int_{(S^1)^{2g} } e^{\Omega}
\tpsi_{I}^*(\frac{
 \eta e^{- \d_I Y_1} }{\tau_I^{2g-2}( 1-e^{ -Y_1}) })))$$                     
and also to $(r(n-r))^{-2g}$ times the iterated integral
$$\int_{N_{n-r}(c(I_2,-\d_I)) } \Phi_{n-r} (
\nusym_{n-r} e^{\bom_{n-r}} \res_{Y_1^I=0} 
\int_{N_r(c(I_1,\d_I))} \Phi_{r}( 
\nusym_r e^{\bom_r} \int_{(S^1)^{2g} } e^{\Omega}
\tpsi_{I}^*(\frac{ 
 \eta e^{- \d_I Y_1} }{\tau_I^{2g-2}( 1- e^{- Y_1}) }))),$$  
for $c(I_1,\d_I)$ and $c(I_2,-\d_I)$ defined as in Proposition
\ref{p7.3} with $\d_I$ as in Remark \ref{bog} and $\tau_I$
as in Remark \ref{feb}. Here $\Omega \in H^2((S^1)^{2g})$ satisfies
$$\tpsi_{I}^*(\bom) = \bom_r + \bom_{n-r} + \Omega - \d_I \hat{e}_1$$
as in Remark \ref{indbasrem}.

We need to consider separately those $I$ containing $n$ and those
for which $n$ is not an element of $I$; first let us suppose that
$n$ is not an element of $I$.  Note that
$$(-1)^{(r(n-r)+ \frac{1}{2}r(r-1) + \frac{1}{2}(n-r)(n-r-1))(g-1)}
=(-1)^{\frac{1}{2}n(n-1)(g-1)},$$
and
$$\frac{e^{-\d_I Y_1}}{1-e^{-Y_1}} =
\frac{e^{(1-\d_I) Y_1}}{e^{Y_1}-1} .$$
The finite cover $\r_I:S^1 \times T_r \times T_{n-r} \to \tor$
is $r(n-r)$ to one, so that it induces an $(r(n-r))^{2g}$ to one
surjection from $(S^1)^{2g} \times T_r^{2g} \times T_{n-r}^{2g}$
to $\tor^{2g}$ and we have
$$\int_{T_n^{2g}} \eta e^{\omega} = \int_{(S^1)^{2g}\times T_r^{2g}
\times T_{n-r}^{2g}} \eta e^{\omega_r + \omega_{n-r} + \Omega}.$$
Moreover this finite cover $\r_I:S^1 \times T_r \times T_{n-r} \to \tor$
takes the coordinate $Y_1 = X_1 - X_2$ 
on $\liet$ to the coordinate $Y^I_1$ on the Lie algebra of $S^1$.
Since $\tpsi_I^*$ was defined using $\r_I$ (see Remark \ref{big}),
we deduce using Remark \ref{indbasrem}  and Remark  \ref{p5.11}
 and induction on $n$
that $(-1)^{r(n-r)(g-1)}$ times the
iterated integral
$$ \int_{N_r(c(I_1,\d_I)) } 
\Phi_r (\nusym_r e^{\bom_r} \res_{Y_1^I=0} 
\int_{N_{n-r}(c(I_2,-\d_I))} \Phi_{n-r}(\nusym_{n-r}
e^{\bom_{n-r}} \int_{(S^1)^{2g} } e^{\Omega}
\tpsi_{I}^*(\frac{
 \eta e^{- \d_I Y_1} }{\tau_I^{2g-2}( 1- e^{- Y_1}) })))$$                     
equals $(-1)^{n_+(g-1)} (r(n-r))^{2g} $ times the iterated residue
$$\res_{X_{i_{1}} - X_{i_{2}}=0} \dots \res_{X_{i_{r-1}} -
X_{i_{r}} =0} \res_{X_1-X_2 =0} \res_{X_{i_{r+1}} - X_{i_{r+2}} =0}
\dots \res_{X_{i_1} - X_{i_2} =0} $$
$$\sum_{w_1\in W_{r-1}}\sum_{w_2\in W_{n-r-1}} \frac{ e^{
\inpr{\tildarg{ {w}_1\tc(I_1,\d_I) } ,Y_{I_1}}} 
e^{ \inpr{ \tildarg{ {w}_2 \tc(I_2,-\d_I ) } ,Y_{I_2}} } 
e^{(1-\d_I) Y_1}
\int_{T_{n}^{2g}} \eta e^{\omega}  }
{ \nusym_{n}^{2g-2}
( e^{Y_1}-1) \prod_{j\neq r}
( \exp (X_{i_{j}}-X_{i_{j+1}})-1 ) } \Biggr ) $$
where $Y_{I_1}$ and $Y_{I_2}$ are the projections of $X$ onto
the Lie algebras of the maximal tori $T_r$ and $T_{n-r}$ of
$SU(r)$ and $SU(n-r)$ embedded in $SU(n)$ via the decomposition of
$\{1,\ldots,n\}$ as $I_1 \cup I_2$, and $W_{r-1}$ and $W_{n-r-1}$
are the Weyl groups of $SU(r-1)$ and $SU(n-r-1)$ embedded
in $SU(r)$ and $SU(n-r)$ using all but the last coordinates. 

There is no need to assume that
$i_1 < i_2 < \dots <i_r$ and
$i_{r+1} < i_{r+2} < \dots < i_n$ here. We simply
need that $I_1 = I \cup \{1\} = \{i_1,\ldots,i_r\}$
and $I_2 = \{1,\ldots,n\} - I_1 = \{i_{r+1},\ldots, i_n\}$.
So let us assume that $$i_r=1$$ and $$i_{r+1} =2.$$
We are also supposing that $n$ is not an element of $I$
(i.e. that $n\in I_2$) so we may assume in
addition that $i_n = n$. Then we can apply the Weyl transformation
$w_I^1 \in W_{n-1}$ given by the permutation
$$\left \lbrack \begin{array}{ccccc}
1 & \dots & r & \dots & n-1 \\
i_1 & \dots & i_r & \dots & i_{n-1}
\end{array} \right \rbrack $$
together with Remark \ref{c} to identify
the iterated residue above with
$$(-1)^{n_+(g-1)} (r(n-r))^{2g}
\res_{Y_{1}=0} \dots \res_{Y_{n-1} =0} \sum_{w_1\in W_{r-1}}
\sum_{w_2\in W_{n-r-1}}  
\frac{ e^{\inpr{\tildarg{w^1_I w_1 w_2(\tc) },X}} 
\int_{T_{n}^{2g}}\eta
e^{\omega} } 
{ \nusym_{n}^{2g-2} \prod_{1\leq j \leq n-1}
(1- \exp (-Y_j)) }  . $$    

When $n\in I$ the argument is similar but we apply induction to
$(-1)^{r(n-r)(g-1)}$ times 
$$ \int_{N_{n-r}(c(I_2,-\d_I)) } 
\Phi_{n-r} (\nusym_{n-r} e^{\bom_{n-r}} \res_{Y_1^I=0} 
\int_{N_{r}(c(I_1,\d_I))} \Phi_{r}(\nusym_{r}
e^{\bom_{r}} \int_{(S^1)^{2g} } e^{\Omega}
\tpsi_{I}^*(\frac{
 \eta e^{ -\d_I Y_1} }{\tau_I^{2g-2}(1 - e^{- Y_1}) })))$$
and observe that
$$\res_{X_1-X_2=0} \frac{e^{-\d_I(X_1-X_2)}}{1-e^{-(X_1-X_2)}} =  
- \res_{X_2-X_1=0} \frac{e^{\d_I(X_2-X_1)}}{1-e^{X_2-X_1}} $$
(see the Remark after Corollary \ref{c4.2}). 
As $I_1 = I \cup \{1\} = \{i_1,\ldots,i_r\}$
and $I_2 = \{1,\ldots,n\} - I_1 = \{i_{r+1},\ldots, i_n\}$
and $n \in I$ we can assume that $i_1 =1$, $i_r = n$ and $i_n=2$.
Then we use  the Weyl transformation
$w^2_I \in W_{n-1}$ given by the permutation
$$\left \lbrack \begin{array}{ccccccc}
1 & \dots & n-r & n-r+1 & \dots & n-1 \\
i_{r+1} & \dots & i_n & i_1 & \dots & i_{r-1}
\end{array} \right \rbrack $$
together with Remark \ref{c} to equate
the iterated integral above with
$$  
\res_{Y_{1}=0} \dots \res_{Y_{n-1} =0} \sum_{w_1\in W_{r-1}}
\sum_{w_2\in W_{n-r-1}}  
\frac{ e^{\inpr{ \tildarg{  w^2_I w_1w_2 (\tc) },X}} 
\int_{T_{n}^{2g}}\eta
e^{\omega} } 
{ \nusym_{n}^{2g-2} \prod_{1\leq j \leq n-1}
( \exp (Y_j)-1) }  . $$         
Thus it suffices to prove

\begin{lemma} For each subset $I$ of $\{3,\ldots,n\}$ with $r-1$
elements, let us fix $i_1,\ldots,i_n$ such that
$I\cup \{1\} = \{i_1,\ldots, i_r\}$ and
$\{2,\ldots,n\} - I = \{i_{r+1},\ldots, i_n\}$ and
also $i_r=1$, $i_{r+1}=2$ and $i_n=n$ (if $n \not\in I$)
or $i_1=1$, $i_r=n$ and $i_n=2$ (if $n \in I$). Define
permutations $w_I^1$ (for $I$ such that $n \not\in I$) and
$w^2_I$ (for $I$ such that $n \in I$) as above. 
Then as 

(i) $r$ runs over $\{1,\ldots,n-1\}$,

(ii) $w_1$ runs over permutations
of $\{1,\ldots,r\}$ fixing $r$,

(iii) $w_2$ runs over permutations
of $\{r+1,\ldots,n\}$ fixing $n$ and

(iv) $I$ runs over subsets of
$\{3,\ldots,n\}$ with $r-1$ elements not containing $n$, 

\noindent the product $w^1_I w_1 w_2 $
runs over the set of permutations
$w$ of $\{1,\ldots,n\}$ fixing $n$ such that
$$w^{-1}(1)<w^{-1}(2).$$
Moreover if instead of (iv) $I$ runs over subsets of
$\{3,\ldots,n\}$ with $r-1$ elements containing $n$, 
then the product $w^2_I w_1 w_2 $
runs over the set of permutations
$w$ of $\{1,\ldots,n\}$ fixing $n$ such that
$$ w^{-1}(1)>w^{-1}(2).$$
\end{lemma}
\Proof If $w \in W_{n-1}$ satisfies $w^{-1}(1)<w^{-1}(2)$ let
$r=w^{-1}(1)$ and $I=\{j:w^{-1}(j)<r\}$. On the other hand if
$w \in W_{n-1}$ satisfies $w^{-1}(1) > w^{-1}(2)$ let
$r=n-w^{-1}(2)$ and $I=\{j>1:w^{-1}(j) > n-r\} \cup \{n\}$.
In each case it is easy to check
that there exist unique choices of $w_1$ and $w_2$
such that $w^1_I w_1 w_2  =w$ or $w^2_I w_1 w_2 =w.$ 

This completes the proof of the lemma and hence of Proposition
\ref{beg}.

\begin{rem} \label{naive} It is 
shown in Proposition 3.4 of \cite{JK3} that the 
multivariable residue (multiplied by the constant
$C_K$) of Theorem \ref{t4.1} and formula (\ref{1.7}) can be 
replaced by the iterated one-variable residue
$$\res^+_{Y_{1}=0} \dots \res^+_{Y_{n-1}=0} $$
multiplied by the Jacobian (in this case $1/n$) of the change
of coordinates from an orthonormal system to
$(Y_1,\ldots,Y_{n-1})$. Here, if
$\res_{y=0}g(y)$ denotes the coefficient of $y^{-1}$ in the
Laurent expansion about 0 of a meromorphic function $g(y)$ of one
complex variable $y$, then $\res^+$ is defined 
for meromorphic functions of the special form
$\sum_{1\leq i\leq s} e^{\l_i y} q_i(y)$,
where $\l_1,\ldots,\l_s$ are real numbers and $q_1,\ldots,q_s$
are rational functions of one variable with complex coefficients, by
$$\res^+_{y=0} (\sum_{1\leq i\leq s} e^{\l_i y} q_i(y))
= \sum_{1\leq i\leq s, \l_i>0} \res_{y=0}(e^{\l_i y} q_i(y)).$$
Since 
$$\frac{e^{\g y}}{e^{y}-1}$$
can be formally expanded as
$$-\sum_{m\in \ZZ, m+\g >0} e^{(m+\g)y}$$
when $0 < \g<1$, the formula (\ref{1.7}) can be formally rewritten as
$$\prod_{r=2}^n a_r^{m_r} \exp(f_2) [\mnd] = \frac{(-1)^{n_+(g-1)}}{n}
\res^+_{Y_{1}=0} \dots \res^+_{Y_{n-1}=0} 
\frac{
e^{\inpr{ \tc,X} } 
\int_{T_n^{2g}} e^{\omega} \prod_{r=2}^n \tau_r^{m_r}}
{\nusym_n^{2g-2} \prod_{1\leq j\leq n-1} (e^{Y_j} -1)}.$$
Moreover the multivariable residue Res is
invariant under the action of the Weyl group,
as are all the other ingredients of the right
hand side of (\ref{1.7}) except for $\tc$. Thus by averaging
(\ref{1.7}) over
 the Weyl group we obtain a special case of Theorem \ref{mainab}.
\end{rem}

\renorm
\section{Residue formulas for general intersection pairings}
\nc{\abk}{{\Phi} } 
\nc{\tq}{\tilde{q}}
\nc{\hattau}{\hat{\tau}}
\nc{\srj}{{s_r^j} }
\nc{\ssjg}{{s_s^{j+g} } }
\nc{\psitexp}{ {\check{\psi}_{X,B}} }
\nc{\tfq}{{ \tilde{f}_{ (q) } }}
\nc{\xsharp}{{X^\#}}
\nc{\free}{\FF}
\nc{\epmb}{{\epc}}
\nc{\thetsimp}[1]{\theta^{(#1)}  } 
\nc{\spart}[2]{  \frac{\partial}{\partial s_{#1}^{#2} }  }

\nc{\mom}{J}

\nc{\conn}[1]{ {\Theta^{(#1)} } }

\nc{\qsgn}{{q}}
\nc{\qought}{{q_{(o)}  }}
\nc{\tsign}{{\tau  }}
\nc{\Ssign}[1]{{ S^{#1}}}
\nc{\signr}{ { (-1)^r }}
\nc{\signs}{ { (-1)^s }}

\nc{\signq}{**sign?**}
\nc{\gen}{{\z }}
\nc{\hu}[1]{ \hat{u}_{#1} }
\nc{\zz}[1]{ Z_{#1} }
\nc{\pqab}{ (\partial^2 \qsgn)(\hu{a}, \hu{b})  }
\nc{\pqabx}{ (\partial^2 \qsgn)_X(\hu{a}, \hu{b})  }
\nc{\simp}{\bigtriangleup}
\nc{\quest}{**??**}

\nc{\tarnox}{{\tilde{a}_r} }
\nc{\tbrnox}{{\tilde{b}_r^j}}
\nc{\hess}{{H_\liet}}

In order to obtain explicit formulas for all the pairings, Witten
 observes that they can be obtained from those for the 
$a_r$ and $f_r$ via his formula
\cite{tdgr}  (5.20). In this section we shall generalize our version 
of his formula (\cite{tdgr} (4.74), which
is our Theorem \ref{mainab} via the results of Section 2)
   to give formulas for 
$\int_{\mnd} \abk (\eta e^\bom)$
where $\eta$ is an equivariant cohomology class that
does not simply involve the 
$\tar$ but also involves the $\tbrj $ and 
the 
$\tfr$  (see Theorems \ref{t9.5} and \ref{t9.6} below).
The key step in the proof is 
Lemma \ref{l9.2}, combined with the argument used in Sections 5-8
to prove Theorem \ref{mainab}. 

In the next section we shall see that  Theorem \ref{t9.5} yields 
Witten's formula \cite{tdgr} (5.20).
This will follow from certain 
equations satisfied by 
the formula given in 
Theorem \ref{t9.4} (Propositions \ref{p9.1} and \ref{p9.2}).

 The next lemma (from \cite{J2}) will 
give an explicit formula for an  equivariant cohomology 
class $\tfr$ on $M(c)$  such that
$\abk(\tfr) = f_r$ (cf. 
Proposition \ref{abftil}). 
In order to state it,
we introduce the  following notation.

\begin{definition} {\bf (The moment)}
If  $\theta$ is the Maurer-Cartan form on $K$,
the moment 
$\mom(\theta) \in \Bigl ( \Omega^1 (K) \otimes \lieks \Bigr )^K $ 
is defined for  $X \in \liek$  by 
\beq \label{9.09} 
\mom(\theta) (X)_k  =  - \iota_{\xsharp} \theta =  -\Ad (k^{-1})
X, \eeq
where $\xsharp$  is the vector field on $K$ given by the left
action of $X$ on $K$.
\end{definition}

\begin{rem} See \cite{BGV}, Chapter 7 for an explanation of the role of the 
moment in the construction of equivariant characteristic 
classes, via an equivariant version of Chern-Weil theory. 
Given a  principal bundle over a $K$-manifold
equipped with a compatible 
 action of $K$ on the total space of the bundle, 
the moment $\mom$ 
plays the same role as  the symplectic
 moment 
map plays for a principal $U(1)$ bundle $\call$
 over a Hamiltonian $K$-manifold
with  $c_1(\call) = [\omega]$
(and with a lift of the action of $K$ to the total space of $\call$). 
In particular, the appropriate notion of ``equivariant curvature''
is the sum of the  usual curvature and the moment $\mom$.
\end{rem}

In the next few paragraphs we provide a brief outline
of the use of the Bott-Shulman construction 
(see for instance \cite{BSS} and other references given 
in \cite{J2}) to obtain
equivariant differential forms representing the equivariant
characteristic classes $\tfr$.  This material is summarized
from  \cite{J2}, which gives a construction of 
 de Rham representatives for equivariant
characteristic classes giving rise to the
characteristic classes of the universal bundle over $\mnd \times \Sigma$.
This was accomplished by regarding this bundle (and the classifying
space for it) as {\em simplicial manifolds}. For more details see 
\cite{J2}.

Let $\simp^2  =  \{ (t_0, t_1, t_2) \in [0,1]^3: 
t_0 + t_1 + t_2 = 1 \}  $ be the standard 2-simplex.
There is a 
principal 
$K$-bundle 
$$ \simp^2 \times K^3 \stackrel{\pi_2}{\longrightarrow} \simp^2 \times
K^2 $$ for which the bundle projection 
 $\pi_2: K^3 \to K^2$ is given by 
$$ \pi_2 (g_0, g_1, g_2) = (g_0 g_1^{-1}, g_1 g_2^{-1}) 
~~~\mbox{ (\cite{J2}, (3.9)) }. $$
We define a connection $\conn{2} $ on the total space of this bundle by
$$ \conn{2}  =  \sum_{i = 0}^2 t_i \thetsimp{i}
~~ \in \Omega^1 (\simp^2 \times K^3) \otimes \liek,  $$ 
where $\thetsimp{i}$ $  \in \Omega^1(K^3) \otimes \liek$ is the 
Maurer-Cartan form on the $i$-th copy of $K$.
The curvature 
$$
F_{\conn{2}} \in \Omega^2 (\simp^2 \times K^3) \otimes \liek $$
of the bundle is 
\beq \label{curv}
F_{\conn{2}} = \sum_id ( t_i \thetsimp{i}) + [\conn{2}, \conn{2}] . \eeq
We use this connection and curvature and the Chern-Weil theory of 
equivariant characteristic classes (see for instance 
Chapter 7 of \cite{BGV}) to define an equivariant form on the total  space
$\simp^2 \times K^3$  of the bundle, which represents the 
equivariant characteristic class associated to $\tau_r$ in equivariant
cohomology. We then integrate this equivariant form 
over the simplex $\simp^2.$  Finally, we may pull this form back to
the
base space $K^2$ via a section 
$\sigma_2: K^2 \to K^3  $ 
given by 
$$ \sigma_2(k_1, k_2) = (k_1 k_2, k_2, 1)
~~\mbox{(\cite{J2}, (4.3)) } $$

Explicitly, we make the following definition:
\begin{definition}
Let 
$\Phi_2^K (\tau_r) = \sigma_2^* \bar{\Phi}_2^K (\tau_r) $
$\in \Omega^{2r-2}_K (K \times K  ) $ 
(see \cite{J2}, above (4.3)) 
where the section $\sigma_2$ was defined above, and 
\beq \label{8.9} \bar{\Phi_2}^K (\tau_r) = \int_{\simp^2} 
\tau_r (F_{\theta(t)} + \mom(\theta(t) ) ) . \eeq
\end{definition}

Let 
$\simp^1 = 
 \{ (t_0, t_1) \in [0,1]^2: 
t_0 + t_1 = 1 \} $ $ \cong [0,1]$ be
the standard 1-simplex.
We shall perform a similar construction using a principal $K$-bundle 
$$\simp^1 \times K^2 \stackrel{\pi_1}{\longrightarrow} \simp^1 \times
K.$$
The bundle projection $\pi_1: K^2 \to K$ is defined by 
$$\pi_1(g_0, g_1) = g_0 g_1^{-1}. $$
A section 
$\sigma_1: K \to K^2 $ of the bundle is given by 
$\sigma_1(k) = (k,1). $ 

On the total space $\simp^1 \times K^2$ we define  a connection
$$\conn{1}  = \sum_{i = 0}^1  t_i \thetsimp{i}  \in 
\Omega^1 (\simp^1 \times K^2) \otimes \liek, $$
where $\thetsimp{i}$ $ \in \Omega^1(K^2) \otimes \liek$ is the 
Maurer-Cartan form on the $i$-th copy of $K$.
The definition of the  curvature 
$$F_{\conn{1}}  \in \Omega^2 (\simp^1 \times K^2) \otimes \liek$$ is
similar to ({\ref{curv}}). 
As before, 
we evaluate the invariant polynomial $\tau_r$ on the 
equivariant curvature and integrate over the simplex $\simp^1$ to get 
an equivariant  form over $K \times K$, and finally we pull this form back to 
$K$ using the section $\sigma_1$: explicitly, we make the following
\begin{definition}
We define $$\Phi_1^K(\tau_r) = \sigma_1^* \bar{\Phi}_1^K(\tau_r) 
\in \Omega^{2r-1}(K),$$ where
\beq \label{9.072}
 \bar{\Phi}_1^K(\tau_r) = \int_{\simp^1} \tau_r(F_{\theta(t)} + 
\mom(\theta(t) ) )  \in \Omega^*_K (K \times K). \eeq
\end{definition}

\begin{definition} {\bf (Equivariant chain homotopy)} 
We define a chain homotopy 
$$I_K: \Omega^{* +1}_K (\liek) \to \Omega^{* }_K (\liek) $$ 
as follows:
when  $v \in \liek $, we have
\beq 
\label{9.2a} (I_K \beta)_v = \int_0^1 F_t^* (\iota_{\bar{v} }  \beta) dt, 
\eeq
where $F_t: \liek \to \liek$ is multiplication by $t$ and $\bar{v}$ is 
the vector field on $\liek$ which takes the constant value $v$. 
\end{definition}

\begin{lemma} \label{l9.1} { \bf(\cite{J2}, Theorem 8.1) } 
The equivariant cohomology class of the equivariant differential form
$$\tfr = \proj_1^* \tfr_1 + \proj_2^* \tfr_2 $$
is a lift of $f_r \in H^{2r-2} (\mnd)$ to 
$H^{2r-2}_K (M(c))$. Here, 
the maps $\proj_1$ and $\proj_2$ are the 
projection maps from $M(c)$ to 
$K^{2g}$ and $\liek$ defined at (\ref{4.4}).
Also, from (\cite{J2}, (7.13)), 
\beq \label{9.3} \tfr_1   = \Bigl ( 
\sum_{j=1}^g (- {\rm ev}_{\g_{j}^1} \times 
{\rm  ev}_{x_{j}}  + 
{\rm ev}_{\g_{j+g}^0} \times {\rm ev}_{x_{j+g}} \Bigr )^*  \Phi_2^K(\tau_r)
 (X) + \eeq
$$\Bigl ( 
\sum_{j=1}^g (- {\rm ev}_{\g_{j+g}^1} \times 
{\rm  ev}_{x_{j+g}}  + 
{\rm ev}_{\g_{j}^0} \times {\rm ev}_{x_{j}} \Bigr )^*  \Phi_2^K(\tau_r)
 (X) \in \Omega^*_K (K^{2g}) $$
and
\beq \label{9.4} 
\tfr_2 = - I_K (\epc^* \Phi_1^K (\tau_r)) \in \Omega^*_K(\liek) \eeq 
where $\g_j^\a$ (for $\a = 0,1$ and 
$j = 1, \dots, 2g$) are certain elements of $\free^{2g}$ (the free 
group on $2g$ generators $x_1, \dots, x_{2g}$, 
as in Section 4),
whose definition is given in (7.12) of \cite{J2}, and for any 
$z \in \free^{2g}$,
$ {\rm ev}_z: K^{2g} \to K$ denotes the evaluation map on $z$.
Here, $\epc: \liek \to K$ is defined by $\epc(\Lambda) = 
c \exp \Lambda$ where the central element 
$c = e^{2 \pi i d/n } \diag(1, \dots, 1) $ was defined at (\ref{1.p1}).
\end{lemma}

By (\ref{9.072}) we have
\beq \label{9.005} \bar{\Phi_1}^K (\tau_r)(-X) = 
\int_{t \in [0,1]} 
\tau_r\Bigl (dt (\thetsimp{0} - \thetsimp{1}) + t d\thetsimp{0} + 
(1-t) d \thetsimp{1} + \half [t \thetsimp{0} + (1-t) \thetsimp{1}, 
t \thetsimp{0} + (1-t) \thetsimp{1}] +  \eeq
$$ t \Ad(g_0^{-1} ) X  +  (1-t) \Ad(g_1^{-1})X  \Bigr ). $$
Now $$  \bar{\Phi}_1^K (\tau_r)|_{T \times T}(-X) = 
\int_{t \in [0,1]}  \tau_r(dt (\thetsimp{0} - \thetsimp{1}) + X) $$
since $$d \thetsimp{i}+ \half [\thetsimp{i}, \thetsimp{i}] = 0 $$ and 
the restrictions of $[\thetsimp{i}, \thetsimp{i}] $ to $T$ vanish.
Further $$\sigma_1^* \bar{\Phi}_1^K (\tau_r)|_T (-X) = 
\int_{t \in [0,1]}   \tau_r(dt \theta + X), $$
where $\theta$ is the Maurer-Cartan form on $T$.
If $\tau_r (Z_1, \dots, Z_{n-1}) = \sum_I (\tau_r)_I 
Z^I $ where $I = (i_1, \dots,
 i_{n-1})$ 
is a multi-index and $Z^I = Z_1^{i_1} \dots Z_{n-1}^{i_{n-1}} $ 
(in terms of a coordinate system $\{Z_a = \inpr{\hu{a},X}, \:
a = 1, \dots, n-1 \} $
 on $\liet$, specified by an oriented orthonormal 
 basis $\hu{a}$ for $\liet $ for which 
$\theta_a, a = 1, \dots, n-1$ are the corresponding components of the 
Maurer-Cartan form $\theta \in \Omega^1 (T) \otimes \liet$),
 then 
we have 
$$   {\Phi_1^K} (\tau_r)|_{T}(-X) = 
\sum_I \int_{t \in [0,1]}  (\tau_r)_I (dt \theta_1 + Z_1)^{i_1}
\dots (dt \theta_{n-1}+ Z_{n-1})^{i_{n-1}} $$
\beq 
\label{9.7} = \sum_{a= 1}^{n-1} 
\theta_a  \partial \tau_r /\partial Z_a .  \eeq

\begin{lemma} \label{l9.3} We have for $\L \in \liet$  
(in terms of the Maurer-Cartan form 
$\theta \in \Omega^1(T) \otimes \liet$) that 
$$I_K (\epc^* \theta)_\Lambda = \Lambda. $$
\end{lemma}
\Proof We have 
$$ I_K (\epc^* \theta)_\Lambda = \int_0^1 F_t^* (\epc^* \theta
(\bar{\Lambda} ) )dt 
 = \L $$ since $\epc^* \theta(\bar{\L} ) : \liet \to \RR $ is the 
function with constant value $\L$. \hfill $\square$ 

Let 
$q \in S(\lieks)^K$ be an 
invariant polynomial which is given in terms of the elementary 
symmetric polynomials $\tau_j$ by 
\beq \label{9.0077}
q(X) = \tau_2(X) + \sum_{r = 3}^n \delta_r \tau_r(X). \eeq
The associated
element $\tfq $ of $\hk(M(c))$ 
is defined by 
\beq \label{9.001} \tfq = 
\tf_2 + 
 \sum_{r = 3}^n \delta_r \tf_r. \eeq
Here, the 
 $\delta_r$ are formal 
  nilpotent parameters: we 
expand
 $\exp \tfq$ as a formal power series in the 
$\delta_r$.  We can alternatively
regard the $\delta_r$ as real parameters and $\exp \tfq$
 as a formal equivariant cohomology class: the integral 
$$ \int_{\mnd} \abk (\exp \tfq) $$ and the 
integral appearing in 
 (\ref{9.2}) are well defined and are polynomial functions of the 
$\delta_j$, since $\int_{\mnd} \abk(\eta) = 0 $ unless
$2  {\rm deg}(\eta) = \dim \mnd$. 

Note that by Lemma \ref{l9.1} we can write
$\tfq(X) = \proj_1^* \tfq(X)_1 
+ \proj_2^* \tfq(X)_2$ where 
$\tfq(X)_1  \in \Omega^*_K(K^{2g} ) $
and 
$\tfq(X)_2  \in \Omega^*_K(\liek). $

Then  we have 
\begin{lemma} \label{l9.3a} For $X \in \liet$, the  restriction of
$\tfq(-X)_2 $ to $\mu^{-1}(\liet)$ is given at 
$(h_1, \dots, h_{2g}, \Lambda) $ $ \in \mu^{-1} (\liet) \subset 
K^{2g} \times \liek$ by 
$$  \tfq(-X)_2|_{\mu^{-1}(\liet)} (h_1, \dots, h_{2g}, \Lambda) 
 = - (d \qsgn)_X (\Lambda). $$
\end{lemma}
\Proof 
$$\tfq(-X)_2|_{\mu^{-1}(\liet)} 
= - I_K \epmb^* \Phi_1^K (q)(-X) $$
$$ =  -  I_K (\sum_{a = 1}^{n-1} 
\theta^a\partial \qsgn /\partial Z^a ) ~~  \mbox{by (\ref{9.7}) } $$
$$ = - (d \qsgn)_X (\Lambda) ~~\mbox{by Lemma 
\ref{l9.3}}. $$ 
\hfill $\square$

\begin{lemma} \label{l9.2} 
Assume that $X \in \liet$. Let $\Lambda = 
\sum_{a=1}^{n-1} m_a \he{a} \in \intlat$  for $m_a \in \ZZ$ 
(where the simple roots $\he{a} $ were defined in 
(\ref{6.1}))\footnote{Note that the $\he{a}$ are a basis of $\liet$, but not
an orthonormal basis.} and let 
$s_\L$ denote the homeomorphism of $\mtc$ given by Lemma 
\ref{l4.3}. Then we have that  on $\mtc$
\beq \label{9.007} 
s_\L^* \tfq(-X) = \tfq(-X) -  (d\qsgn)_X (\L),\eeq 
or equivalently 
$$
s_\L^* \tfq(X) = \tfq(X) + (d\qought)_{X} (\L),$$
where we have introduced the notation 
$$\qought (X) = q( -X). $$
\end{lemma}

\noindent{\em Remark:}
This result generalizes (\ref{8.4}).

\noindent{\em Proof of Lemma \ref{l9.2}:}
 Since $\tfq(X) = \proj_1^* \tfq(X)_1 
+ \proj_2^* \tfq(X)_2$, we need to prove the formula for 
$s_\L^* \proj_2^* \tfq(X)_2 $ where  $\tfq(X)_2 = - I_K 
\epmb^* \Phi_1^K (q)$ for $\Phi_1^K (q) \in 
\Omega^*_K(K). $  Lemma \ref{l9.2} then follows 
from Lemma \ref{l9.3a}.  \hfill $\square$ 

\begin{theorem} \label{t9.4} Suppose $\eta$ is a polynomial in 
the $\tar$ and $\tbrj$. Let $q \in S(\lieks)^K$.
Then for any 
$X \in \liet$ we have 
$$ \int_{N_n(V) } \abk \Biggl (\eta e^{\tfq} 
\Bigl ( e^{ (d \qought)_X (\he{1})}  - 1 \Bigr ) \alpha  \Biggr ) = 
-   \sum_{F\in\calf: -||\he{1}||^2  < \inpr{\he{1},\mu(F)} < 0}
\res_{Y_{1} = 0 }  \int_F \frac{\phicom ( \eta  e^{\tfq}  \alpha ) }
{e_{ F}} . $$
Here, we sum over the components $F $ of the fixed point set of 
$\hat{T_1} $ in $P^{-1}(V) \cap \mu^{-1} (\hat{\liet}_1)/T_{n-1} $;
 the notation is as in the 
statement of 
Lemma \ref{old5.17}. The notation 
$\qought$ was introduced in the statement of 
Lemma \ref{l9.2}. 
We have defined the map $\phicom$ 
 in Proposition \ref{stages}, and after (\ref{snow}). 
\end{theorem}
\Proof This follows from the same proof as for Lemma \ref{l6.4},
replacing (\ref{8.4}) by its generalization Lemma \ref{l9.2}. \hfill $\square$

We aim to prove the following result by induction:
\begin{theorem} \label{t9.5} {\bf (a)} For the particular
$q$ defined in (\ref{9.0077}), we have 
$$ \int_{N_n(c) } \abk 
 ( e^{\tfq} \nusym_n \eta ) 
=\frac{(-1)^{n_+(g-1)} }{n!}\sum_{w \in W_{n-1} }
   \res_{Y_{1} = 0 } \dots \res_{Y_{n-1} = 0 }  
\frac { \int_{T^{2g}\times  \{ -\tildarg{w\tc}  \}  )}\Bigl ( e^{ \tfq(X) }   \eta(X)
\Bigr )  }{\nusym(X)^{2g-2} 
\prod_{j = 1}^{n-1} (\exp -B(-X)_{j} - 1) 
 }, $$
where $\eta$ is a polynomial in the $\tarnox$ and $\tbrnox$ and 
$B(X)_j = -(d\qsgn)_X (\he{j}) $. 
Here we have used the fixed invariant inner product
on $\liek$ to identify $d\qsgn_X: \liet \to \RR$ with 
an element of $\liet$ and thus define the  map $B: \liet \to \liet.$ 
The notation $\tildarg{\gamma}$ was introduced in Definition 
\ref{bracedef}.
\end{theorem} 
Substitution of $-X$ for $X$ on the right hand side of the
equation in Theorem \ref{t9.5} (a) gives the equivalent formulation

\noindent{\bf Theorem \ref{t9.5} (b)}  {\em In the notation
of Theorem \ref{t9.5} (a) we have 
$$ \int_{N_n(c) } \abk 
 ( e^{\tfq} \nusym_n \eta ) 
=\frac{(-1)^{n_+(g-1)}  }{n!} \sum_{w \in W_{n-1} } 
  \res_{Y_{1} = 0 } \dots \res_{Y_{n-1} = 0 }  
\frac { \int_{T^{2g}\times \{ - \tildarg{w\tc}  \} }
\Bigl (  e^{\tfq(-X)}  \eta(-X) \Bigr ) }{\nusym(X)^{2g-2} 
\prod_{j = 1}^{n-1} (1 - \exp -B(X)_{j} )
 }. $$
}
Finally we may use 
 Lemma 
\ref{l9.8} and Lemma \ref{l9.9'} (a)  where the restrictions 
to $T^{2g}$ of the equivariant cohomology classes $\tfr$ and $\tbrj$
are expressed in terms of the basis
 $\zeta_a^j$ for $H^1 (T^{2g})$  (for $a = 1, \dots, n-1$ and 
$j = 1, \dots, 2g$). We also use Lemma \ref{l9.03}, where
the symplectic volume of $T^{2g}$ is calculated. 
These lemmas enables us to compute 
$\int_{T^{2g}} e^{ \tfq(-X) } \eta(-X)$ and rephrase
Theorem \ref{t9.5} (b) as follows.
(Here we have also reformulated the left hand side of Theorem \ref{t9.5}
(b) in terms of the pairings on $\mnd$, using Lemma \ref{l3}.)

\begin{theorem} \label{t9.6} In the notation of Theorem \ref{t9.5} 
we have 

\noindent{\bf (a)} 

$$
\int_{\mnd} \exp (f_2 + \delta_3 f_3 + \dots + \delta_n f_n) 
\prod_{r = 2}^n a_r^{m_r} \prod_{k_r = 1}^{2g} 
(b_r^{k_r})^{p_{r,k_r} }  = $$
\beq \frac{ (-1)^{n_+(g-1)} }{n!} \sum_{w \in W_{n-1} } 
   \res_{Y_{1} = 0 } \dots \res_{Y_{n-1} = 0 }  
 \Biggl ( \frac {
e^{ d \qsgn_X (\tildarg{w\tc} ) }  \Bigl ( \prod_{r = 2}^n 
\tau_r(X)^{m_r}   \Bigr ) 
 }{\nusym(X)^{2g-2} 
\prod_{j = 1}^{n-1} (1 - \exp -B(X)_{j} ) 
 } \times  \eeq
$$
 \int_{T^{2g} } 
\exp \Bigl \{ -  \sum_{a,b}  \sum_{j = 1}^g \zeta_a^j \zeta_b^{j + g}  
 \partial^2 \qsgn_X(\hu{a}, \hu{b})  \Bigr \} 
\prod_{r = 2}^n \prod_{k_r = 1}^{2g}
\Bigl  (    \sum_{a = 1}^{n-1} (d\tau_r)_X (\hu{a}) \gen_a^{k_r}
\Bigr )^{p_{r,k_r}}  \Biggr ).  $$ 

\noindent{\bf (b)} In particular we have that 
$$
\int_{\mnd} \exp (f_2 + \delta_3 f_3 + \dots + \delta_n f_n) 
\prod_{r = 2}^n a_r^{m_r}  = $$
\beq (-1)^{n_+(g-1)}
  \frac{n^{g}}{n!} \sum_{w \in W_{n-1} } 
  \res_{Y_{1} = 0 } \dots \res_{Y_{n-1} = 0 }  
\frac {
e^{  d \qsgn_X (\tildarg{w\tc} )}   \prod_{r = 2}^n 
\tau_r(X)^{m_r}  (\det H_\liet (X))^g  }{\nusym(X)^{2g-2} 
\prod_{j = 1}^{n-1} (1 - \exp -  B(X)_{j}) 
 }. \eeq

\end{theorem}

\begin{rem}
In the preceding Theorem, we have used the following notation.
The $a_r, f_r$ and $b_r^j$ (for $r = 2, \dots, n$ and
$j = 1, \dots, 2g$) are generators of the cohomology ring, introduced 
in Section 2. The $\tau_r$ are the elementary symmetric 
polynomials, and the $\delta_r$ are formal nilpotent parameters
which were introduced in (\ref{9.0077}). The polynomial $q = 
\tau_2 + \sum_{r = 3}^n \delta_r \tau_r$ was introduced in 
(\ref{9.0077}). Its derivative $dq_X: \liet \to \RR$ is
identified with an element of $\liet$ via the inner product
$\inpr{\cdot,\cdot}$ 
on $\liet$, and hence $dq: \liet \to \liet^*$ is identified with 
a map $B: \liet \to \liet$ (see the statement of Theorem 
\ref{t9.5}(a)). 
If $\gamma \in \liet$,
the notation $\tildarg{\gamma}$ was introduced
in Definition \ref{bracedef}: it is the unique element 
 in the fundamental domain defined by
the simple roots for the translation action
on $\liet$ of the integer lattice which is 
 equivalent to $\gamma$ under translation by the integer
lattice.
The $\zeta_a^j$ are derived from the components of the
Maurer-Cartan form $\theta \in \Omega^1(T) \otimes \liet$ in terms of an
orthonormal basis $\{ \hu{a}, a = 1, \dots, n -1\}$ of 
$\liet$: they have been identified with 
a basis of $H^1 (T^{2g})$. (See 
Definition \ref{zedef} below.) Finally $\det H_\liet(X)$ is the 
determinant of the Hessian 
of  $q: \liet \to \RR,$  in terms of the coordinates
on $\liet$ given by the orthonormal 
basis $\{\hu{a} \}$: it is independent of the choice of orthonormal
basis. 

We note that in Theorem \ref{t9.6} 
the orthonormal basis introduced above could
be replaced by a general basis, 
provided one defines the $\zeta_a^j$ using that basis, and
multiplies the Hessian by a factor due to the change of basis: see 
 Remark \ref{r10.1} below.
\end{rem}

\begin{rem}  We can replace $f_2$  by any nonzero constant 
scalar multiple $\epsilon f_2$ provided we replace the 
polynomial $q$ by $q^\epsilon$ where
$$q^\epsilon(X) = \epsilon \tau_2(X) + \delta_3 \tau_3(X) + \dots
+ \delta_n\tau_n(X) $$
(cf. Remark \ref{rem8.3}  (b)). 
\end{rem}

In order to prove Theorem \ref{t9.5} and hence Theorem
\ref{t9.6} we follow the  proof of Theorem \ref{mainab} using

\begin{lemma} \label{l9.7} Suppose $\eta$ is a polynomial 
in the $\tar$ and $\tbrj$. Then for any $X \in \liet$ and 
$q \in S(\lieks)^K$ chosen as in (\ref{9.0077}), we have
\beq \label{9.02} 
\int_{N_n(V)} \abk \Bigl (\eta
 \nusym  e^{\tfq }  \alpha \Bigr )  = 
-   \sum_{F\in\calf: -||\he{1}||^2  < \inpr{\he{1},\mu(F)} < 0}
\res_{Y_{1} = 0 }  \int_F \frac{\phicom (
 \nusym \eta(X) e^{\tfq} \alpha  ) }
{e_{F}  
\Bigl ( e^{  (d\qought)_X (\he{1})} - 1 \Bigr ) },\eeq
where the notation is as in 
the statement of Lemma \ref{old5.17}.
\end{lemma}
\Proof This follows from Theorem \ref{t9.4} by replacing 
$\eta $ by 
\beq \label{e9.16} \frac{\eta \nusym  }
{\Bigl ( e^{(d\qought)_X (\he{1})} - 1 \Bigr ) } = 
\frac{\eta \bar{\nusym } Y_{1} }
{\Bigl ( e^{  (d\qought)_X (\he{1}) } - 1 \Bigr ) }, \eeq
where we have defined 
$\bar{\nusym} = \nusym/Y_1. $ 
Notice that $(d\qought)_X (\he{1})$ is divisible by $Y_1$: to see this,
we observe that if we define 
the generating functional 
$$P(X_1, \dots, X_n) = \prod_{j = 1}^n 
(1 + t X_j) = \sum_{r = 0}^n \tau_r(X_1, \dots, X_n) t^r $$
(where the $\tau_r$ are the elementary symmetric polynomials) then
$$ dP = \left ( (1+ tX_1) tdX_2 + (1 + t X_2) tdX_1 \right) 
\prod_{j = 3}^n (1 + t X_j) + {\cal P}$$
where ${\cal P}$ is a collection of terms involving 
$dX_3, \dots, dX_n$. Evaluating $dP$ on 
$\he{1} = (1,-1, 0, \dots, 0)$ we thus obtain
$$ t^2 (-Y_1) \prod_{j=3}^n (1 + tX_j) = 
\sum_{r = 0 }^n t^r  (d\tau_r)_X(\he{1}). $$
It follows that the $(d\tau_r)_X(\he{1})$ (and hence
$(d \qought)_X (\he{1}) $) are  divisible by 
$Y_1$.
Thus $-(d\qought)_X(\he{1}) = - Y_{1}(1   + {\nu})$ where
${\nu} \in H^*_{T} $ has degree at least $ 1$, so we have 
$$\Bigl ( e^{  (d \qought)_X (\he{1})} - 1 \Bigr )  = 
-   Y_{1}(1   - \tilde{\nu}) $$ where 
$\tilde{\nu} = \sum_{j \ge 1} \tilde{\nu}_j $ is a formal sum of 
classes $\tilde{\nu}_j$ with degree 
at least $1$ in (a completion of) $H^*_T$. Then the 
expression 
$$ \abk \Biggl ( \frac{\eta \nusym  }
{\Bigl ( e^{  (d\qought)_X (\he{1})} - 1 \Bigr ) }  \Biggr ) $$ 
(which appears on the left hand side of the equation in 
Theorem \ref{t9.4})
 is well defined. On the right hand side, we may replace
$$\frac{  Y_{1} }
{\Bigl ( e^{+ (d\qought)_X (\he{1})} - 1 \Bigr ) }$$
by $- (1 - \tilde{\nu})^{-1} =  -\sum_{s \ge 0 } 
 {\tilde{\nu}}^s. $  \hfill $\square$ 

We now use Lemma \ref{l9.7} to prove Theorem \ref{t9.5} (a)
by induction on $n$. The proof follows the outline of the 
proof of Theorem \ref{mainab} 
when $q = q_2$, with the following modifications:
\begin{enumerate}
\item $e^\bom$ is replaced by $e^{\tfq}$
(and $e^\omega$ replaced by $e^{\tfq - \tilde{f_2}} e^\omega
$).
\item$(d\qought)_X (\he{1})$ replaces $-Y_1$, 
so 
 $e^{ (d\qought)_X (\he{1}) }
  - 1 $ replaces $e^{-\normcon Y_1} - 1$.
\item In particular, $-(d \qought)_X (\he{1})$ $ = -B (-X)_1 $ replaces
$  Y_1$ in the identity 
$$ \frac{e^{- \delta_I Y_1} }{1 - e^{-Y_1} } = 
 \frac{e^{(1- \delta_I )Y_1} }{e^{Y_1}  - 1 } $$
which is used in the proof of Proposition \ref{beg}.
\end{enumerate}
We also use the elementary fact that 
$(d \qought)_X = - (dq)_{-X} = B (-X). $

\renorm\
\section{Witten's formulas for general intersection pairings}

In this section  we state and prove
 Witten's formulas (Propositions \ref{p9.1} and
\ref{p9.2} below; cf. \cite{tdgr}, Section 5, 
in particular the calculations
(5.11)-(5.20)),
which enabled him to calculate general intersection pairings
in terms of those of the form 
$$\int_{\mnd} \prod_r a_r^{m_r} e^{f_2} . $$
We shall prove these formulas starting from our explicit formulas for
the 
general intersection pairings (see Theorem \ref{t9.5}). 

Some of the  notation in the statement of Propositions \ref{p9.1} and 
\ref{p9.2} was introduced at the beginning of Section 9.
The invariant polynomial $q$ was defined by (\ref{9.0077}). 
Using the 
invariant metric on $\liek$, the map $-d\qsgn:
\liek \to \lieks$ may be regarded as a map 
$B = B^{(2)}  + 
\sum_{r \ge 3} \delta_r B^{(r)} :\liek \to \liek$, 
where we have written $B^{(r)} =    -d \tau_r:
  \liek \to \liek$;
we find $B^{(2)} = -d \tau_2 = id: \liek \to \liek$. (Note that 
we have put $\tau_2(X) = -\half \inpr{X,X}$ in terms 
of the inner product 
$\inpr{\cdot,\cdot} $ defined at (\ref{1.02}).)
The other maps $B^{(r)}$ are not
linear.

The Hessian of $-\qsgn$ is 
$H$; it is a function from $\liek$ to symmetric bilinear forms 
on $\liek$. 
If $k, l$ run over an orthonormal basis $\{\hat{v}_k \} $ of 
$\liek$ then the Hessian at $X$ is the matrix 
\beq \label{hessdef}
H(X)_{kl} = -  (\partial^2 \qsgn)_X (\hat{v}_k,
\hat{v}_l) . \eeq

\begin{rem} \label{r10.1}
In most places in Sections 9 and 10, 
 the orthonormal basis
 $\{\hu{a}\} $  for $\liet$ may be replaced by any basis 
for $\liet$ (including the basis $\{ \he{a}, a = 1, \dots, n-1 \}$,
which is of course not orthonormal), and 
similarly for the orthonormal basis $\{ \hat{v}_l \} $ for 
$\liek$. However it is more convenient
to define the determinant of the Hessian (given  in (\ref{hessdef}))
 in terms of an {\em orthonormal} basis, since
one must otherwise include a normalization factor proportional
to the square of the determinant of a matrix whose columns 
are the basis elements.
The second place where it is useful to introduce  an orthonormal 
basis is in the definition of the symplectic form in terms
of the generators $\zeta_a^j$ for the cohomology of $T^{2g}$: the
symplectic form is defined using the inner product
$\inpr{\cdot, \cdot} $ on $\liet$, and the
formula (Lemma \ref{l9.7a}) for the  restriction of the
symplectic form to $T^{2g}$ is cleaner  in terms of an 
orthonormal basis.

For these reasons we have chosen to use an orthonormal basis 
for $\liet$ throughout
Sections 9 and 10, although in many specific instances this basis
may be replaced by a general basis.
In particular in the statement of our main theorem 
Theorem \ref{t9.6}, it is easy to check that the orthonormal basis 
may be replaced by a general basis, provided that the 
$\zeta_a^j$ are also defined using this basis,
and that the Hessian is multiplied by the
appropriate
factor.
\end{rem}

\nc{\winv}{{ { B^{-1} }}}
\nc{\bwinv}{ (\winv) } 

 We assume the 
$\delta_r$ are formal {\em nilpotent}
 parameters: then the invertibility of  $B $ 
is guaranteed. 
We write $ B^{-1}: \liek \to \liek$ as the inverse of $B$. 
(If the $\delta_r$ are nilpotent, the inverse of $B$ may be written
as a formal power series in the $\delta_r$.)

\begin{prop} 
\label{p9.1} 
 For any invariant polynomial 
$\tau \in S(\lieks)^K$, the integral 
\beq \label{10.1a}  \int_{\mnd} \abk ( \tau(-X) \exp \tfq   )  \eeq
is equal to the integral 
\beq \label{9.1}
 \int_{\mnd} \abk \Biggl (\tau(\winv(-X)   ) \Bigl ( \det H(\winv(-X)) 
\Bigr )^{g-1}   \Biggr ) \exp f_2 \eeq
which is of the form that  may be calculated by 
Theorem \ref{mainab}.
\end{prop}
\begin{prop} \label{p9.2} 
Let $\tau \in S(\lieks)^K$ be an invariant polynomial, so that
$$\tau = \sum_{m_2, \dots, m_n} c_{m_2, \dots, m_n} \prod_{r = 2}^n 
\tarnox^{m_r}$$
 is a polynomial in the $\tilde{a}_2, \dots, \tilde{a}_n$. 
Let $s_r^j$ be real parameters
(for $r = 2, \dots, n$ and $j = 1, \dots, 2g$). Then we have 
\beq \label{9.2} \int_{\mnd} \abk \Bigl ( \tau(X) 
\exp (\sum_{r = 2}^n 
\sum_{j=1}^{2g} \srj
 \tbrj ) \exp \tfq \Bigr )
 = \int_{\mnd} \abk\Bigl 
(\tau(X) \exp \hattau(X) \exp \tfq
\Bigr ). \eeq
Here, the invariant polynomial $\hattau$ on 
$\liek$ is defined (for $X \in \liet$) by 
$$\hattau(-X)  = - \sum_{a,b=1}^{n-1}
\sum_{r,s  = 2}^{n}  \sum_{j=1}^{g}   
   s^j_r s^{j+g}_s (d\tau_r)_X(\hu{a}) 
 (d\tau_s)_X(\hu{b})  (\partial^2 \qsgn )^{-1}_{ab}, $$
where  $\{ \hu{a}: a = 1, \dots, n-1 \} $ denotes an oriented orthonormal 
basis of $\liet$: see 
(\ref{hattaudef}) for the definition.
\end{prop}

\begin{rem} Notice that in our conventions on the 
equivariant cohomology differential and the moment, the construction of 
\cite{J2} described at the beginning of Section 9 yields 
$$\tarnox(X) = \tau_r(-X) . $$
Thus $\tau(X) = 
 \sum_{m_2, \dots, m_n} c_{m_2, \dots, m_n} \prod_{r = 2}^n 
\tau_r(-X)^{m_r}. $
\end{rem}

Proposition \ref{p9.1} is proved by comparing Theorem 
\ref{t9.6} (b) (applied to 
(\ref{10.1a})) with Theorem \ref{mainab} (applied to 
(\ref{9.1})). Proposition 
\ref{p9.2} is obtained by applying Theorem \ref{t9.5} (b) to 
both sides of (\ref{9.2}) and examining the restrictions to 
$T^{2g}$ (which are computed in  Lemmas \ref{l9.8} and \ref{l9.10}).

Propositions \ref{p9.1} and \ref{p9.2} enable us to extract
formulas for all pairings, by differentiating the formulas (\ref{9.1})
and (\ref{9.2}) with respect to the parameters $\delta_r$ and 
$\srj$ and then setting these parameters equal to zero.
In fact, for any nonnegative integers $n_r$ 
(for $r \ge 3$) we have
$$ \Biggl (\prod_{r = 3}^n \Bigl ( \frac{\partial}{\partial \delta_r}\Bigr
)^{n_r}
\int_{\mnd} \abk ( \tau(X) \exp \tfq   )  \Biggr )_{\delta_3 = \dots =
  \delta_n = 0 } = 
\int_{\mnd} \prod_{r = 3}^n f_r^{n_r} 
\abk \Bigl ( \tau(X) \Bigr ) \exp f_2 , $$
and likewise for any nonnegative integers
$n_r$ (with $n_2 = 0 $) and any choices of $p_{r, j_r} = 0, 1$ we have
$$ \Biggl ( \prod_{r = 2}^n \Bigl ( \frac{\partial}{\partial \delta_r}
\Bigr)^{n_r}  \prod_{j_r = 1}^{2g}
\Bigl ( \frac{\partial}{\partial s_{r}^{ j_r} }\Bigr
)^{p_{r,j_r} } 
 \int_{\mnd} \abk \Bigl ( \tau(X) 
\exp (\sum_{r = 2}^n 
\sum_{j=1}^{2g} \srj
 \tbrj ) \exp \tfq \Bigr )  \Biggr)_{\delta_r = 0 , 
~s_{r,j} = 0}
 $$
$$
=  \int_{\mnd} \abk \Bigl ( \tau(X) \Bigr )  \exp f_2 
\prod_{r = 2}^n f_r^{n_r} \prod_{j_r = 1}^{2g} \Bigl ( b_{r}^{j_r} 
  \Bigr )^{p_{r,j_r}}
$$
(where the parameters $\delta_r$ and 
$\srj$ on the left hand side run over 
$r = 2, \dots, n$ and $j = 1, \dots, 2g$).

We can use Proposition 
\ref{p9.2} to  give an explicit formula for pairings of the form 
\beq \label{9.31}  \int_{\mnd} \Phi \Bigl ( \prod_{r=2}^n \prod_{k_r = 1}^{2g} 
(\tilde{b}_r^{k_r}(X) )^{p_{r,k_r} } \tau(X)  \Bigr ) e^{f_2}  \eeq
where $p_{r,k_r} = 0 $ or $1$. 
We note that by Proposition \ref{p9.2}  this equals 
\beq \label{9.32}
\prod_{r = 2}^n \prod_{k_r = 1}^{2g} 
\Bigl ( \frac{\partial}{\partial s_r^{k_r} } \Bigr )^{p_{r,k_r} }  
\int_{\mnd} \Phi \Bigl ( \tau(X) \exp \hattau(X)  \Bigr ) e^{f_2} 
\mid_{s_r^j = 0 ~\forall ~ r,j} \eeq 
where 
\beq \label{9.32a} 
\hattau(-X) = 
- \sum_{a, b = 1}^{n-1} 
\sum_{j = 1}^g 
\sum_{r,s = 2}^n 
      s_r^j s_s^{j+g} (d \tau_r)_X (\hu{a} ) 
(d \tau_s)_X (\hu{b} ) (\partial^2 q)^{-1}_{ab} . \eeq
Here $\hu{a}$ are an oriented orthonormal basis of $\liet$. We introduce
$T_{rs}: \liek \to \RR$ given by\footnote{Notice that 
$T_{rs}$ is an invariant polynomial on $\liek$.}
$$ T_{rs} (-X) = 
- \sum_{a, b = 1}^{n-1} 
       (d \tau_r)_X (\hu{a} ) 
(d \tau_s)_X (\hu{b} )(\partial^2 q)^{-1}_{ab}.   $$
Thus we may rewrite (\ref{9.32a}) as 
\beq \label{9.33}
\hattau(X) = 
 \sum_{r , s = 2}^n T_{rs}(X) (\sum_{j = 1}^g s_r^j s_s^{j+g} ). \eeq
We observe that in order for the pairing (\ref{9.31}) to be nonzero, one 
requires $p_{r,j}  = 0 $ or $1$ for all $r$ and $j$ (since the $b_r^j$ 
are of odd degree). Further, in order for the expression 
(\ref{9.32}) to yield a nonzero answer, we require for each 
$j = 1, \dots, g$ that 
$$p_{2,j} + \dots + p_{n,j} = p_{2,j+g} + \dots + p_{n,j+g}  = l_j $$
for some $l_j$. We may then rewrite (\ref{9.32}) as 
\beq \label{9.34} 
\Biggl ( \prod_{j = 1}^g
 \Bigl ( \spart{r_1}{j}\dots \spart{r_{l_j}}{j} \Bigr ) 
\Bigl ( \spart{s_1}{j+g}\dots \spart{s_{l_j}}{j+g} \Bigr )
\int_{\mnd} \Phi \Bigl ( \tau(X) \exp \hattau(X)  \Bigr ) e^{f_2} 
\Biggr )_{s_r^j = 0 ~\forall ~ r,j}. \eeq 
Because $\hattau$ is quadratic in the $s_r^j$ and we are setting all the 
$s_r^j$ to zero in the end, for each $j$ we may represent the symbols
$\spart{r}{j} $ and $\spart{r}{j+g} $
as  1-valent vertices (labelled by $r$) 
in a bipartite graph: there must be exactly one edge coming out of each of 
these vertices, and these edges must connect  the symbol 
$\spart{r}{j} $ with a symbol $\spart{s}{j+g} $ for some $s$. Such
bipartite graphs of course correspond to permutations
$\sigma_j$ of 
$\{ 1, \dots, l_j \}$. 

It follows from (\ref{9.33}) that 
$$ \spart{r}{j} \spart{s}{j+g} \hattau(X) = T_{rs} (X)  ~\mbox{for any $j$} $$
so that 
\beq \label{9.35} 
\prod_{j = 1}^g 
\Bigl ( \spart{r_1}{j}\dots \spart{r_{l_j}  }{j} \Bigr ) 
\Bigl ( \spart{s_1}{j+g}\dots \spart{s_{l_j}}{j+g} \Bigr )
\int_{\mnd} \Phi \Bigl ( \tau(X) \exp \hattau(X)  \Bigr ) e^{f_2}   \eeq
$$ = \int_{\mnd} 
\Phi \Biggl ( \prod_{j = 1}^g 
\sum_{ \sigma_j} T_{r_1 s_{\sigma_j(1)} } (X)  \dots 
T_{r_{l_j} s_{\sigma_j(l_j) } } (X)  \tau(X)   \Biggr) e^{f_2} $$ 
where we sum over all permutations $\sigma_j$ of $ \{ 1, \dots, l_j \}$. 
Hence we obtain by Remark \ref{rem8.3} (a) and 
Lemma \ref{l3}
\begin{theorem}
\beq \label{9.36} 
\int_{\mnd}  
\prod_{j = 1}^g b_{r_1}^j \dots b_{r_{l_j} }^j 
b_{s_1}^{j+g}  \dots b_{s_{l_j} }^{j +g}  \Phi (\tau(X) ) e^{f_2}  \eeq
$$ 
= \int_{\mnd} 
\Phi \Biggl ( \prod_{j = 1}^g 
\sum_{ \sigma_j} T_{r_1 s_{\sigma_j(1)} }(X)  \dots 
T_{r_{l_j} s_{\sigma_j(l_j) } } (X)  \tau(X)   \Biggr) e^{f_2} $$
which equals 
$(-1)^{n_+(g-1) }  \frac{n^{g}}{n!}$ times the iterated residue
$$ \sum_{w \in W_{n-1} } 
 \res_{Y_{1} = 0 } \dots 
\res_{Y_{n-1} = 0 } 
\frac{
\prod_{j = 1}^g \sum_{\sigma_j} 
T_{ {r_1} s_{\sigma_j(1) } }  (-X)  \dots
T_{ {r_{l_j} } s_{\sigma_j(l_j) } }  (-X)  \tau(-X) 
e^{ - \inpr{\bracearg{w\tc}, X} } }{\nusym^{2g-2}(X) (1 -  \exp (-Y_1) ) 
\cdots  ( 1 - \exp (-Y_{n-1} )  ) } . $$

\end{theorem}

Let $\hu{a} $ $(a = 1, \dots, n-1)$ denote an oriented orthonormal basis on 
$\liet$.
For $X \in \liet$ 
define coordinates $\zz{a}$ by  $\zz{a} = 
(X, \hu{a}) $ so that 
$X = \sum_a \zz{a} \hu{a}.$
Write the Maurer-Cartan form $
\theta$ on $T$ as $\theta = \sum_a \theta_a \hu{a}; $ then 
the $\theta_a $ form a set of generators of $H^1(T). $  

\begin{definition} \label{zedef}
A set of  generators $\{\zeta^j_a \} $ ($j = 1, \dots, 
2g; $ $a = 1, \dots, n-1$) for  $H^1(T^{2g})$ 
is defined by by specifying that
$\zeta^j_a = \pi_j^* \theta_a$ where $\pi_j: T^{2g} \to T$ is the 
projection onto the $j$'th copy of $T$.
\end{definition}

\begin{lemma} \label{l9.3b} We have 
$ \int_T \theta_1 \wedge \dots \wedge \theta_{n-1} = 
  \vol(T). $ 
Here, $\vol(T)$ is the Riemannian volume of $T = \liet/\intlat$ 
in the metric $\langle \cdot , \cdot \rangle$: in other words it is
given by $(\det E)^\half$  $ = \sqrt{n} $ 
where $E$  is the 
$(n-1) \times (n-1) $ 
matrix (known as the Cartan matrix)
 given by ${E}_{ab} = \inpr{\he{a}, \he{b}} $ in terms 
of the  basis for the integer lattice 
$\intlat \subset \liet$  over $\ZZ$ 
given by the 
simple roots  $\{\he{a} \}$, $a = 1, \dots, n-1$.
\end{lemma}

\begin{lemma} \label{l9.7a} The restriction of $\tfq(X)_1$ to $T^{2g}$ 
is given in terms of the generators $\gen_a^j$ of $H^1(T^{2g})$ by
$$ \tfq(X)_1 \mid_{T^{2g} } =  \half \sum_{a,b} \partial^2 \qsgn (\hu{a}, 
\hu{b}) \sum_{j = 1}^g (-\gen_a^j \gen_b^{j+g}  + \gen_a^{j+g}
\gen_b^j ), $$
where $\{ \hu{a} \}$ are an oriented  orthonormal basis of $\liet$. 
\end{lemma}

\Proof We need to understand the restriction 
of $\tfq$ to $T^{2g}$. As in (\ref{9.3}), we have 
$$
 \tfq(X)_1  = \Bigl ( 
\sum_{j=1}^g (- {\rm ev}_{\g_{j}^1} \times 
{\rm  ev}_{x_{j}}  + 
{\rm ev}_{\g_{j+g}^0} \times {\rm ev}_{x_{j+g}} \Bigr )^*  \Phi_2^K(q)
 (X) + $$
$$\Bigl ( 
\sum_{j=1}^g (- {\rm ev}_{\g_{j+g}^1} \times 
{\rm  ev}_{x_{j+g}}  + 
{\rm ev}_{\g_{j}^0} \times {\rm ev}_{x_{j}} \Bigr )^*  \Phi_2^K(q)
 (X) \in \Omega^*_K (K^{2g}) $$
 where (after 
restricting to $T \times T \times T$) 
\beq \label{9.002}  \bar{\Phi}_2^K (q) |_{T\times T \times T}(-X)  = 
\int_{(t_0, t_1, t_2) \in \simp^2} 
\qsgn ( \sum_{k = 0 }^2 dt_k \thetsimp{k} + X) \in \hht(T\times T
\times T)  \eeq
(\cite{J2}, above (5.6)) 
and $\Phi_2^K(q) = \sigma_2^* \bar{\Phi}_2^K(q) $  where 
$\sigma_2: (g_1, g_2) \mapsto
(g_1 g_2, g_2, 1). $ 
By (\ref{9.002})  we have 
\beq \label{9.003} 
 \Phi_2^K(q)|_{T \times T}(-X) = 
- \half \sum_{a,b} \partial^2 \qsgn(\hu{a}, \hu{b}) \gen_a^1 
\gen_b^2  \in \hht (T \times T).\eeq

For the purposes of evaluation on  $T^{2g}$ the generators 
$\g_j^\tau$ in (\ref{9.3}) reduce to 
$$\g_{j}^0  = \g_{j+g}^1  = 1,~~ \g_{j}^1  = x_{j+g}, 
~~\g_{j+g}^0  = x_{j}, $$
where $x_1, \dots, x_{2g}$ are the chosen generators of $\free^{2g}$. 
So we get from (\ref{9.3}) 
$$ \tfq(X)_1 \mid_{T^{2g} } = \sum_{j=1}^g \Bigl ( 
-{\rm ev}_{x_{j+g}} \times {\rm ev}_{x_j}  + 
{\rm ev}_{x_{j}} \times {\rm ev}_{x_{j+g}  }    \Bigr )^* 
\Phi_2^K(q) (X).  $$

We find that
$$ \tfq(-X)_1|_{T^{2g}} = - \half \sum_{a,b} \partial^2 \qsgn (\hu{a}, 
\hu{b}) \sum_{j = 1}^g (\gen_a^j \gen_b^{j+g}  - \gen_a^{j+g} \gen_b^j
)
= - \sum_{a,b} \partial^2 \qsgn (\hu{a}, 
\hu{b}) \sum_{j = 1}^g \gen_a^j \gen_b^{j+g}.
\square$$

Similarly for $\L \in \liet \subset \liek$ we have 
$$ \tfq(-X)_2(\L) =  - (d\qsgn)_X (\L) $$ (see Lemma \ref{l9.3a}). 

As a result we see immediately that 
\begin{lemma} \label{l9.8} Suppose
 $c \exp \L = 1.$ Then we have 
$$ \int_{T^{2g} \times \{ \Lambda \}  } \exp \tfq (-X) = 
  \int_{T^{2g} \times \{\Lambda \}  }  \exp \tf_{q_2}(-X) (\det \hess(X))^g, $$
$$ = e^{-(d \qsgn)_X (\L) }
\int_{T^{2g} }  (\det \hess(X))^g \exp \omega $$
where $\omega$ is the standard
symplectic form on $T^{2g}$ and the quadratic form
 $\hess(X)  $ is the Hessian of the restriction of 
$-\qsgn$ to $\liet$ (evaluated on an oriented  {\em orthonormal basis}
of $\liet$). In other words, 
$$\hess(X)_{ab}   = -(\partial^2 \qsgn)_X (\hu{a}, \hu{b}) $$
where $\{ \hu{a}: a = 1, \dots, n-1 \} $ is an oriented orthonormal 
basis for $\liet$. 
\end{lemma}

\Proof 
This follows by integrating 
$$ \exp -\sum_{a,b} 
 \partial^2 \qsgn (\hu{a}, 
\hu{b}) \sum_{j = 1}^g (\gen_a^j \gen_b^{j+g}) $$
over $T^{2g}$. (Notice that $\partial^2 \qsgn (\hu{a}, 
\hu{b})$ is symmetric in $a$ and $b$.)  \hfill $\square$ 

\begin{lemma} \label{l9.03}
We have that 
$$ \int_{T^{2g} } \exp \omega = n^g. $$
\end{lemma}

\Proof This follows from Lemmas \ref{l9.3b} and 
\ref{l9.7a}.   \hfill $\square$

In order to prove Proposition \ref{p9.1}, note that by Theorem 
\ref{t9.5}(b), we have that
$$\int_{N_n(c)} \Phi ( e^{\tfq} \nusym_n  \eta) $$ equals 
$\frac{(-1)^{n_+(g-1) } }{n!}$ times the iterated 
residue 
$$ \sum_{w \in W_{n-1} } 
    \res_{Y_{1} = 0 }  \dots 
\res_{Y_{n-1} = 0 } \frac{
   \int_{T^{2g} \times \{ -\tildarg{w\tc}  \}  }
\Bigl ( e^{ \tfq (-X)}   \eta(-X) \Bigr ) 
  }{\nusym(X)^{2g-2}
(1 - e^{ -B( X)_{n-1}} ) \dots (1 - e^{ -B(X)_1 }) }. $$
This applies in particular when $\eta (X) = \tau(-X) $ is a linear combination 
of monomials $\prod_r \tarnox^{m_r}$ in the $\tarnox$ which does not
involve the $\tbrnox$; since 
$\tarnox(X)  =\tau_r(-X),$ it is natural to write
$\eta(-X) = \tau(X). $ For $\eta$ of this form, 
the expression above  equals
$$ \frac{(-1)^{n_+(g-1) } }{n!} \sum_{w \in W_{n-1} } 
 \res_{Y_{1} = 0 }  \dots \res_{Y_{n-1} = 0 }  
\frac{ \int_{T^{2g} \times \{ - \tildarg{w \tc} \} }
 e^{  \tf_{q_2}(-X) }
\Bigl (\det \hess(X) \Bigr )^g  \tau(X) } {\nusym(X)^{2g-2} 
(1 - e^{ -B(X)_{n-1} } ) \dots (1 - e^{ -B(X)_1} ) } $$
by Lemma \ref{l9.8}.

We now replace $X$ by $\winv(X)$ (where the 
 transformation $\winv:$  $\liek \to \liek$ was 
defined above Proposition \ref{p9.1}). This change of variables produces
a Jacobian $\Bigl ( \det \hess(\winv(X)) \Bigr )^{-1} $. Thus we obtain that 
\beq \label{9.004} \int_{N_n(c)}  \abk (\nusym e^{\tfq} \eta) = 
\frac{(-1)^{n_+(g-1) } }{n!} \sum_{w \in W_{n-1} } 
\res_{Y_{1} = 0 }  \dots \res_{Y_{n-1} = 0 }   \Biggl ( 
 \Bigl (\det \hess
(\winv(X)) \Bigr
 )^{g  -1} \times \eeq
$$ \frac{  \int_{T^{2g} \times \{  - \tildarg{w \tc} \} } 
e^{\tf_{q_2}   (-\winv(X)) }
\tau (\winv(X) ) }{  \nusym^{2g-2}(\winv(X)) 
(1 - e^{ -   Y_{n-1}} ) \dots (1 - e^{  -  Y_1})  } \Biggr ). $$ 
Now we have 
\begin{lemma} \label{l9.9} 
$$  \nusym^2 (\winv(X) ) = \nusym^2(X)
 (\det (-\partial^2 \qsgn)_\lietp )^{-1} , $$
where $ (\partial^2 \qsgn)_\lietp $ denotes the restriction of the
symmetric bilinear form 
$ (\partial^2 \qsgn)_X$  on $\liek$ to 
a symmetric bilinear form on $\lietp$, which is
then identified with a linear map from 
$\lietp$ to itself using the fixed invariant inner product.
\end{lemma} 
This lemma will be used in establishing Proposition 
\ref{p9.1} since the Hessian $H$ appearing in that proposition is 
the Hessian of the ($K$-invariant) function $-q: \liek \to \RR$, 
which is block diagonal with one block being the
Hessian $H_{\liet}$  of the restriction of
this function  to $\liet$ and the other block being 
$ -(\partial^2 \qsgn)_\lietp $.  

\Proof 
We  introduce the (orthonormal) basis $X_\g, Y_\g$  for $\lietp$
corresponding to the positive roots $\g$, and a corresponding 
system of coordinates $x_\g, y_\g$ on $\lietp$:  we have
$$ [X_\g, X] = \g(X) Y_\g, ~~~[Y_\g, X] = -\g(X) X_\g.$$ 
We observe that  the map $B$ and its inverse 
$\winv$  on $\liek$ are $K$-equivariant, and map $\liet$ to 
$\liet$ and $\lietp$ to $\lietp$. 
Hence 
\beq (d \winv)_X ([X_\g, X] ) = 
[X_\g, \winv(X) ] \eeq
and
$$\winv \Biggl (  {\rm Ad} \exp (X_\g)  ( X ) \Biggr )  = 
{\rm Ad} \exp (X_\g) (\winv (X) ), $$
and similarly for $Y_\g$. We find 
\beq  \frac{\g(\winv ( X))}{\g(X)} = 
(d\bwinv_{y_\g})_X(Y_\g) 
= (d\bwinv_{x_\g})_X(X_\g)
\eeq
where
$\bwinv_{x_\g}, \bwinv_{y_\g}:$ $\lietp \to \RR$
are the coordinate functions in the directions
$ x_\g$ and  $y_\g$. 
Thus, 
\beq \label{6.10} \cald^{2}(\winv(X)) = \cald^{2}(X) 
 (\det d \bwinv_\perp)  
 \eeq
$$ \onebl \onebl = \cald^{2}(X) (\det -\partial^2 \qsgn)_{\lietp}^{-1} $$
where $d \bwinv_\perp$ is the square matrix of partial 
derivatives of the $\lietp$ components of $\winv$ in 
the directions along $\lietp$.  This completes the proof
of Lemma \ref{l9.9}. \hfill $\square$

Proposition \ref{p9.1}
now follows immediately by using (\ref{9.004}) and Lemma \ref{l9.9}
to express (\ref{10.1a}) as an iterated residue, and observing that
Theorem \ref{mainab} (in the version given by 
Remark \ref{rem8.3} (a)) yields the same iterated residue for (\ref{9.1}). 

Let us  now
consider the proof of Proposition \ref{p9.2}.
For the rest of  this section 
let $a = 1, \dots, n-1$ index an oriented {\em orthonormal} basis 
$\{\hu{a}\}$  
of $\liet$. 
We have 
$\tbrj  = \proj_1^*  \tilde{b}_r^{j,1} $ where
$ \tilde{b}_r^{j,1}  = \ev{x_j}^* \Phi_1^K (\tau_r) $
and $\Phi_1^K (\tau_r) = 
\sigma_1^* \bar{\Phi_1^K }(\tau_r) $ where 
$\bar{\Phi_1^K }(\tau_r)$ was defined by (\ref{9.005}).
Also, $x_j$ (for $j=1, \dots, 2g$) are our
 chosen set of generators of $H_1 (\Sigma)$.

Theorem \ref{mainab} applies when 
$\eta(X) = \tau(X) \exp 
\sum_{j=1}^{2g} \sum_{r \ge 2}  s_r^j \tbrj $  for 
$s_r^j \in \CC$ and $\tau \in S(\lieks)^K$.  Define 
$S^j \in S(\lieks)^K$ by 
$S^j (X) = \sum_{r \ge 2} s_r^j \tau_r(X)$; we then define
$\tilde{b}_{S^j}^j$ by 
$\tilde{b}_{S^j}^j (X)  = 
\sum_{r \ge 2} s_r^j \tbrj$. 
\begin{lemma} \label{l9.9'}  

(a) The restriction to $T^{2g}$ of 
$\tbrnox(-X) $ is 
$    \sum_{a = 1}^{n-1}  (d \tau_r)_X (\hu{a} ) \gen_a^j $ where 
$\gen_a^j$ (for $a = 1, \dots, {n-1}$  and 
$j = 1, \dots, 2g$) are the elements of the basis of 
$H^1 (T^{2g}) $ corresponding to  an oriented orthonormal basis $\{ \hu{a}
\} $ for $ \liet$. 

(b)  The restriction to $T^{2g}$ of 
$\tilde{b}_{S^j}^j(-X)   $ is 
$ \sum_{a = 1}^{n-1}  (d \Ssign{j})_X (\hu{a} ) \gen_a^j $.
\end{lemma}
\Proof We have by (8.21) of \cite{J2} that 
$$\tbrj_1 = \ev{x_j}^* \Phi_1^K (\tau_r). $$
By (\ref{9.7}) we have that 
$$ \Phi_1^K (\tau_r)|_{T^{2g}}(-X)  = 
\sigma_1^*  \bar{\Phi}_1^K (\tau_r)|_{T^{2g}}(-X)  =  
    \sum_{a = 1}^{n-1} (d\tau_r)_X (\hu{a}) \theta_a    $$
so 
$$ \tbrnox(-X)_1 \mid_{T^{2g}}  = 
\ev{x_j}^*  \sigma_1^*  \bar{\Phi}_1^K (\tau_r)|_{T^{2g}}(-X)  = 
    \sum_{a = 1}^{n-1} (d\tau_r)_X (\hu{a}) \gen_a^j    $$
since the generators $\gen_a^j$ of $H^1(T^{2g}) $ become identified 
with the components $\theta_a$ of the Maurer-Cartan form on the $j$-th
copy of $T$ in $T^{2g}$.\hfill $\square$

\begin{lemma} \label{l9.10} In the notation introduced just before
Lemma \ref{l9.9'}, we have 
\beq \label{eq10.13} \int_{T^{2g} }   \exp  \tfq (-X) 
 \exp \sum_{j,r} s_r^j  \tbrnox(-X)    = 
\int_{T^{2g} }    
\exp  \tfq(-X) 
 \exp \hattau(X) \eeq 
where $$\hattau(X) = - 
\sum_{a,b=1}^{n-1} \sum_{j=1}^{g} (d\Ssign{j})_X(\hu{a}) 
 (d\Ssign{j+g})_X(\hu{b})  (\partial^2 \qsgn)^{-1}_{ab} $$
\beq \label{hattaudef} = - \sum_{a,b=1}^{n-1}
\sum_{r,s  = 2}^{n}  \sum_{j=1}^{g}   
   s^j_r s^{j+g}_s (d\tau_r)_X(\hu{a}) 
 (d\tau_s)_X(\hu{b})  (\partial^2 \qsgn )^{-1}_{ab}. \eeq
Here, $\{ \hu{a}: a = 1, \dots, n-1 \} $ denotes an oriented orthonormal 
basis of $\liet$.
\end{lemma}
\Proof We need to consider 
the left hand side of (\ref{eq10.13}), which is 
\beq \label{9.11} 
\int_{T^{2g} } 
\exp  \tfq(-X)
 \exp \sum_{j,r} s_r^j  \tbrnox(-X) . \eeq
By Lemma \ref{l9.7a} the restriction of $\tfq(-X)$ to 
$T^{2g}$ is 
$\exp -\half \sum_{a,b} \sum_{j=1}^g(\partial^2 \qsgn )_X(\hu{a}, \hu{b})  
(\gen_a^j \gen_b^{j+g} -
\gen_a^{j+g} \gen_b^{j} ) $
while $
 \exp \sum_{j,r} s_r^j  \tbrnox(-X)  $ restricts on $T^{2g}$ (by 
Lemma \ref{l9.9'}) to 
$$\exp  \sum_{r,j}    s_r^j 
\sum_{a = 1}^{n-1} (d \tau_r)_X(\hu{a}) \z_a^j. $$
Thus for any given $j= 1, \dots, g$ we must compute the integral 
\beq \label{9.14} 
\int_{T^2} \exp
\Bigl (  \sum_{\sigma, \tau} y_\sigma A^{\sigma \tau} y_\tau /2 
+ \sum_\sigma y_\sigma B_j^\sigma
\Bigr ) \eeq
where $\sigma$ runs over pairs $(a,i)$ for 
$a = 1, \dots, n-1$ and $i = 0, 1$ (where
$i = 0$ corresponds to $j$ and $i = 1$ to $j+g$) and 
$y_{a,i}  = \z_a^{j+gi} $. Here, 
  the matrix $A$ is given by 
\beq 
A^{a 0, b1} = -\pqabx  = - A^{a 1, b0}; \; \; 
A^{a 0, b0}   =  A^{a 1, b1} = 0; \eeq
thus the  Pfaffian of $A$  (whose square is $\det A$) is given by 
$${\rm Pf}  (A) = \det \Bigl  (- \partial^2 q|_\liet\Bigr). $$
For $j = 1, \dots, g$ the vector $B_j$ is 
\beq \label{9.16} 
B_j^{a0} = -(\partial \Ssign{j})_X (\hu{a})\, ;  \onebl \onebl \onebl 
B_j^{a1} = -(\partial \Ssign{j+g})_X (\hu{a}). \eeq
The result is that 
\beq \label{9.17}
\int_{T^2} \exp
\Bigl (  \sum_{\sigma, \tau} y_\sigma A^{\sigma \tau} y_\tau /2 
+ \sum_\sigma y_\sigma B_j^\sigma
\Bigr )  = 
\int_{T^2}  \exp 
\Bigl (  \sum_{\sigma, \tau} y_\sigma A^{\sigma \tau} y_\tau /2 \Bigr ) 
\exp ( -  B_j^t A^{-1} B_j )/2  \eeq
$$ = {\rm Pf} (A) \exp ( -  B_j^t A^{-1} B_j )/2 , $$
where  $B_j^t$ denotes the transpose of the vector 
$B_j$.
Thus we find that (\ref{9.11}) becomes
\beq \label{9.12}
\det (-\partial^2 \qsgn |_\liet  )^g  \exp
\Biggl (   \sum_{a,b = 1}^{n-1}
\sum_{j = 1}^g   (d\Ssign{j} )_X (\hu{a}) 
(d\Ssign{j+g} )_X (\hu{b}) ( \partial^2  \qsgn )^{-1}_{a,b}
\Biggr ),  \eeq
which equals the right hand side of (\ref{eq10.13}).
This completes the proof of Lemma \ref{l9.10}.\hfill$\square$

Proposition \ref{p9.2}
 follows from Theorem \ref{t9.5} once we have shown that 
$$ \int_{T^{2g}} \exp \tfq(-X)  \exp \sum_{r \ge 2} s_r^j \tbrnox(-X)  = 
\int_{T^{2g} } \exp \tfq(-X)  \exp \hattau(X)$$ where $\hattau$ is given by 
(\ref{9.32a}). This is now clear from 
Lemma \ref{l9.10}. \hfill $\square$

\renorm 
\section{The Verlinde formula}

\nc{\tilc}[1]{ { \tilde{c}_{#1} } }
\nc{\wtilc}[1]{ { (w\tilde{c})_{#1} } }
\nc{\bracewtilc}[1]{ { \tildarg{w\tilde{c}}_{#1} } }

\nc{\dgk}{{ D_{n,d}(g,k) }} 
\nc{\vgk}{{ V_{n,d}(g,k) }} 
\nc{\lineb}{L}
\nc{\resid}[1]{\res_{#1 = 0 } }
\nc{\residone}[1]{\res_{#1 = 1 } }
\nc{\kmod}{r}
\nc{\gmax}{{\g_{\rm max} }}
\nc{\constq}{ {\frac{1}{n}  }}
\nc{\sol}{ { S_{0 \l} }}
\nc{\ch}{ { \rm ch} } 
\nc{\td}{ { \rm td} } 

\nc{\tg}{{\tilde{\gamma} } }
\renewcommand{\ell}{{\call}}

The Verlinde formula is a formula for the dimension $\dgk$ 
of the space of holomorphic sections of 
powers of $\ell$ , where $\ell$ 
is a particular line bundle over $\mnd$: it has been proved
by Beauville and Laszlo \cite{BL}, Faltings \cite{F}, Kumar,
Narasimhan and Ramanathan \cite{KNR} and Tsuchiya, Ueno
and Yamada \cite{TUY}. In this section 
we show how the Verlinde formula follows from our formula
(Theorem \ref{mainab}) for intersection pairings in $\mnd$.

A line bundle $\ell$ over $\mnd$ may be defined for which
$c_1 (\ell) = n f_2$, since $n f_2 \in H^2(\mnd, \ZZ)$
(see \cite{DN}).  As described in Section 1, this 
bundle is the {\em determinant line bundle}. 
 Whenever $k$ is a 
positive integer divisible by $n$, we then define
\beq \label{10.1} \dgk = \dim H^0(\mnd, \ell^{k/n} ). \eeq

Let us introduce
 $$\kmod = k + n;$$
let us also introduce 
the highest root $\gmax$, which is given by 
$\gmax(X) = X_n - X_1$ or 
$\gmax = e_1 + e_2 + \dots + e_{n-1}$. We then make the following
definition:
\begin{definition} \label{d10.1} The {\em Verlinde function}
$\vgk$ is given by 
$$\vgk = \sum_{\l \in \weightl_{\rm reg} \cap \liet_+: 
\inpr{\l, \gmax} < \kmod}
\frac{e^{-\itwopi \inpr{ \l - \rho ,\tc}  } } {(\sol(k) )^{2g-2} } $$
where $\rho$ is half the sum of the positive roots and 
$$ \sol(k) = \frac{1 }{\sqrt{n} \kmod^{(n-1)/2} }
\prod_{\g > 0 } 2 \sin \pi \inpr{\g,\l} /\kmod . $$
\end{definition}
(See \cite{GW} (A.44) and \cite{qym} (3.16).)
Verlinde's conjecture 
 says that the Verlinde function 
specifies the dimension of the space of holomorphic sections
of $\ell^{k/n}$:
\begin{theorem}  \label{verl} {\bf (Verlinde's conjecture)} 
$$ \dgk = \vgk.$$  \end{theorem}
We shall show how to extract  Verlinde's conjecture from our 
previous results: an outline  of the  method we use was given by 
Szenes \cite{Sz} (Section 4.2).

In fact  $H^i (\mnd, \ell^m ) = 0 $ for all $i > 0$ and
$m> 0$ by an argument using the Kodaira vanishing theorem and the
facts that
 $\call$  is a positive line 
bundle and 
the canonical bundle of $\mnd$ is equal to 
$\ell^{-2} $ 
(see \cite{Beauv1} Section 5 and Th\'eor\`eme
F of \cite{DN}),   so 
$\dgk$ is given for  $k>0 $ by the Riemann-Roch formula:

\beq \label{10.2} \dgk = \int_{\mnd }\ch \ell^{k/n} \td \mnd. \eeq
We use the following results to convert (\ref{10.2}) into a form
to which we may apply our previous results.
\begin{lemma} \label{p10.1} 
For any complex manifold 
$M$ the Todd class of $M$ is given by 
$$\td(M) = e^{c_1(M)/2} \hat{A} (M) $$ 
where $c_1(M)$ is the first Chern class of the holomorphic tangent 
bundle of $M$, and $\hat{A}(M)$ is the $A$-roof genus of $M$.
\end{lemma}
\Proof See for example \cite{Gilkey}, pages 97-99. \hfill $\square$

\begin{prop} \label{p10.2}
We have 
$$\hat{A} (\mnd) = \abk  \Bigl ( \prod_{\g > 0} 
\frac{ \g(X)/2 }{\sinh \g(X)/2 } \Bigr )^{2g-2}. $$
\end{prop}
\Proof This is proved by Newstead\footnote{Newstead writes the 
details of the proof 
only for $n=2$ but the same proof  yields the result  for general $n$.}  in 
 \cite{Newstead}.
\begin{lemma} \label{p10.3} 
We have 
$$c_1(\mnd) = 2 n f_2. $$
\end{lemma}
\Proof This is proved in  \cite{DN}, 
Th\'eor\`eme F  .\hfill $\square$

Of course the Chern character of $\ell^{k/n}$ is given by 
$\ch \ell^{k/n} = e^{k f_2}. $  
Thus we obtain
\begin{corollary} \label{p10.4}
The quantity $\dgk$  is given by 
$$ \dgk = \int_{\mnd} e^{(k + n)f_2} \abk 
  \Bigl ( \prod_{\g > 0} 
\frac{ \g(X) }{e^{\g(X)/2}  - e^{-\g(X)/2} } \Bigr )^{2g-2}. $$
\end{corollary} 
\Proof This follows immediately from (\ref{10.2}), Lemmas 
\ref{p10.1} and \ref{p10.3} and Proposition \ref{p10.2}.\hfill $\square$

\begin{theorem} \label{t10.5} We have
$$\dgk = \frac{(-1)^{n_+(g-1)}}{n!}\sum_{w \in W_{n-1} } 
 \resid{Y_{1}}
\dots \resid{Y_{n-1} }  \Biggl ( e^{\kmod \inpr{\tildarg{w\tc},X} } 
\int_{T^{2g}} e^{\kmod \omega} \times $$
\beq \prod_{\g > 0} 
\Bigl (\frac{ \g(X)}{e^{\g(X)/2} - e^{-\g(X)/2} }\Bigr )^{2g-2}  
\frac{1 } { \prod_{j = 1}^l (e^{ \kmod Y_j} - 1)
 \nusym(X)^{2g-2} } \Biggr ).\eeq
\end{theorem}
\Proof This comes straight from Corollary 
\ref{p10.4} and Theorem \ref{mainab}.
Note that because the factor $e^{f_2}$ in the statement of 
Theorem \ref{mainab} has been replaced by $e^{\kmod f_2} $, it is necessary
to replace 
$  e^{\inpr{\tildarg{w\tc},X} } $ by $e^{\kmod \inpr{\tildarg{w\tc},X}
} $, 
and 
$e^{Y_j} - 1$ by $e^{\kmod Y_j} - 1$
(cf. Remark \ref{rem8.3} (c) ).  \hfill   $\square$

We introduce $Z_j = \exp Y_j$. Since  for any $w \in W_{n-1}$ we have that
$$\tildarg{w\tc} = \bracewtilc{1} \he{1} + 
\bracewtilc{2} \he{2} + \dots + \bracewtilc{n-1} 
\he{n-1}  $$
(as in the statement of Proposition \ref{p:sz})
with $n \bracewtilc{j} \in \ZZ $ 
for all $j$, and $0 \le \bracewtilc{j} < 1$ for 
all  $j$, 
we obtain 
$$ e^{  \kmod \inpr{\bracewtilc,X} } = 
Z_1^{ \bracewtilc{1} \kmod } Z_2^{  \bracewtilc{2} \kmod } 
\dots Z_{n-1}^{ \bracewtilc{n-1} \kmod}. $$
(Recall that $k$ and $r$ are divisible by $n$ so 
$e^{  \kmod \inpr{\tc,X} }  $ is a well defined single valued
function of $Z_1$, $\dots, Z_{n-1}$.)
Thus we can equate $\dgk$ with 
$$ 
 \frac{(-1)^{n_+(g-1) } }{n!}
\sum_{w \in W_{n-1} } \residone{Z_{1}} \dots \residone{Z_{n-1} }   
\Biggl ( \Bigl ( \prod_{j = 1}^{n-1} 
\frac{1}{Z_j}   \Bigr ) 
\int_{T^{2g}} e^{\kmod \omega} \times 
$$ 
$$
\frac{ Z_1^{ \bracewtilc{1} \kmod } Z_2^{  \bracewtilc{2} \kmod} \dots 
Z_{n-1}^{ \bracewtilc{n-1}\kmod} 
}{\prod_{\g > 0 } (\tg^{1/2}- \tg^{-1/2})^{2g-2} (Z_1^{  \kmod} - 1) 
\dots (Z_{n-1}^{ \kmod} - 1 ) } \Biggr )   $$ 
\beq \label{10.5} 
=  \frac{(-1)^{n-1 + n_+(g-1) } }{n!} \sum_{w \in W_{n-1} } 
\residone{Z_{1}} \dots \residone{Z_{n-1} }   
\Biggl ( \Bigl ( \prod_{j = 1}^{n-1} 
\frac{1}{Z_j}  \Bigr ) \times \eeq
$$ \int_{T^{2g}} e^{\kmod \omega} 
\frac{ Z_1^{ - \bracewtilc{1} \kmod}
 Z_2^{ - \bracewtilc{2}\kmod} \dots Z_{n-1}^{ -\bracewtilc{n-1} \kmod} 
}{\prod_{\g > 0 } (\tg^{1/2}- \tg^{-1/2})^{2g-2} (Z_1^{ - \kmod} - 1) 
\dots (Z_{n-1}^{ - \kmod} - 1 ) } \Biggr ). $$
Here, we have introduced $\tg  $ defined (for the root
$\g  = e_r + e_{r+ 1} + \dots + e_{s-1}$) by 
$$\tg(Z_1, \dots, Z_{n-1}) = Z_r \dots Z_{s-1}.$$
We also have 
\begin{lemma} \label{p10.6} 
$$\int_{T^{2}} e^{\omega}  =  n $$
and hence 
$$\int_{T^{2g}} e^{r \omega}  =  r^{(n-1)g} n^g. $$
\end{lemma}
\Proof This follows from 
Lemma \ref{l9.03}.\hfill $\square$

The following may be proved by the same method as in Section 2
(see \cite{Sz}):
\begin{prop} \label{p:vsz}  Suppose $f$ is the 
meromorphic  function on the complexification
$T^{\CC}$ of $T$ defined by 
\beq \label{11.6} 
f(Z) = (-1)^{n-1} (-1)^{n_+(g-1)} \kmod^{(n-1)(g-1)} n^{g -1}
 \frac{Z_1^{-\tilc{1}\kmod} \dots Z_{n-1}^{-\tilc{n-1}\kmod} } 
{ \prod_{\g> 0 } (\tg^{1/2} - \tg^{-1/2} )^{2g-2} }.  \eeq
   Then we have  that
\beq\frac{1}{(n-1)!} \residone{Z_{1}} 
\dots \residone{Z_{n-1} } \sum_{w \in W_{n-1} } 
\prod_{j = 1}^{n-1} 
\Bigl ( \frac{\kmod }{Z_j} \Bigr )   \frac{\bracearg{ w f} (Z)} {\prod_{j = 1}^{n-1} 
(Z_j^{ - \kmod} - 1 ) } 
= \sum_{\l \in  \weightl_{\rm reg} \cap \lietpl : \inpr{\l,\gmax} < r} 
f(\exp \itwopi \l/\kmod). \eeq
Here, $W_{n-1} $ is the permutation group on 
$\{ 1, \dots, n-1 \} $ which is (isomorphic to) the Weyl group of
$SU(n-1)$, and $\bracearg{w f} $ is the function
\beq \label{11.6'} 
\bracearg{wf}(Z) = (-1)^{n-1} (-1)^{n_+(g-1)} \kmod^{(n-1)(g-1)} n^{g -1}
 \frac{Z_1^{-\bracewtilc{1}\kmod} \dots Z_{n-1}^{-\bracewtilc{n-1}\kmod} } 
{ \prod_{\g> 0 } (\tg^{1/2} - \tg^{-1/2} )^{2g-2} }.  \eeq

\end{prop}
\noindent{\em Remark:}
Notice that 
we have 
$$ \sum_{ \lambda \in \weightl_{\rm reg} \cap \lietpl : \inpr{\l,\gmax} < r
  }  f(\exp \itwopi \lambda/r) = \frac{1}{n-1} 
\sum_{m_j = 1}^{r-1}  f\Bigl (e^{
\itwopi ( \sum_j m_j w_j )/r } \Bigr ) . $$ (Here,
the $w_j$ are the fundamental weights, which
are dual to the simple roots.) The set 
$\{ X \in  \liet: \; $ $ X = \sum_j \lambda_j \he{j}, 0 \le \lambda_j
< 1, \; j = 1, \dots, 
n-1 \} $
 is a 
fundamental domain for the action of the integer lattice $\Lambda^I$ 
on $\liet$, while the set 
$\{ X \in  \lietpl \subset \liet: \gamma_{\rm max} (X) < 1 \} $
 is a fundamental domain for the {\em affine Weyl group}
$W_{\rm aff}$ 
(the semidirect product of the Weyl group and the integer lattice),
and  $\intlat$
 has index $(n-1)! $ (rather 
than $n!$) in $W_{\rm aff} $
(in other words a fundamental domain 
for $\intlat$ contains $(n-1)!$ fundamental
domains for $W_{\rm aff}$). This difference
accounts
for the  factor $1/(n-1)!$ in Proposition 
\ref{p:vsz} which replaces the factor $1/n!$ in  its analogue Proposition 
\ref{p:sz}.

Applying  Proposition  \ref{p:vsz} 
we find (recalling from Section 2 
that  $(-1)^{n-1} = c^\rho$ when $n$ and $d$ are coprime) that
\beq \label{10.10} 
\dgk =  (-1)^{n_+(g-1)} \kmod^{(n-1)(g-1)} n^{g -1} c^\rho 
\sum_{\l \in \weightl_{\rm reg} \cap \lietpl: 
\inpr{\l, \gmax} < \kmod } 
\frac{ e^{ - \itwopi \inpr{\tc, \l}  } }
 { \prod_{\g> 0 } (e^{\itwopi \inpr{\frac{\g}{2\kmod}, \l} } - 
e^{-\itwopi \inpr{\frac{\g}{2\kmod}, \l} } )^{2g-2}}. \eeq
This gives 
\beq \label{10.11}
\dgk = (-1)^{n_+(g-1)}  r^{(n-1)(g-1)} n^{g-1} 
\sum_{\l \in \weightl_{\rm reg} \cap \lietpl: 
\inpr{\l, \gmax} < \kmod } 
\frac{ e^{ - \itwopi \inpr{\tc, \l - \rho } }  }
{\prod_{\g > 0 } (2 i \sin \pi \inpr{\g,\l}/\kmod )^{2g-2} } \eeq
\beq \label{10.12}
 ~~~= \kmod^{(n-1)(g-1)} 
n^{g-1} 
\sum_{\l \in \weightl_{\rm reg} \cap \liet_+: 
\inpr{\l, \gmax} < \kmod } 
\frac{ e^{ - \itwopi \inpr{\tc, \l - \rho } }  }
{\prod_{\g > 0 }\Bigl (2  \sin 
\frac{\pi \inpr{\g,\l} }{\kmod }\Bigr )^{2g-2} }. \eeq
Comparing with Definition \ref{d10.1}, we see that $\dgk = \vgk$. 
This completes the proof of Theorem \ref{verl}. \hfill $\square$

\end{document}